\documentclass[fleqn,usenatbib]{mnras}

\usepackage{newtxtext,newtxmath} 
\usepackage{microtype} 

\usepackage[T1]{fontenc}
\usepackage{ae,aecompl}
\usepackage{pdflscape}

\usepackage{graphicx}	
\usepackage{amsmath}	
\usepackage{amssymb}	

\newcommand{\mbf}[1]{\mathbf{#1}}			%
\newcommand{\x}{\mbf{x}}

\newcommand{\Q}{\mathbf{Q}}
\newcommand{\D}{\mathcal{D}}
\newcommand{\A}{\mathbf{A}}
\newcommand{\0}{\mathbf{0}}
\renewcommand{\u}{\mathbf{u}}

\newcommand{\q}{\mathbf{q}}
\newcommand{\F}{\mathbf{F}}

\newcommand{\Bn}{\mathcal{B}}
\newcommand{\Vn}{\mathcal{V}}
\newcommand{\f}{\mathbf{f}}
\newcommand{\g}{\mathbf{g}}

\newcommand{\halb}{\frac{1}{2}} 

\newcommand{\B}{\mathbf{B}}
\newcommand{\de}[2]{\frac {\partial #1}{\partial#2}}

\renewcommand{\S}{\mathbf{S}}
\newcommand{\be}{\begin{equation} \begin{aligned} }
\newcommand{\ee}{\end{aligned} \end{equation}}

\newcommand{\uphi}{u_\phi}
\newcommand{\sphi}{\sin \phi \,}
\newcommand{\cphi}{\cos \phi \,}

\newcommand{\apriori}{\textit{a priori}}

\renewcommand{\epsilon}{\varepsilon}
\renewcommand{\phi}{\varphi}

\newtheorem{remark}{Remark}[section]

\usepackage{soul}

\title[Well balanced schemes for the Euler equations with gravity]{Well balanced Arbitrary-Lagrangian-Eulerian finite volume schemes on moving nonconforming meshes for the Euler equations of gasdynamics with gravity}

\author[E. Gaburro et al.]{
	Elena Gaburro,$^{1}$
	Manuel J. Castro,$^{2}$
	Michael Dumbser,$^{1}$\thanks{E-mail: michael.dumbser@unitn.it}
	\\
	$^{1}$Department of Civil, Environmental and Mechanical Engineering, University of Trento, Via Mesiano, 77 - 38123 Trento, Italy.\\
	$^{2}$Department of Mathematical Analysis, Statistics and Applied Mathematics, University of M\'alaga, Campus de Teatinos, 29071 M\'alaga, Spain. \\	
}

\pubyear{2017}

\begin{document}
	\label{firstpage}
	\pagerange{\pageref{firstpage}--\pageref{lastpage}}
	\maketitle
	
\begin{abstract} 
	In this work we present a novel second order accurate well balanced Arbitrary-Lagrangian-Eulerian (ALE) finite volume scheme on moving nonconforming meshes for the Euler equations of compressible gasdynamics with gravity in cylindrical coordinates. The main feature of the proposed algorithm is the capability of preserving many of the physical properties of the system exactly also on the discrete level: besides being conservative for mass, momentum and total energy, also \textit{any} known steady equilibrium between pressure gradient, centrifugal force and gravity force can be exactly maintained up to machine precision. Perturbations around such equilibrium solutions are resolved with high accuracy and with minimal dissipation on moving contact discontinuities even for very long computational times. 
	This is achieved by the novel combination of well balanced path-conservative finite volume schemes, that are expressly designed to deal with source terms written via nonconservative products, with ALE schemes on moving grids, which exhibit only very little numerical dissipation on moving contact waves. 
	In particular, we have formulated a new HLL-type and a novel Osher-type flux that are both able to guarantee the well balancing in a gas cloud rotating around a central object. Moreover, to maintain a high level of quality of the moving mesh, we have adopted a \textit{nonconforming} treatment of the sliding interfaces that appear due to the differential rotation. 
	A large set of numerical tests has been carried out in order to check the accuracy of the method close and far away from the equilibrium, both, in one and two space dimensions. 
\end{abstract}




\section{Introduction}

The main goal of this article is to develop a new family of numerical methods that allow to study problems in computational astrophysics connected with the rotation of gas clouds 
around a central object for very long computational times and with high accuracy. 

The physical situation we want to study is described by the Euler equations of compressible gas dynamics with an externally given gravitational field generated by a central object. A very important family of stationary solutions of the governing equations 
is characterized by the \textit{equilibrium} between pressure gradient, centrifugal force and gravity force. We suppose these equilibrium solutions to be known and want to design numerical methods  
that are able to preserve a rather wide class of such equilibria \textit{exactly} also on the discrete level (i.e. up to machine precision), so that small physical perturbations around the equilibrium solution can be solved with high accuracy and are not hidden by spurious numerical oscillations. However, at the same time our new numerical schemes are able to deal with situations far from the 
equilibrium, hence they do not fall into the class of perturbation methods. 

To preserve the equilibria in a system of equations with source terms, following \cite{Pares2006,castro2007well,MuellerWB}, we decide  to rewrite some of them in terms of non-conservative products  obtaining a system that can be cast in the following general form
\be
\label{eq.generalform}
\de{\Q}{t} + \nabla \cdot \F(\Q)  + \B(\Q) \cdot \nabla \Q = \mathbf{S}(\Q), \quad \x \in \Omega(t) \subset \mathbb{R}^2.
\ee
In this system, $\x$ is the spatial position vector, $t$ represents the time, $\Omega(t)$ is the computational domain at time $t$, $\Q = (q_1,q_2, \dots, q_{\nu})$ is the vector of the conserved variables defined in the space of the admissible states $\Omega_{\Q} \subset \mathbb{R}^{\nu}$, $ \F(\Q) = (\,\f(\Q), \g(\Q)\,) $ is the non linear flux tensor, $\B(\Q) = (\, \B_1(\Q), \B_2(\Q) \,) $ is a matrix collecting the non-conservative terms, and $\mathbf{S}(\Q)$ represents a non linear algebraic source term. The system \eqref{eq.generalform} can also be written in the following quasi-linear form
\be 
\label{eq.quasilinear}
\de{\Q}{t}  + \A(\Q) \cdot \nabla \Q = \mathbf{S}(\Q), \quad \x \in \Omega(t) \subset \mathbb{R}^2, 
\ee 
with the system matrix $\A(\Q) = \partial \F / \partial \Q + \B(\Q)$. The system is \textit{hyperbolic} if for any normal direction $\mathbf{n} \neq 0$ the matrix $\A(\Q) \cdot \mathbf{n}$ has $\nu$ real 
eigenvalues and a full set of $\nu$ linearly independent eigenvectors for all $\Q \in \Omega_{\Q}$. 
PDE systems like \eqref{eq.generalform} include as particular cases systems of conservation laws ($\B = \0$, $\mbf{S} = \0$), systems of conservation laws with source terms or balance laws ($\B = \0$), and  even non-conservative hyperbolic systems ($ \B \ne \0$). They appear in many fluid flow models in different contexts: shallow water models, multiphase flow models, compressible gas dynamics, etc.  

The main difficulty of systems written in this form, both from the theoretical and the numerical points of view, comes from the presence of non-conservative
products that do not make sense in the standard framework of distributions when the solution $\Q$ develops discontinuities.
Another difficulty is related to the numerical computation of stationary solutions: standard methods that solve correctly systems of conservation
laws can fail in this case when approaching equilibria or when simulating phenomena close to equilibrium solutions.

From the theoretical point of view, in this paper we assume the definition of non-conservative products as Borel measures given in \cite*{DalMaso1995}. 
This definition, which depends on the choice of a family of paths in the phase space $\Omega_{\Q}$, allows one to give a rigorous definition of weak solutions of (\ref{eq.generalform}).

We consider here the discretization of system (\ref{eq.generalform}) by means of numerical schemes which are \textit{path-conservative} in the sense introduced in \cite{Pares2006}. 
The concept of a path-conservative method, which is also based on a prescribed family of paths, provides a generalization of conservative schemes introduced by Lax for systems of conservation laws.
Moreover, the idea of constructing numerical schemes that preserve some equilibria, which are called \textit{well balanced} schemes, has been studied by many authors. 
The design of numerical methods with good properties is a very active front of research: see, for instance, \cite{Bermudez1994, audusse2004fast, bouchut2004nonlinear, castro2001q, rebollo2003family, rebollo2004asymptotically, castro2007wellb, gosse2000well, gosse2001well, greenberg1996well, greenberg1997analysis,  leveque1998balancing, pares2004well, perthame2001kinetic, perthame2003convergence, tang2004gas, toro2001shock}.

In particular, in the context of the Euler equations with gravity, in which the pressure forces are balanced by the gravitational forces, there is a growing interest in the community to construct new  numerical schemes that are able to achieve this precise balancing exactly even at the discrete level. 
At this point it has to be emphasized that conventional numerical schemes are in general \textit{not} able to preserve such stationary solutions, especially on coarse meshes, although the source term is discretized in a consistent manner, but 
consistency alone is not enough to achieve good results on coarse grids. 
This leads to erroneous numerical solutions especially when trying to compute small perturbations around the steady states necessitating the need for very fine meshes.  
Many recent papers have been devoted to this topic, in particular we refer to \cite{BottaKlein, Kapelli2014, KM15_630, Klingenberg2015, schaal2015astrophysical, desveaux2014well, desveaux2016well,  bermudez2016numerical} and the references therein.

An additional problem is usually given by the numerical dissipation on moving contact discontinuities which, in another context, is typically addressed by employing either pure Lagrangian schemes 
\cite{Despres2005,Maire2007,Maire2009b,carre2009cell} or indirect Arbitrary-Lagrangian-Eulerian (ALE) methods, see \cite{ShashkovRemap3,ShashkovRemap4,MaireMM2,scovazzi1,bochev2013fast}. 

However, all Lagrangian schemes are generally affected by a common problem that is  the severe mesh distortion or the mesh tangling that happens in the presence of shear flows and that may even destroy the computation.  
Hence, to reach long computational times, all Lagrangian methods must be in general combined with an algorithm to (locally) rezone the mesh at least from time to time and to remap the solution from the old mesh to the new mesh in a conservative manner. Lagrangian remesh and remap ALE schemes are very popular and some recent work on that topic can be found in the references on indirect ALE schemes 
listed above. In contrast to indirect ALE schemes (purely Lagrangian phase, remesh and subsequent remap phase) there are the so-called direct ALE schemes, where the local rezoning is performed before the computation of the numerical fluxes, that is, changing directly the chosen mesh velocity of the ALE approach, see for example \cite{boscheri2013arbitrary,Lagrange3D} for recent work in that direction based on high order Lagrangian ADER-WENO schemes, as well as \cite{Springel} for a powerful ALE framework on moving polygonal and polyhedral meshes. 

Our ALE scheme is based directly on a space-time conservation formulation of the governing PDE system, hence it fits in the framework of direct ALE schemes. 
Moreover, in order to avoid the typical mesh distortion caused by the shear flows, the sliding element interfaces are automatically detected during the computation,
and nodes along such sliding edges are allowed to move in a \textit{nonconforming} way by the insertion and deletion of new nodes and new edges.  
This strategy allows to maintain the quality of the moving mesh even for long computational times.  The robustness and efficiency of this approach has been tested in \cite{gaburro2016direct}  for the case of sliding interfaces lying over straight lines. In particular, this method is interesting when the mesh slides along circumferences and cylindrical coordinates are used, which 
is the case  here. For further references on the treatment of slide lines in Lagrangian schemes, the interested reader is referred to 
\cite{Caramana2009,LoubereSL2013,Clair201356,Clair2014315,bertoluzza2016conservative,gaburro2016direct}. 

To the very best of our knowledge, this is the first time that \textit{well balanced} numerical schemes are coupled with Arbitrary-Lagrangian-Eulerian schemes on moving nonconforming grids for the Euler equations with gravity.

The rest of the paper is organized as follows. 
First, in Section \ref{sec.euler.grav} we derive, from the standard Euler equations with gravity written in Cartesian coordinates, the equations written in cylindrical coordinates  where new source terms appear. 
Then in Section \ref{sec.NumMethod1d} we describe the details of the well balanced method for the one dimensional case and in Section \ref{sec.Results1d} we present some 1D numerical results. 
Later, in Section \ref{sec.NumMethod2d} we extend the method to two space dimensions and to moving nonconforming meshes.  Section \ref{sec.Results2d} is devoted to check the efficiency of the method with some nontrivial 2D test problems in a rotating Keplerian gas disk with variable density. In particular, the numerical results show that the proposed method significantly reduces the numerical dissipation on moving contact discontinuities in comparison with a standard non-well balanced Eulerian method on a fixed grid.

\section{Euler equations with gravity}
\label{sec.euler.grav} 

The Euler equations with gravity in two space dimensions represent a strongly hyperbolic system  that can be cast in the form of a system of balance laws by taking in (\ref{eq.generalform})
\begin{equation}  
	\label{eq.EulerCartCons}
	\begin{aligned}
		& \Q =  \left( \begin{array}{c} \rho \\[3pt]  \rho u_x \\[3pt] \rho u_y \\[3pt] \rho E \end{array} \right),  \
		\f(\Q) = \left( \begin{array}{c} \rho u_x \\[3pt] \rho u_x^2 + P \\[3pt] \rho u_x u_y \\[3pt] u_x( \rho E + P )   \end{array} \right), \
		\g(\Q) = \left( \begin{array}{c} \rho u_y \\[3pt] \rho u_x u_y \\[3pt] \rho u_y^2 + P \\[3pt] u_y( \rho E + P)    \end{array} \right), \!\!\!\!\!\!\!\!\!\!\!\!\!\! \!\!\!\!\!\!\!\!\!\!\!\!\!\! \!\!\!\!\!\!\!\!\!\!\!\!\!\! \!\!\!\!\!\!\!\!\!\!\!\!\!\! \!\!\!\!\!\!\!\!\!\!\!\!\!\! \!\!\!\!\!\!\!\!\!\!\!\!\!\! \!\!\!\!\!\!\!\!\!\!\!\!\!\! \\
		& \B(\Q) = 0, \quad
		\S(\Q) = \left( \begin{array}{c} 0 \\[3pt] -\cphi \rho \, \frac{G \, m_s}{r^2} \\[3pt] -\sphi \rho \, \frac{G \, m_s}{r^2} \\[3pt] -  (u_x \cphi + u_y \sphi) \rho \frac{G  m_s}{r^2}  \end{array} \right).
	\end{aligned}
\end{equation}
Here $\rho$ is the density, $u_x$ and $u_y$ are respectively the velocities along the $x$ and $y$ directions, $r = \sqrt{x^2+y^2}$, $\phi = \arctan(y/x)$,  $E$ is the specific total energy 
(excluding the gravitational energy), $m_s$ is the mass of the central object, $G$ is the gravitational constant and the pressure $P$ is given by 
\be
P = (\gamma - 1) \left ( \rho E  - \frac{1}{2} \rho \,\left (u_x^2 + u_y^2 \right ) \right ), \quad \gamma = \frac{c_p}{c_v} > 1, 
\ee
where $\gamma$ is the ratio of the specific heats at constant pressure and at constant volume, and which is supposed to be constant. 

Now we are interested in studying \textit{rotational phenomena} affected by sheared vortex flows, so we decide to rewrite the Euler equations in \textit{cylindrical coordinates} $(r,\phi)$ according to the usual relations  
\be 
	x = r \cphi, \qquad y = r \sphi.  
\ee 

Let $u_r$ and $\uphi$ be respectively the radial and the angular component of the velocity, linked to $u$ and $v$ by 
\be
\label{eq.PolarVelocities}
	u_x =  \cphi u_r - \sphi  \uphi,    \qquad 
	u_y =  \sphi u_r + \cphi \uphi \ , 
\ee
and consider the map 
\be
\label{eq.changeOfDer}
	\de{}{x} =  \, \cphi  \de{}{\rho} - \frac{\sphi}{\rho} \de{}{\phi}, \qquad 
	\de{}{y} = \, \sphi  \de{}{\rho} + \frac{\cphi}{\rho} \de{}{\phi}\ .  
\ee
To shorten the notation, from now on, we denote the radial velocity $u_r$ by $u$, and the angular velocity $u_\phi$ by $v$. 

By substituting into \eqref{eq.EulerCartCons} the expressions given in \eqref{eq.PolarVelocities} and \eqref{eq.changeOfDer}, after some calculations, we derive a new set of hyperbolic equations
that still takes the form \eqref{eq.generalform} with 
\be 
\label{eq.EulerPolar}
\begin{aligned}
	& \Q \!=\! \left( \begin{array}{c} \!r \rho\!  \\[3pt] \!\!r \rho u\!  \\[3pt]  \!r \rho v \!  \\[3pt] \!r \rho E\! \\[3pt] r  \end{array} \right), \ 
	\f(\Q) \!=\!  \left( \begin{array}{c} \!r \rho u\! \\[3pt] \!\!r \rho u^2 + rP\!\! \\[3pt] \!r \rho u v\! \\[3pt] \!\!r u( \rho E + P )\!\! \\[3pt] 0 \end{array} \right),  \ 
	\g(\Q) \!=\!  \left( \begin{array}{c}  \!\rho v\! \\[3pt] \!\!\rho u v\!\! \\[3pt] \!\!\rho v^2 + P \!\!\\[3pt]  \!\!v( \rho E  + P)\!\! \\[3pt] 0  \end{array} \right), \!\!\!\!\!\!\!\!\!\!\!\!\!\! \!\!\!\!\!\!\!\!\!\!\!\!\!\! \!\!\!\!\!\!\!\!\!\!\!\!\!\! \!\!\!\!\!\!\!\!\!\!\!\!\!\! \!\!\!\!\!\!\!\!\!\!\!\!\!\! \!\!\!\!\!\!\!\!\!\!\!\!\!\! \!\!\!\!\!\!\!\!\!\!\!\!\!\! \!\!\!\!\!\!\!\!\!\!\!\!\!\! \\
	& \B(\Q)  = 0, \quad 
	\mbf{S}(\Q) = 
	\left( \begin{array}{c} 0 \\[3pt] -\, \rho \frac{G m_s}{r} + P + \rho v^2  \\[3pt] - \, \rho u v \\[3pt] -\, \rho u \frac{G m_s}{r} \\[3pt] 0 \end{array} \right) .
\end{aligned}
\ee
Note that the system is written in terms of conserved variables, which is made possible by the insertion of an additional trivial equation 
\be 
\label{eq.TrivialEq}
\de{r}{t} = 0,
\ee
which implies that the radius $r$ is both a coordinate and a conserved variable.

\medskip 
The \textit{goal} of our work is to construct a finite volume scheme that is second order accurate in general situations, and, at the same time, can solve exactly (i.e. up to machine precision) a class of stationary solutions given by 
\be
\label{eq.EquilibriaConstraint}
	\rho = \rho(r), \qquad 
	u = 0, \qquad 
	\de{v}{\phi} = 0. 
\ee 
Looking at the second equation in (\ref{eq.EulerPolar}) and the equilibrium constraints in (\ref{eq.EquilibriaConstraint}), we notice that equilibria should balance the pressure and gravitational forces. More precisely
\be
\label{eq.PressureAndGravForces}
\de{rP}{r} = - \rho \left (  \frac{Gm_s}{r} - v^2 \right) + P.
\ee

This relation has to be achieved also at the discrete level in order to preserve these stationary solutions.
In standard finite volume schemes, fluxes and sources are typically discretized in different ways and therefore, the balancing between them is usually lost.  In order to construct a numerical scheme that exactly preserves those stationary solutions, here we first rewrite the equations in the following way, where both, pressure and gravitational forces (\ref{eq.PressureAndGravForces}) are treated as non-conservative terms.Thus,  by exploiting some trivial equalities as
\be
\de{rP}{r} = P + r \de{P}{r} \quad \text{and} \quad \de{r}{r} = 1,
\ee
the forces in (\ref{eq.PressureAndGravForces}) can be rearranged as
\be
\label{eq.PressureAndGravForcesWellBalanced}
r \de{P}{r} + \left ( \rho \frac{G m_s}{r} - \rho v^2 \right ) \de{r}{r} = 0,
\ee
and finally the Euler equations with gravity in polar coordinates can be cast in form~(\ref{eq.generalform}) with non trivial non-conservative terms and with zero algebraic source term as 
\be 
\label{eq.EulerPolarWellBalanced}
\begin{aligned}
	& \Q = \left( \begin{array}{c} \!r \rho\!  \\[3pt] \!r \rho u\!  \\[3pt]  \!r \rho v\!  \\[3pt] \!r \rho E\! \\[3pt] r  \end{array} \right), \ 
	\f(\Q) = \left( \begin{array}{c} \!r \rho u\! \\[3pt] \!\!r \rho u^2 \!\! \\[3pt] \!r \rho u v\! \\[3pt] \!\!r u( \rho E + P )\!\! \\[3pt] 0 \end{array} \right),  \ 
	\g(\Q) =  \left( \begin{array}{c}  \!\rho v\! \\[3pt] \!\!\rho u v \!\! \\[3pt] \!\!\rho v^2 + P \!\!\\[3pt]  \!\!v( \rho E  + P)\!\! \\[3pt] 0  \end{array} \right), \\
	& \mbf{S}(\Q) \! =\! \0 , \ \qquad  \ 
	\B(\Q) \! \cdot \!  \nabla \Q \! = \! \left ( \begin{array}{c} 0 \\[3pt] \!\!\!\! r \de{P}{r} + \left ( \rho \frac{G m_s}{r} - \rho v^2  \right ) \de{r}{r} \!\!\!\! \\[3pt] \left( \rho u v \right ) \de{r}{r}	\\[3pt] \rho u \frac{G m_s}{r}\!  \de{r}{r}	\\[3pt] 0\\[3pt] \end{array} \right ),  
	\end{aligned}
\ee
i.e. 
\be
\begin{aligned}
	& \B_1  = \left ( \begin{array}{ccccc}
		0 & 0 & 0 & 0 & 0 \\[3pt] 
		\!\!\! r \de{P}{q_1} & \!r \de{P}{q_2} & \!r \de{P}{q_3} & \!r \de{P}{q_4} & \!r \de{P}{q_5} \!+\! \rho\frac{Gm_s}{r} \!-\! \rho v^2 \!\!\!  \\[3pt] 
		0 & 0 & 0 & 0 & \!\! \rho u v \\[3pt]
		0 & 0 & 0 & 0 & \!\! \rho u \frac{G m_s}{r}  \\[3pt]
		0 & 0 & 0 & 0 & 0 \\
	\end{array}  \right ),  \\
		& \B_2 = \0, 
\end{aligned} 
\ee
where $q_i$, $i=1, \cdots, 5$ denotes the $i$-th component of vector $\Q$.
Notice that it is possible to write the source terms as non-conservative products thanks to the introduction of the coordinate $r$ also as conserved variables (see the added equation in (\ref{eq.TrivialEq})), which is the typical strategy adopted in \cite{greenberg1996well, greenberg1997analysis,gosse2000well, gosse2001well, castro2007well}.

In the following, we first focus on the one dimensional version of the previous system (where $\g$ and $\B_2$ are not considered) achieving an exact balancing in the radial direction $r$. Then, we will extend the method to two space dimensions and moving nonconforming meshes.
In both cases the key point of our new numerical method is the discretization of the terms in (\ref{eq.PressureAndGravForcesWellBalanced}).	 

\section{Numerical method in one dimension}
\label{sec.NumMethod1d}

For the numerical approximation of the one dimensional system, the spatial domain is discretized by $N$ cells (or finite volumes) $I_i=[r_{i-1/2}, r_{i+1/2}]$ of regular size $\Delta r = r_{i+1/2}-r_{i-1/2}$, $ i=1,\dots,N$. 
After having integrated Eq.(\ref{eq.generalform}) over a cell $I_i$, we approximate the time derivative of the cell averages $\Q_i(t)$ at each time $t$ by  a path-conservative scheme: 
\be 
\label{eq.SemidiscreteScheme}
\begin{aligned} 
	\frac{d \Q_i}{d t}(t) = 
	& - \frac{\Delta t}{\Delta r} \!\left (\! \D_{i\!-\!\frac{1}{2}}^+\left( \q_{i-\frac{1}{2}}^{-}(t), \q_{i-\frac{1}{2}}^{+}(t)  \! \right) +  \D_{i\!+\!\frac{1}{2}}^-\left( 
	\q_{i+\frac{1}{2}}^{-}(t),  \q_{i+\frac{1}{2}}^{+}(t)\!  \right) \! \right  )  \!\!\!\!\!\!\!\!\!\!\!\!\!\! \!\!\!\!\!\!\!\!\!\!\!\!\!\! \!\!\!\!\!\!\!\!\!\!\!\!\!\! \!\!\!\!\!\!\!\!\!\!\!\!\!\! \!\!\!\!\!\!\!\!\!\!\!\!\!\! \!\!\!\!\!\!\!\!\!\!\!\!\!\! \!\!\!\!\!\!\!\!\!\!\!\!\!\!  \\
	& -  \frac{\Delta t}{\Delta r} \int \limits_{r_{i-\frac{1}{2}}}^{r_{i+\frac{1}{2}}} \de{}{r}\f \left( \q_i(r,t) \right)  dr  \\
	& -  \frac{\Delta t}{\Delta r} \int \limits_{r_{i-\frac{1}{2}}}^{r_{i+\frac{1}{2}}} \B_1\left( \q_i(r,t) \right)\de{}{r}\left( \q_i(r,t) \right) dr .
\end{aligned} 
\ee 
In the scheme, $\q_i(r,t)$ is the approximation of the conserved variables inside cell $I_i$ at time $t$, computed via a reconstruction operator from the conserved variables $\Q_i(t)$ in a given stencil,  while $ \q_{i-\frac{1}{2}}^{+}(t)=\q_i(r_{i-1/2},t)$ and  $ \q_{i+\frac{1}{2}}^{-}(t)=\q_i(r_{i+1/2},t)$ denote the evaluation of $\q_i(r,t)$ at the left and right boundaries of cell $I_i$. According to  \citet{Pares2006} and   \citet{Castro2012} $\D^{\pm}_{i+\frac{1}{2}}$ can be defined as follows: 
\be 
\label{eq.FIrstOrderDgeneric}
\D^{\pm}_{i+\frac{1}{2}} \left( \q_{i+\frac{1}{2}}^{-}, \q_{i+\frac{1}{2}}^{+} \right) & = \frac{1}{2} \biggl ( \ \f(\q_{i+\frac{1}{2}}^{+}) - \f(\q_{i+\frac{1}{2}}^{-}) \ + \\
&  \Bn_{i+\frac{1}{2}} \left ( \q_{i+\frac{1}{2}}^{+}-\q_{i+\frac{1}{2}}^{-} \right ) \pm \Vn_{i+\frac{1}{2}} \left ( \q_{i+\frac{1}{2}}^{+}-\q_{i+\frac{1}{2}}^{-} \right ) \Biggr),
\ee 
where $\f(\q)$ is the physical flux, $\Bn_{i+\frac{1}{2}} \left (\q_{i+\frac{1}{2}}^{+}-\q_{i+\frac{1}{2}}^{-} \right )$ is the discretization of the non-conservative terms and  $\Vn_{i+\frac{1}{2}} \left (\q_{i+\frac{1}{2}}^{+}-\q_{i+\frac{1}{2}}^{-}\right )$ is the viscosity term, that characterizes the method. In \eqref{eq.FIrstOrderDgeneric}, the dependency on $t$ has been dropped for simplicity.

$\Bn_{i+\frac{1}{2}} \left (\q_{i+\frac{1}{2}}^{+}-\q_{i+\frac{1}{2}}^{-}\right )$ and $\Vn_{i+\frac{1}{2}} \left (\q_{i+\frac{1}{2}}^{+}-\q_{i+\frac{1}{2}}^{-}\right )$ are defined in terms of a family of paths $\Phi(s, \q_{i+\frac{1}{2}}^{-},\q_{i+\frac{1}{2}}^{+})$, $s \in [0,1]$.  

In general,  according to  the theory of \cite{DalMaso1995}, the family of paths should be a Lipschitz continuous family of functions  $\Phi(s,\Q_L, \Q_R), \ s \in [0,1]$  satisfying some regularity and compatibility conditions, in particular,
\be 
\label{eq.Pathproperties}
\begin{aligned}
	\Phi(0,\Q_L, \Q_R) = \Q_L, \quad \Phi(1; \Q_L, \Q_R) = \Q_R, \quad \Phi(s, \Q, \Q) = \Q.
\end{aligned}
\ee
Moreover, according to \citet{Pares2006}, $\D^{\pm}$  should satisfy the following properties:
\be 
\begin{aligned}
	\D^{\pm} (\Q, \Q) = \0 \quad \forall \Q \in \Omega_{\Q},
\end{aligned}
\ee
being $\Omega_{\Q}$ the set of admissible states for the problem, and, for every $\Q_L$, $\Q_R \ \in \Omega_{\Q}$, 
\be  
\label{PC}
\begin{aligned}
	\D^-\!(\Q_{\!L},\! \Q_{\!R}) \!+\! \D^+\!(\Q_{\!L},\! \Q_{\!R}) \!=\!\!\! \int_0^1 \!\!\!\!\!\! \mbf{A}(\Phi(s;\Q_{\!L},\!\Q_{\!R}))\de{\Phi}{s}(s;\Q_{\!L}, \!\Q_{\!R})ds,
\end{aligned}
\ee
where
\be 
\label{eq.ExtendedJacobian}
\mbf{A}(\Q) = \mbf{J_f}(\Q) + \B_1(\Q), \  \qquad \  \mbf{J_f}(\Q) = \de{\f(\Q)}{\Q}
\ee 
with $\mbf{J_f}$ denoting the Jacobian of the flux function $\f$.  
Note that, in this particular case equation \eqref{PC} could also be rewritten as follows:
\be  
\begin{aligned}
	\D^-\!(\Q_{\!L},\! \Q_{\!R}) \!+\! \D^+\!(\Q_{\!L},\! \Q_{\!R}) \!= \f(\Q_R)-\f(\Q_L) + \Bn_{LR}(\Q_R-\Q_L),
\end{aligned}
\ee
where
\be
\label{prop-B}
\Bn_{LR}(\Q_R-\Q_L)= \!\! \int_0^1 \!\!\! \B_1(\Phi(s;\Q_{\!L},\!\Q_{\!R}))\de{\Phi}{s}(s;\Q_{\!L}, \!\Q_{\!R})ds. 
\ee
The interested reader is referred to \cite{DalMaso1995} and \cite{Pares2006} for a rigorous and complete presentation of this theory.

In this paper, the family of paths will be chosen so that stationary solutions given by \eqref{eq.EquilibriaConstraint}-\eqref{eq.PressureAndGravForces} are preserved.

The rest of this section is organized as follows: we start by proposing two different first order well balanced schemes, the first one is denoted by Osher-Romberg scheme, and the second one is a well balanced HLL scheme. Next we propose a second order scheme constructed using the previous first order schemes in combination with a second order well balanced reconstruction operator.

\subsection{First order well balanced schemes}

Let us remark first, that the scheme \eqref{eq.SemidiscreteScheme} reduces to
\be 
\label{eq.SemidiscreteScheme_1}
\begin{aligned} 
	\frac{d \Q_i}{d t}(t) = 
	& - \frac{\Delta t}{\Delta r} \!\left (\! \D_{i\!-\!\frac{1}{2}}^+\left( \q_{i-\frac{1}{2}}^{-}(t), \q_{i-\frac{1}{2}}^{+}(t)  \! \right) +  \D_{i\!+\!\frac{1}{2}}^-\left( 
	\q_{i+\frac{1}{2}}^{-}(t),  \q_{i+\frac{1}{2}}^{+}(t)\!  \right) \! \right  ), \\
\end{aligned} 
\ee 
if $\q_i(r,t) = \Q_i(t)$ is constant within each cell, for every time $t$ and coincides with the cell average $\Q_i(t)$. The time derivative is discretized by the first order explicit Euler method. Thus, the resulting  scheme will be first order accurate in space and time. Moreover,  $\q_{i+\frac{1}{2}}^{-}=\q_i=\Q_{i}$ and $\q_{i+\frac{1}{2}}^{+}=\q_{i+1}=\Q_{i+1}$.

Therefore, to determine the numerical scheme,  $\Bn_{i+1/2}\left(\q_{i+1}-\q_{i}\right)$ and $\Vn_{i+\frac{1}{2}} \left (\q_{i+1}-\q_{i} \right )$ should be defined.  In order to define $\Bn_{i+1/2}\left(\q_{i+1}-\q_{i}\right)$, a family of paths should be prescribed, so that the resulting scheme is well balanced. Note that if the standard segment path is prescribed, that is 
\be
\Phi(s; \q_i, \q_{i+1}) = \q_i + s(\q_{i+1} - \q_i),
\ee
then, the resulting scheme is not well balanced for this set of stationary solutions.  Here we propose the following family of paths.  Let $\Phi^E(s,\Q^E_i, \Q^E_{i+1})$ be a reparametrization of a stationary solution given by \eqref{eq.EquilibriaConstraint}-\eqref{eq.PressureAndGravForces} that connects the state $\Q^E_i$ with $\Q^E_{i+1}$, where $\Q^E_i$ is the cell average of the given stationary solution in the cell $I_i$. Note that in the case of first and second order schemes $\Q^E_i$ could be approximated by the evaluation of the stationary solution at the center of the cell.  Then we define $\Phi(s,\q_i, \q_{i+1})$ as follows
\be
\label{eq:path}
\Phi(s,\q_i,\q_{i+1})=\Phi^E(s,\Q^E_i, \Q^E_{i+1})+\Phi^f(s,\q^f_i, \q^f_{i+1}), 
\ee
where $\q^f_i=\q_i- \Q^E_i$ and $\q^f_{i+1}=\q_{i+1}- \Q^E_{i+1}$ and
\be 
\Phi^f(s,\q^f_i, \q^f_{i+1})=  \q^f_i + s(\q^f_{i+1} - \q^f_i).
\ee
That is, $\Phi^f$ is a segment path on the {\it fluctuations} with respect to a given stationary solution. With this choice, it is clear that if $\q_i$ and $\q_{i+1}$ lie on the same stationary solution satisfying \eqref{eq.EquilibriaConstraint}-\eqref{eq.PressureAndGravForces}, then $\q^f_i=\q^f_{i+1}=\0$ and $\Phi$ reduces to $\Phi^E$. In such situations we have that $\f(\q_{i+1})=\f(\q_i)=\0$ and
\be
\Bn_{i+\frac{1}{2}}(\q_{i+1}-\q_i)=\!\!\!\int_0^1 \!\!\!\!\! \B_1(\Phi^E(s,\q_i,\q_{i+1})) \de{\Phi^E}{s}(s;\q_i, \q_{i+1})ds\! =\!\0.
\ee
Therefore
\be 
\f(\q_{i+1})-\f(\q_i)+\Bn_{i+1/2}(\q_{i+1}-\q_i)=\0. 
\ee
For the sake of simplicity, in the following we will use the notation $\Phi(s)$ instead of $\Phi(s; \q_i, \q_{i+1})$ when there is no confusion.

Let us now define $\Bn_{i+1/2}(\q_{i+1}-\q_i)$ in the general case, where $\q_{i+1}$ and $\q_i$ do not lie on a stationary solution. In this case we have that
\begin{equation}
	\label{eq-B}
	\Bn_{i+1/2}(\q_{i+1}-\q_i)= \left( b^{i+1/2}_1 \ b^{i+1/2}_2 \ b^{i+1/2}_3 \ b^{i+1/2}_4 \ b^{i+1/2}_5 \right)^T.
\end{equation}
It is clear from the definition of $\B_1$ that 
\begin{equation} \label{eq-b15}
	b^{i+1/2}_1=b^{i+1/2}_5=0,  
\end{equation}

\be
b^{i+1/2}_2 = \int_0^1 \Phi_r(s) \de{\Phi_P}{s}(s) + \Phi_{(r\rho)}(s) \Phi_{\zeta_r}(s) \de{\Phi_r}{s}(s) ds,
\ee
where $\Phi_r(s)=\Phi_r(s; r_i,r_{i+1})=r_i + s(r_{i+1}-r_i)$, $\Phi_P(s)=\Phi^E_P(s)+\Phi^f_P(s)$, $\Phi_{(r\rho)}(s)(s)=\Phi^E_{(r\rho)}(s)+\Phi^f_{(r\rho)}(s)$ and, finally, $\Phi_{\zeta_r}(s)=\Phi^E_{\zeta_r}(s)+\Phi^f_{\zeta_r}(s)$ where
\be
\label{eq.zeta}
\zeta_r(r) = \left ( \frac{G m_s}{r^2} - \frac{v^2}{r} \right ), \quad\text{with } \ \zeta(r) = \!\int \!\zeta_r(r) dr.
\ee

Taking into account that 
\be
\int_0^1 \Phi_r(s) \de{\Phi^E_P}{s}(s) + \Phi^E_{(r\rho)}(s) \Phi^E_{\zeta_r}(s) \de{\Phi_r}{s}(s) ds=0,
\ee
$b^{i+1/2}_2$ can be rewritten as follows:
\begin{eqnarray} 
b^{i+1/2}_2 &=&  \int_0^1 \!\! \Phi_r(s) \de{\Phi^f_P}{s}(s) ds \nonumber \\ 
            &&   + \int_0^1  \left(\!\Phi^E_{(r\rho)}\!(s) \Phi^f_{\zeta_r}\!(s) + \Phi^f_{(r\rho)}\!(s)\Phi_{\zeta_r}\!(s)\right) \!\!\de{\Phi_r}{s}(s) ds. \nonumber \\ 
\end{eqnarray} 
Note that,  $\de{\Phi^f_P}{s}(s)=P^f_{i+1}-P^f_i$ and $\de{\Phi_r}{s}(s)=r_{i+1}-r_i=\Delta r_{i+1/2}$. Observe that in uniform meshes $\Delta r_{i+1/2} =\Delta r$.  With the previous notation $b^{i+1/2}_2$ reduces to
\begin{eqnarray} 
b^{i+1/2}_2 & = & r_{i+1/2}\Delta P^f_{i+1/2}  \nonumber \\ 
            &&    + \left(\int_0^1\!\!\! \left( \Phi^E_{(r\rho)}\!(s) \Phi^f_{\zeta_r}\!(s)\! +\! \Phi^f_{(r\rho)}\!(s)\Phi_{\zeta_r}\!(s)\right) ds \right)\! \Delta r_{i+1/2}, \nonumber \\ 
\end{eqnarray} 
where $r_{i+1/2}=\frac{r_i+r_{i+1}}{2}$ and
$	\Delta P^f_{i+1/2}= P^f_{i+1}-P^f_i$.

In general, the integral term could be difficult to compute, therefore we propose to use a numerical quadrature formula. Here the mid-point rule is used. In this case, we define $b^{i+1/2}_2$ as follows:
\begin{eqnarray} 
\label{b2}
b^{i+1/2}_2 & = & \left( (r\rho)^E_{i+1/2} (\zeta_r)^f_{i+1/2} \!+\! (r\rho)^f_{i+1/2} (\zeta_r)_{i+1/2}  \right)\!\Delta r_{i+1/2} \nonumber \\ 
               && + r_{i+1/2}\Delta P^f_{i+1/2}, 
\end{eqnarray} 
where
\be 
(r\rho)^E_{i+1/2}= \Phi^E_{(r\rho)}(1/2), \ (\zeta_r)^f_{i+1/2} = \frac{(\zeta_r^f)_i + (\zeta_r^f)_{i+1}}{2},
\ee
\be 
(r\rho)^f_{i+1/2}= \frac{(r\rho)^f_i + (r\rho)^f_{i+1}}{2}, \text{ and } \ (\zeta_r)_{i+1/2}=\Phi_{\zeta_r}(\frac{1}{2}).
\ee
It is clear from the definition that $b^{i+1/2}_2=0$ if $\q_i$ and $\q_{i+1}$ lie on the same stationary solution as $\Delta P^f_{i+1/2}=0$, $(r\rho)^f_{i+1/2}=0$ and $(\zeta_r)^f_{i+1/2}=0$.

Finally, terms $b^{i+1/2}_3$ and $b^{i+1/2}_4$ could be approximated in the same way. Nevertheless, as those terms explicitly depend on $u$ and we are interested in preserving equilibria with $u=0$,  
a simpler approach can be used. Thus, $b^{i+1/2}_3$ is defined as 
\begin{equation}
	\label{eq-b3}
	b^{i+1/2}_3= \frac{(r \rho u)_{i+1/2}}{r_{i+1/2}} v_{i+1/2} \Delta r_{i+1/2},
\end{equation}
where 
\be
(r \rho u)_{i+1/2}= \frac{(r \rho u)_i +(r \rho u)_{i+1}}{2}, \ v_{i+1/2}= \frac{v_i + v_{i+1}}{2},
\ee
and $b^{i+1/2}_4$ as
\begin{equation}
	\label{eq-b4}
	b^{i+1/2}_4=(r \rho u)_{i+1/2} \frac{G m_s}{r_{i+1/2}^2} \Delta r_{i+1/2}.
\end{equation}
Note that both terms vanish when $u=0$.

As pointed out in \cite{Pares2006}, a sufficient condition for a first order path-conservative scheme to be well balanced is that $\D^{\pm}_{i+1/2} \left(\q_i, \q_{i+1} \right) =\0$,
if $\q_i$ and $\q_{i+1}$ lie on the same stationary solution. Therefore, with the previous choice of paths, $\D^{\pm}_{i+1/2}=0$ if $\Vn_{i+1/2}(\q_{i+1}-\q_i)=\0$. In the next paragraph we are going to present two different schemes defined in terms of two different viscosity terms, both of them verifying that $\Vn_{i+1/2}(\q_{i+1}-\q_i)=0$ for stationary solutions \eqref{eq.EquilibriaConstraint}-\eqref{eq.PressureAndGravForces}.

\subsubsection{Well-balanced Osher-Romberg scheme}
\label{ssec.OsherRomberg}
A path-conservative Osher-type scheme following \cite{OsherUniversal,OsherNC,ApproxOsher} can be cast in form (\ref{eq.FIrstOrderDgeneric}) with $\Vn(\q_{i+1}-\q_i)$ being defined as follows: 
\be
\label{eq.OsherMatrixFormal}
\Vn_{i+1/2}(\q_{i+1}-\q_i) = \int_0^1  \left | \mbf{A} \left ( \Phi(s)  \right)  \right | \partial_s \Phi(s) ds, \quad   0 \le s \le 1, 
\ee 
with $| \A | = \mathbf{R} |\boldsymbol{\Lambda}| \mathbf{R}^{-1}$ being the usual definition of the matrix absolute value operator given in terms of the right eigenvector
matrix $\mathbf{R}$, its inverse $\mathbf{R}^{-1}$ and the diagonal matrix of the absolute values of the eigenvalues $|\boldsymbol{\Lambda}| = \text{diag}(|\lambda_1|, |\lambda_2|, ..., |\lambda_\nu|)$. 
For the numerical approximation of the viscosity matrix, first we notice that it can be written as
\be
\Vn_{i+1/2}(\q_{i+1}-\q_i) = \int_0^1  \text{sign} \left( \mbf{A} \left ( \Phi(s)  \right) \right )  \mbf{A} \left ( \Phi(s)  \right)  \partial_s \Phi(s) ds,  
\ee 
with $\text{sign}( \A ) = \mathbf{R} \, \text{sign} ( \boldsymbol{\Lambda} ) \mathbf{R}^{-1}$ and $\text{sign} ( \boldsymbol{\Lambda} )$ the diagonal matrix containing the signs of all
eigenvalues of $\A$. Then, we approximate the previous expression by a quadrature formula as follows:
\be 
\Vn_{i+1/2}(\q_{i+1}-\q_i) =\sum_{j=1}^l \omega_j \text{sign} \left( \mbf{A}(\Phi(s_j)\right) \mbf{A}(\Phi(s_j))\partial_s \Phi(s_j).
\ee
Now, we propose to approximate $\mbf{A}(\Phi(s_j))\partial_s \Phi(s_j)$ by the following expression:
\be 
\mbf{A}(\Phi(s_j)) \partial_s \Phi(s_j) \approx \frac{\mbf{A}_{\Phi_j}}{2\epsilon_j}\left(\Phi(s_j+\epsilon_j)-\Phi(s_j-\epsilon_j)\right),
\ee
where $\mbf{A}_{\Phi_j}=A(\Phi(s_j-\epsilon_j),\Phi(s_j+\epsilon_j))$ is a Roe-matrix associated to the system (see  \cite{Pares2006} for details), that is a matrix satisfying
\be 
\begin{array}{l} 
\mbf{A}_{\Phi_j}\left(\Phi(s_j+\epsilon_j)-\Phi(s_j-\epsilon_j)\right)= \f(\Phi(s_j+\epsilon_j))-\f(\Phi(s_j-\epsilon_j)) \\ 
	\quad  +\Bn_{\Phi_j}\left(\Phi(s_j+\epsilon_j)-\Phi(s_j-\epsilon_j)\right),  
\end{array} 
\ee 
where $\Bn_{\Phi_j}\left(\Phi(s_j+\epsilon_j)-\Phi(s_j-\epsilon_j)\right)$ is defined as in the previous section using the states $\Phi(s_j-\epsilon)$ and $\Phi(s_j+\epsilon)$.
Therefore, the viscosity term reads as follows:
\be
\label{viscosity_n_osher}
\Vn_{i+1/2}(\q_{i+1}-\q_i) =\sum_{j=1}^l \omega_j \text{sign}\left( \mbf{A}(\Phi(s_j)\right) \frac{\mathcal{R}_j}{2\epsilon_j},
\ee
where 
\begin{eqnarray} 
\mathcal{R}_j &=& \f(\Phi(s_j+\epsilon_j))-\f(\Phi(s_j-\epsilon_j)) \nonumber \\
              & &  + \Bn_{\Phi_j}\left(\Phi(s_j+\epsilon_j)-\Phi(s_j-\epsilon_j)\right). 
\end{eqnarray} 
Note that if $\q_{i}$ and $\q_{i+1}$ lie on the same stationary solution we have $\Phi(s)=\Phi^E(s)$ and $\mathcal{R}_j=\0$, $j=1, \dots, l$ and $\Vn_{i+1/2}(\q_{i+1}-\q_{i})$ vanishes. Therefore, the numerical scheme \eqref{eq.SemidiscreteScheme_1} with \eqref{eq.FIrstOrderDgeneric}, where $\Bn_{i+1/2}(\q_{i+1}-\q_i)$ is defined as  \eqref{eq-B}, \eqref{eq-b15}, \eqref{b2}, \eqref{eq-b3} and \eqref{eq-b4}  and $\Vn_{i+1/2}(\q_{i+1}-\q_{i})$ is defined by \eqref{viscosity_n_osher} is exactly well balacend for stationary solutions given by \eqref{eq.EquilibriaConstraint}-\eqref{eq.PressureAndGravForces}.

Here we propose the Romberg method with $l=3$ and
\be 
\begin{array}{l}
	s_1=1/4, \ s_2=3/4, \ s_3=1/2, \\
	\omega_1=2/3, \ \omega_2=2/3, \ \omega_3=-1/3, \\
	\epsilon_1=1/4, \ \epsilon_2=1/4, \ \epsilon_3=1/2.
\end{array}
\ee
With this choice, the  viscosity term $\Vn_{i+1/2}(\q_{i+1}-\q_i)$ of the Osher-Romberg method reads as follows:
\be
\label{eq.ViscosityOR}
\begin{array}{l}
	\Vn_{i+\halb} \left( \q_{i+1}-\q_i \right) = \\
	\frac{4}{3} \text{sign}(\mbf{A}(\Phi(\frac{1}{4}))) \left( \f(\Phi(\frac{1}{2}))-\f(\q_i) +\Bn_{i+\frac{1}{4}}\left(\Phi(\frac{1}{2})-\q_i\right)\right) + \\
	\frac{4}{3}\text{sign}(\mbf{A}(\Phi(\frac{3}{4}))) \left(\f(\q_{i+1})-\f(\Phi(\frac{1}{2})) +\Bn_{i+\frac{3}{4}} \left( \q_{i+1}-\Phi(\frac{1}{2})\right) \right) \\
	-\frac{1}{3}\text{sign}(\mbf{A}(\Phi(\frac{1}{2}))) \left(\f(\q_{i+1})-\f(\q_i) +\Bn_{i+1/2} \left(\q_{i+1}-\q_i) \right)\right).
\end{array}
\ee
Note that the major drawback in the previous expression is that the complete eigenstructure of the matrix $\mbf{A}$ (\ref{eq.ExtendedJacobian}) is required since 
$ \text{sign}(\mbf{A}) = \mbf{R}\, \text{sign}(\mbf{\Lambda})\, \mbf{R}^{-1}$. 
However, on the other hand, the Osher-Romberg method is very little dissipative and is stable under the standard CFL condition.

\subsubsection{Well-balanced HLL scheme}
\label{ssec.hll1d}

Following \cite{Castro2012}, the standard HLL scheme can be written in the form \eqref{eq.SemidiscreteScheme_1} with \eqref{eq.FIrstOrderDgeneric}, where the numerical viscosity term is given by
\be
\label{eq.hll_matrixform}
\Vn_{i+\halb}(\q_{i+1}-\q_i)= \alpha^0_{i+\halb} \mathbf{I}_{i+\halb} (\q_{i+1}-\q_i) + \alpha^1_{i+\halb} \mathcal{R}_{i+\halb}, 
\ee
where $\mathbf{I}_{i+\halb}$ is the identity matrix, 
\be 
\mathcal{R}_{i+\halb}=\f(\q_{i+1})-\f(\q_i)+\Bn_{i+\halb}(\q_{i+1}-\q_i)
\ee
and
\be 
\label{eq.alphaHll}
\begin{aligned}
	\alpha^0_{i+\halb} \!=\! \frac{S^R_{i+\halb} |S^L_{i+\halb}| - S^L_{i+\halb} |S^R_{i+\halb}| }{S^R_{i+\halb} - S^L_{i+\halb}}, \,  
	\alpha^1_{i+\halb} \!=\! \frac{|S^R_{i+\halb}| - |S^L_{i+\halb}| }{S^R_{i+\halb} - S^L_{i+\halb}}. 
\end{aligned}
\ee 
Here, $S^L_{i+\halb} \leq 0$ and $S^R_{i+\halb} \geq 0$ denote the minimum and the maximum of the wave speeds of the Riemann problem associated with the states $\q_i$ and $\q_{i+1}$.
To compute $S^L_{i+\halb}$ we take the minimum of zero and the eigenvalues associated to $\q_i$ and $\frac{\q_i + \q_{i+1}}{2}$;  
to compute $S^R_{i+\halb}$ we take the maximum of zero and the eigenvalues associated to $\frac{\q_i + \q_{i+1}}{2}$ and $\q_{i+1}$. 

It is clear that $\Vn_{i+\halb}(\q_{i+1}-\q_i)$ does not vanish if $\q_{i+1}$ and $\q_i$ lie on a stationary solution: $\mathcal{R}_{i+\halb}$ vanishes, but it is not the case for the term  $\alpha^0_{i+\halb} \mathbf{I}_{i+\halb} (\q_{i+1}-\q_i)$.

Here, we follow the ideas described in \cite{castro2010some} and \cite{Castro2012} to modify the viscosity term such that the resulting scheme is exactly well balanced for the stationary solutions \eqref{eq.EquilibriaConstraint}-\eqref{eq.PressureAndGravForces}. In particular $\mathbf{I}_{i+\halb}(\q_{i+1}-\q_i)$, will be replaced by $\widetilde{ \mathbf{I}}_{i+\halb}(\q_{i+1}-\q_i)$ that vanishes when a stationary solution is considered. Here we consider the following expression for $\widetilde{\mathbf{I}}_{i+\halb}(\q_{i+1}-\q_i)$: 
\be
\label{viscosity_HLL}
\widetilde{\mathbf{I}}_{i+\halb}(\q_{i+1}-\q_i) = \left(
\begin{array}{c}
	b^{i+\halb}_2  \left( \frac{ \rho}{\gamma P} \right)_{i+\halb} \\[3pt] 
	\Delta \left(r \rho u \right )_{i+\halb} \\[3pt] 
	b^{i+\halb}_2  \left( \frac{  \rho}{\gamma P} \right)_{i+\halb} \left(v\right)_{i+\halb} \\[3pt] 
	b^{i+\halb}_2   \left( \frac{ \rho}{\gamma P} \right)_{i+\halb} \left(z\right)_{i+\halb} \\[3pt] 
	0
\end{array}
\right),
\ee
where $b^{i+\halb}_2$ is given in \eqref{b2}, $\left( \frac{ \rho}{\gamma P} \right)_{i+\halb}= \frac{\rho_{i+1}+\rho_i}{\gamma (P_{i+1}+P_i)}$, 
$\Delta \left(r \rho u \right )_{i+\halb}= (r \rho u)_{i+1}-(r \rho u)_{i}$, $\left(v\right)_{i+\halb}=\frac{v_{i+1}+v_i}{2}$,  $\left(z\right)_{i+\halb}=\frac{z_{i+1}+z_i}{2}$, being
$z \!=\!\! \de{\!\left(r u \left ( \rho E + P \right) \right )}{q_2}$.

Following \cite{castro2010some}  and  \cite{Castro2012} $\widetilde{\mathbf{I}}_{i+\halb}(\q_{i+1}-\q_i)$ is obtained as follows: we start by computing the eigenstructure of the extended Jacobian matrix $\mbf{A}$ at the equilibrium:
\be
\label{eq.Aequilibrium}
\mbf{A}(\Q) = \left (  \begin{array}{ccccc}
	0 & 1 & 0 & 0 & 0 \\ 
	\!\! \!r \de{P}{q_1} \!& 0 &\!\! r \de{P}{q_3}  \!& \! r \de{P}{q_4} & \!\! \rho \frac{G m_s}{r} + \rho v^2 \!\! \!\\ 
	0 & v & 0 & 0 & 0 \\ 
	0  & \!\!\de{\left(r u \left ( \rho E + P \right) \right )}{q_2} \!\!& 0  & 0 & 0 \\ 
	0 & 0 & 0 & 0 & 0
\end{array}   \right ).
\ee 
In this situation the eigenstructure is easy to compute: let $\mbf{R}$ the matrix of the right-eigenvectors and $\mbf{\Lambda} = \text{diag}(\lambda_1, \lambda_2, \dots, \lambda_5)$ the diagonal matrix of the eigenvalues of (\ref{eq.Aequilibrium}). In particular we have
\be
\mbf{\Lambda} = \text{diag} \left ( \frac{\rho u + \sqrt{\gamma \rho P}}{\rho}, \frac{\rho u - \sqrt{\gamma \rho P}}{\rho}, u, u, 0 \right ), \ \text{ with } u = 0.
\ee 
Then $\widetilde{\mathbf{I}}_{i+\halb}(\q_{i+1}-\q_i)$ is given by 
\be
\widetilde{\mathbf{I}}_{i+\halb}(\q_{i+1}-\q_i) = \mbf{R_{i+\halb}} \widetilde{\mbf{\Lambda}} (\mbf{R_{i+\halb}})^{-1} (\q_{i+1} -\q_i),
\ee
where
\be 
\widetilde{\mbf{\Lambda}} =\left( 
\begin{array}{ccccc}
	1 & 0 & 0 & 0 & 0 \\ 
	0 & 1 & 0 & 0 & 0 \\ 
	0 & 0 & 0 & 0 & 0 \\ 
	0 & 0 & 0 & 0 & 0 \\ 
	0 & 0 & 0 & 0 & 0
\end{array}     
\right).
\ee
Note that $\widetilde{\mbf{\Lambda}}$ is a diagonal matrix composed of $0$ and $1$, where the $0$ elements on the diagonal correspond to the zero eigenvalues at the stationary solution. The final expression \eqref{viscosity_HLL} is obtained  considering the following relation that it is derived from \eqref{prop-B}:
\be
& r_{i+\halb} \!\left(\! \left(\de{P}{q_1}\right)_{i+\halb} \!\!\Delta q_{1,i+\halb} + \left(\de{P}{q_3}\right)_{i+\halb} \!\!\Delta q_{3,i+\halb} \right.  \\
&  + \left. \left(\de{P}{q_4}\right)_{i+\halb} \!\!\Delta q_{4,i+\halb}  + \left( \rho \frac{G m_s}{r^2} - \frac{\rho v^2}{r} \right)_{i+\halb} \!\!\Delta r_{i+\halb} \!\right)\!= b^{i+\halb}_2 . \!\!\!\!\!\!\!\!
\ee 
Finally, we would like to note that a similar HLL scheme could also be obtained within the framework of path-conservative HLLEM methods recently proposed by \cite{Dumbser2015}, in which according to
\cite{munz91} the intermediate HLL state is assumed to be linear rather than constant.

\subsection{2nd order well balanced reconstruction}
\label{ssec.1d_2ndorder}
Let us recall the numerical scheme presented in \eqref{eq.SemidiscreteScheme} considering 
the space-time conservation form of the PDE
\be 
\label{eq.DiscreteScheme}
\begin{aligned} 
	\Q^{n+1}_i = \Q^{n}_i 
	& - \frac{\Delta t}{\Delta r}  \left(  \D_{i\!-\!\frac{1}{2}}^+\left( \q_{i-\frac{1}{2}}^{n{^+},-}, \q_{i-\frac{1}{2}}^{n{^+},+}   \right)  +  \D_{i\!+\!\frac{1}{2}}^-\left( 
	\q_{i+\frac{1}{2}}^{n{^+},-},  \q_{i+\frac{1}{2}}^{n{^+},+}\!  \right) \, \right )    \!\!\!\!\!\!\!\!\!\!\!\!\\
	& -  \frac{\Delta t}{\Delta r} \int_{r_{i-\frac{1}{2}}}^{r_{i+\frac{1}{2}}} \de{}{r}\f \left( \q_i^{n{^+}}(r) \right) dr  \\
	& -  \frac{\Delta t}{\Delta r} \int_{r_{i-\frac{1}{2}}}^{r_{i+\frac{1}{2}}} \B_1\left( \q_i^{n{^+}}(r) \right)\de{}{r}\left( \q_i^{n{^+}}(r) \right) dr,
\end{aligned} 
\ee 
where $\q_i^n(r,t)$ is the approximation of the conserved variables inside cell $I_i$ at time $t^n$,  $ \q_{i-\frac{1}{2}}^{n^{+},+}(t)=\q_i^n(r_{i-1/2},t^{n+1/2})$ and  $ \q_{i+\frac{1}{2}}^{n^{+}-}(t)=\q_i^n(r_{i+1/2},t^{n+1/2})$, that is the evaluation of $\q_i^n(r,t)$ at the two boundaries of cell $I_i$ at the time-midpoint of $[t^n, t^{n+1}]$.
We would like to underline that in order to obtain a second order scheme $\q_i^n$ should be a second order reconstruction of the cell averages $\Q_{i-1}^n, \Q_{i}^n, \Q_{i+1}^n$.	

According to \cite{Pares2006} and \cite{Castro2006}, scheme \eqref{eq.DiscreteScheme} is well balanced if both, the underlying first order scheme and the reconstruction operator are well balanced, and all the integrals that appear in \eqref{eq.DiscreteScheme} are computed exactly.   Therefore, in order to define a second order scheme, a second order well balanced reconstruction operator should be defined. 	 

The most popular way to define a second order reconstruction  operator  is based on the MUSCL method introduced by van Leer in \cite{leer5} joint with the \textit{minmod} limiter. He proposed to reconstruct $\q_i^n$ using a linear polynomial in space and time as follows
\be
\mathcal{P}_i^{n} (r,t) = \Q_i^n +  \frac{\Delta \Q_i^n}{\Delta r}(r-r_i) + \partial_t\Q_i^n (t -t^n),
\ee 
where
\be
\Delta \Q_i^n = \text{minmod} \left ( \Delta \Q_{i-1/2}^n, \Delta \Q_{i+1/2}^n \right)
\ee 
with {$ \Delta \Q_{i-1/2}^n = \Q_i^n - \Q_{i-1}^n$}, $ \Delta \Q_{i+1/2}^n = \Q_{i+1}^n - \Q_i^n $ and 
\be
\text{minmod}(a,b) = \begin{cases}
	0, \quad \text{if }  ab \le 0\\
	a, \quad \text{if } |a| < |b|\\
	b, \quad \text{if } |a| \ge |b|.\\
\end{cases} 
\ee 
It is clear that the standard MUSCL method is only well balanced for linear stationary solutions, which is not the case here. In this paper we therefore follow the strategy proposed in \cite{Castro2008},  where the reconstruction operator is defined as a combination of a smooth stationary solution together with a standard reconstruction operator to reconstruct the {\it fluctuations} with respect to the  given stationary solution, that is
\be
\label{eq.FinalReconstruction}
\q_i^n(r,t) = \Q^E_i(r) + \mathcal{P}^{f}_i(r,t), \quad r \in I_i, \ t \in [t^n, t^{n+1}],
\ee 
where $\mathcal{P}^{f}_i(r,t)$ is the standard MUSCL reconstruction operator applied to the fluctuations around the stationary solution at every cell of the stencil. Thus, if we define
\be
\label{eq.FluctuationForRec}
\Q^{f\!,n}_i \!= \Q_i^{n} \!- \Q^{E\!}_i\!, \quad  \Q^{f\!,n}_{i-1} \!= \Q^n_{i-1} \!- \Q^{E\!}_{i-1}\!, \quad  \Q^{f\!,n}_{i+1} \!= \Q^n_{i+1} \!- \Q^{E\!}_{i+1}\!,
\ee
then, $\mathcal{P}^{f}_i(r,t)$ is defined as follows:
\be
\label{reconstruction_fluctuation}
\mathcal{P}^{f,n}_i (r,t) = \Q_i^{f,n} +  \frac{\Delta \Q_i^{f,n}}{\Delta r}(r-r_i) + \partial_t\Q_i^n (t -t^n),
\ee 
where
\be
\Delta \Q_i^{f,n} = \text{minmod} \left ( \Delta \Q^{f,n}_{i-1/2}, \Delta \Q^{f,n}_{i+1/2} \right)
\ee 
with
\be
\Delta \Q^{f,n}_{i-1/2} = \Q_i^{f,n} - \Q_{i-1}^{f,n}, \quad  \Delta \Q_{i+1/2}^{f,n} = \Q_{i+1}^{f,n} - \Q_i^{f,n}.
\ee 
Note that we have replaced  $\partial_t\Q^{f,n}_i (t -t^n)$ by $\partial_t\Q_i^n (t -t^n)$ in \eqref{reconstruction_fluctuation} as $\partial_t\Q^E_i=0$.
It is clear from its construction that the reconstruction operator is exactly well balanced, and it is second order accurate for non-stationary solutions as $\Q^E(r)$ is a smooth stationary solution.  The term $\partial_t \Q_i^n$ indicates the time derivative of $\Q$ and  it can be computed using a discrete version of the governing equation
\be
& \partial_t \Q_i^n \!= - \frac{\f(\q_{i+1/2}^{n,-}) - \f(\q_{i-1/2}^{n,+})}{\Delta r} - \frac{\Bn_{i}(\q_{i+1/2}^{n,-} - \q_{i-1/2}^{n,+})}{\Delta r}, \\
& \q_{i\pm 1/2}^{n,\mp} = \q_i(x_{i\pm1/2}^{\mp}, t^n),
\ee 
where the fluxes have been approximated by a central finite difference with respect to the cell center $r_i$, and $ \Bn_{i}(\q_{i+1/2}^{n,-} - \q_{i-1/2}^{n,+})$ is obtained in the same way of \eqref{eq-B},\eqref{b2},\eqref{eq-b3},\eqref{eq-b4}, where by replacing for example $\q_i$ and $\q_{i+1}$ by $\q_{i+1/2}^{n,-}$ and $\q_{i-1/2}^{n,+}$ respectively, and using as central value the cell average one obtains
\be
& \frac{\Bn_{i}(\q_{i+1/2}^{n,-} - \q_{i-1/2}^{n,+})}{\Delta r}= \left( b^{i}_1 \ b^{i}_2 \ b^{i}_3 \ b^{i}_4 \ b^{i}_5 \right)^T \ \text {with} \\
& b^{i}_1=b^{i}_5=0, \quad b^{i}_3= \rho_i u_i v_{i}, \quad b^{i}_4= \rho_i u_{i} \frac{G m_s}{r_{i}}, \\
& b^{i}_2 = r_{i} \left(P_{i+1/2}^{f,n,-} - P_{i-1/2}^{f,n,+} \right) + \left( r_i \rho^E_i (\zeta_r)^f_{i} + r_i \rho^f_{i} (\zeta_r)_{i}  \right). \\
\ee

The last ingredient for a second order scheme is the computation of the integrals in $(\ref{eq.DiscreteScheme})$: the first  one can be computed exactly
\be
\int \limits_{r_{i-1/2}}^{r_{i+1/2}} \de{}{r}\f \left( \q_i(r,t) \right) dr = \f(\q_{i+1/2}^{n^{+},-}) - \f(\q_{i-1/2}^{n^{+},+}).
\ee
Note that this first integral vanishes for stationary solutions with $u=0$.
The second integral is more sophisticated, and it is not easy to compute it exactly, except in some particular situations.  Therefore we will use a quadrature formula to approximate this integral, but this must be done carefully to maintain the well balanced property of the scheme: effectively, a wrong choice in the quadrature formula will destroy all the work we have done up to now in order to define a well balanced scheme.  Here we proceed as follows: first we express the particular form of the reconstruction operator: $\q_i^n(x,t) = \Q^E_i(x) + \mathcal{P}^{f}_i(x,t)$ and we use the fact that
\be 
\int_{r_{i-1/2}}^{r_{i+1/2}} \B_1( \Q^E_i(r) )\de{\Q^E_i(r)}{r} dr =\0.
\ee
Here, we only show the details for the second component of 
\be
\int_{r_{i-1/2}}^{r_{i+1/2}} \B_1( \q_i^n(r) )\de{\q_i^n(r)}{r} dr ,
\ee
\be
& \quad \ \int_{r_{i-1/2}}^{r_{i+1/2}} r \left [ \de{P}{r} + \rho \left ( \frac{G m_s}{r^2} - \frac{v^2}{r}  \right) \right]  dr  \\
&= \int_{r_{i-1/2}}^{r_{i+1/2}} r \left [ \de{(P^E + P^f)}{r} + \left( \rho^E + \rho^f \right ) \left ( \zeta^E + \zeta^f  \right)_r \right]  dr\\
&= \int_{r_{i-1/2}}^{r_{i+1/2}} r \de{P^f}{r} + r\rho^E\zeta_r^f + r \rho^f \zeta_r dr. \\[2pt]
\ee
Now, the mid-point quadrature formula is used to ensure second order accuracy obtaining that	
\be
&\quad \ \int_{r_{i-1/2}}^{r_{i+1/2}} r \left [ \de{P}{r} + \rho \left ( \frac{G m_s}{r^2} - \frac{v^2}{r}  \right) \right] dr \\
& \approx \Delta r \Biggl [ \, r_i \, \left(\Delta P^f\right)_i \, +\, \left(r \rho^E \right)_i \, \left(\zeta_r^f\right)_i \, + \, \left(r \rho^f\right)_i \, \left(\zeta_r\right)_i   \Biggr ]\,, \\
& \text{where } \ \ \left(\Delta P^f\right)_i = \frac{P^{f,-}_{i+1/2} - P^{f,+}_{i-1/2}}{\Delta r}, \\
& \ (\zeta_r)_i = \frac{G m_s}{r_i^2} - \frac{v_i^2}{r_i}, \quad
\ (\zeta_r^f)_i = \frac{v_i^{E^2}}{r_i} - \frac{v_i^2}{r_i}.
\ee 
It is clear that this approximation is second order accurate and, moreover, will vanish for stationary solutions \eqref{eq.EquilibriaConstraint}-\eqref{eq.PressureAndGravForces}.
For the third and fourth component we could perform the same procedure, but, again, as both trivially vanish when $u=0$, we could use directly  the mid-point rule.

Finally,  note  that $r^\pm_{i+1/2}=r_{i+1/2}$ and therefore $\Delta r_{i+1/2} =0$. Therefore $\Bn_{i+1/2}(\q_{i+\frac{1}{2}}^{n{^+},+}-  \q_{i+\frac{1}{2}}^{n{^+},-})$ reduces to
\[
\Bn_{i+1/2}(\q_{i+\frac{1}{2}}^{n{^+},+}-  \q_{i+\frac{1}{2}}^{n{^+},-})=(0,b_2^{i+1/2},0,0,0)^T
\]
where
\[
b_2^{i+1/2}=r_{i+1/2} \Delta P^f_{i+1/2}= r_{i+1/2} \left( P^{f,+}_{i+1/2}-P^{f,-}_{i+1/2}\right).
\]

\section{Numerical results in one dimension}
\label{sec.Results1d}

First of all, we show the ability of both schemes to  preserve a wide class of stationary solutions and we also report the convergence tables for some smooth solutions. Then, we test both methods  with some classical Riemann problems, and finally we study their behavior in capturing small perturbations around the equilibrium.

\subsection{Stationary solutions with constant pressure}

Simple, but non trivial, stationary solutions of the Euler equations can be obtained by considering velocities as in (\ref{eq.EquilibriaConstraint}) and a constant pressure $P$. 
It is easy to verify that under these conditions for any density profile the velocity in the angular direction $v$ must satisfy 
\be
v = \sqrt{\frac{Gm_S}{r}},
\ee 
while $u=0$. 
For the numerical simulations we consider a spatial domain $r \in [1,2]$, $ G = 1$, $m_s =1$, $ \gamma = 1.4$, $P =1$ and two density profiles:
\be
\label{eq.EqPcnst_cont}
\rho_1 = r, 
\ee
\vspace{-4mm}  
\be
\label{eq.EqPcnst_disc}
\rho_2 = \begin{cases} 1, \ &\text{ if } r< 1.5 \\
	0.1, \ &\text{ if } r \ge 1.5.
\end{cases}
\ee  

In Table~\ref{tab.EquilibriumPreservation} we report the errors between the exact equilibrium and the numerical solution obtained with both schemes  using a hierarchy of meshes and long term time integration. We can notice that all the errors are of the order of machine precision and no significant differences can be noticed between the two fluxes.
Moreover the method is perfectly well balanced both with continuous and discontinuous density profiles, as expected.

\begin{table}	
	\centering
	\caption{Constant pressure equilibrium. The following results show the capability of the schemes to preserve equilibria both for a hierarchy of meshes for a fixed time $t=1$ (on the left) and for a fixed mesh ($N=64$ cells) and for increasing computational times. 
		The table on the top refers to the $L_1$-norm error between the continuous $\rho_1$ profile and the table on the bottom refers to the discontinuous $\rho_2$ profile. Data have been obtained using either the Osher or HLL flux (and no significant differences have been noticed). }
	\begin{tabular}{cc|cc}
		\hline 
		\multicolumn{2}{c|}{$ \text{tend} = 1 $} &\multicolumn{2}{c}{  $N = 64$}        \\ 
		\hline
		N & \qquad $E_{\rho}$ - Osher  \qquad \qquad & \quad time\qquad  &  $E_{\rho}$ - Osher   \\ 
		\hline
		64   & 9.54E-17 &  1    & 9.54E-17   \\
		128  & 9.54E-17  &  2    & 2.36E-16   \\  
		256  & 6.49E-16  &  5    & 8.85E-16   \\
		512  & 6.23E-16  & 10   & 1.67E-15  \\	
		1024 & 1.21E-15  & 50   & 6.24E-17    \\		
		\hline 
	\end{tabular}	
	\medskip \hspace{-0.6cm}
	\begin{tabular}{cc|cc} 
		\hline 
		\multicolumn{2}{c|}{$ \text{tend} = 1 $} &\multicolumn{2}{c}{  $N = 64$}        \\ 
		\hline
		N & \qquad $E_{\rho}$ - HLL\quad  \qquad \qquad & \quad time\qquad  &  $E_{\rho}$ - HLL \   \\ 
		\hline
		64   &  8.45E-18  &  1    &  8.45E-18  \\
		128  &  1.38E-16  &  2    &  1.19E-17 \\  
		256  &  5.54E-16  &  5    &  6.71E-16 \\
		512  &  2.64E-15  &  10   &  2.42E-15 \\			
		1024 &  5.05E-16  &  50   &  1.13E-13  \\		
		\hline 
	\end{tabular}
	\label{tab.EquilibriumPreservation}
\end{table}

\subsection{General equilibrium}

Using the equilibrium relation between the pressure and the gravitational forces in (\ref{eq.PressureAndGravForces}) and $\zeta$  given by  (\ref{eq.zeta}), we obtain another class of stationary solutions of the Euler equations
\be
 \rho = \rho_0 e^{-\zeta(r)}, \quad 
 P = \rho + P_0, \quad 
 v = \sqrt { r \left( \frac{G m_s}{r^2} - \zeta_r \right) }. 
\ee 

We have applied  both schemes to two different choices of $\zeta$ obtaining always a well balanced result. 
Table~\ref{tab.EquilibriumPconstant} shows the $L_1$-norm error for the density between the equilibrium and the numerical solution in the case
\be
\zeta = kr, \ k = -1,\ 
\rho = \rho_0 e^{-kr}\!, \ \rho_0 = 1, \
P = \rho + P_0, \ P_0 = 1.
\ee 
Again, both methods are able to exactly preserve these non-trivial equilibria. 

\begin{table}		
	\centering
	\caption{General equilibrium. $L_1$-norm error  for the density between the exact and the numerical solution. On the left we have the error for different meshes at $t=1$ and  on the right we show the error for a given mesh ($N=64$ cells) for different computational times. }
	\begin{tabular}{cc|cc} 
		\hline 
		\multicolumn{2}{c|}{$ \text{tend} = 1 $} &\multicolumn{2}{c}{  $N = 64$}        \\ 
		\hline
		N & \qquad $E_{\rho}$ - OSHER  \qquad \qquad & \quad time\qquad  &  $E_{\rho}$ - HLL   \\ 
		\hline
		64   &  6.28E-15  &  1    & 5.03E-15   \\
		128  &  1.17E-14  &  2    & 1.01E-14   \\ 
		256  &  1.70E-14  &  5    & 2.65E-14\\
		512  &  2.15E-14  & 10   & 7.21E-14   \\	
		1024 &  3.19E-14  &  50   & 3.07E-12   \\		
		\hline 
	\end{tabular}	
	\label{tab.EquilibriumPconstant}
\end{table}

\subsection{Order of convergence}

To study numerically the order of convergence  of both schemes we have considered the following equilibrium situation 
\be
 \rho = 1, \quad u = 0,  \quad P = 1,  \quad 
 v = \sqrt { r \left( \frac{G m_s}{r^2} - \zeta_r \right) },  
\ee 
and at the initial time, we have added a small perturbation (with a Gaussian profile) to the velocity field
\be
\tilde{u} = u + 10^{-5} \text{exp} \left ( \frac{ - 0.5 \left(r-1.5\right)^2 }{0.01} \right ) ,   \\
\tilde{v} = v + 10^{-5} \text{exp} \left ( \frac{ - 0.5 \left(r-1.5\right)^2 }{0.01} \right ) .
\ee 

We have computed a reference solution using our method with the Osher-type flux on a fine mesh ($N = 2^{13} = 8192$ ). 
In Table~\ref{tab.convergenceTable1d} we report the $L_1$ error norms for  the density $\rho$ with respect to our reference solution and both numerical schemes achieve second order of convergence.

\begin{table} 
	\caption{Perturbation around a stationary solution.  The reference solution has been obtained with the second order Osher-type scheme over $2^{13}$ cells. $L_1$-norm errors for  $\rho$ at time $t=0.1$ are shown:  on the left we report the result obtained using the Osher-type flux and on the right using the HLL-type flux.} 
	\begin{center} 	
		\begin{tabular}{ccc|ccc} 
			\hline 
			\multicolumn{3}{c}{Osher $\mathcal{O}2$}& \multicolumn{3}{c}{HLL $\mathcal{O}2$  }       \\
			\hline		
			N   &  $ \epsilon (\rho) $ & $\mathcal{O}(L_1)$ \qquad & N &  $ \epsilon (\rho) $ & $\mathcal{O}(L_1)$  \\ 
			\hline
			16  & 1.59E-07 &   -    &   16    & 1.16E-07 &   -    \\
			32  & 3.82E-08 &  2.06  &	32    & 2.90E-08 &  2.01  \\	
			64  & 9.50E-09 &  2.00  & 	64    & 7.22E-09 &  2.00  \\	
			128 & 2.31E-09 &  2.04  &	128   & 1.77E-09 &  2.03  \\	
			256 & 5.72E-10 &  2.01  &	256   & 4.44E-10 &  1.99  \\	
			512 & 1.45E-10 &  1.97  &   512   & 1.14E-10 &  1.96  \\				
			\hline 
		\end{tabular}		
	\end{center}
	\label{tab.convergenceTable1d}
\end{table}

\subsection{Riemann Problem}

\begin{figure*}
\begin{tabular}{cc} 
	\includegraphics[width=0.4\textwidth]{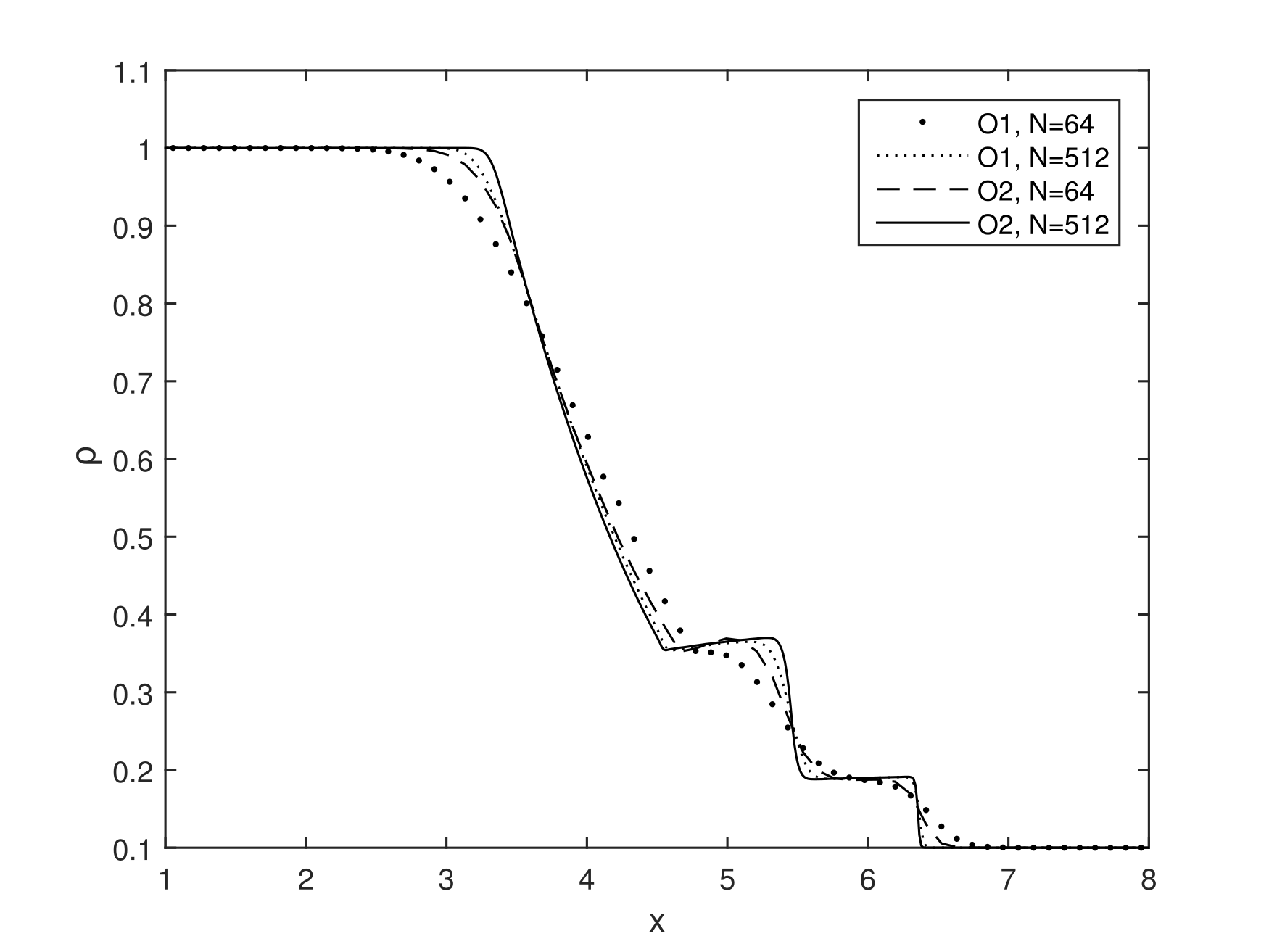} &
	\includegraphics[width=0.4\textwidth]{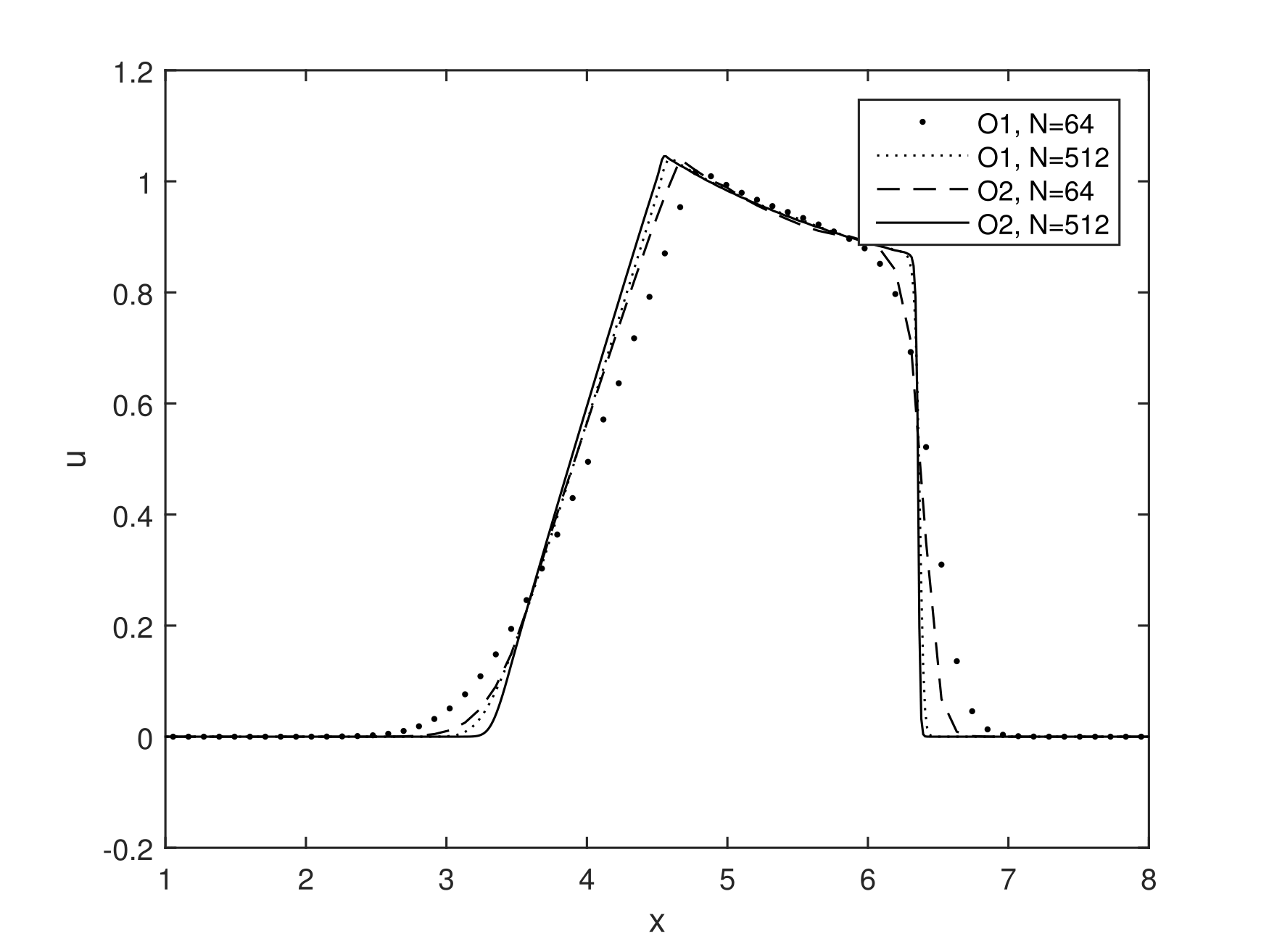} \\
	\includegraphics[width=0.4\textwidth]{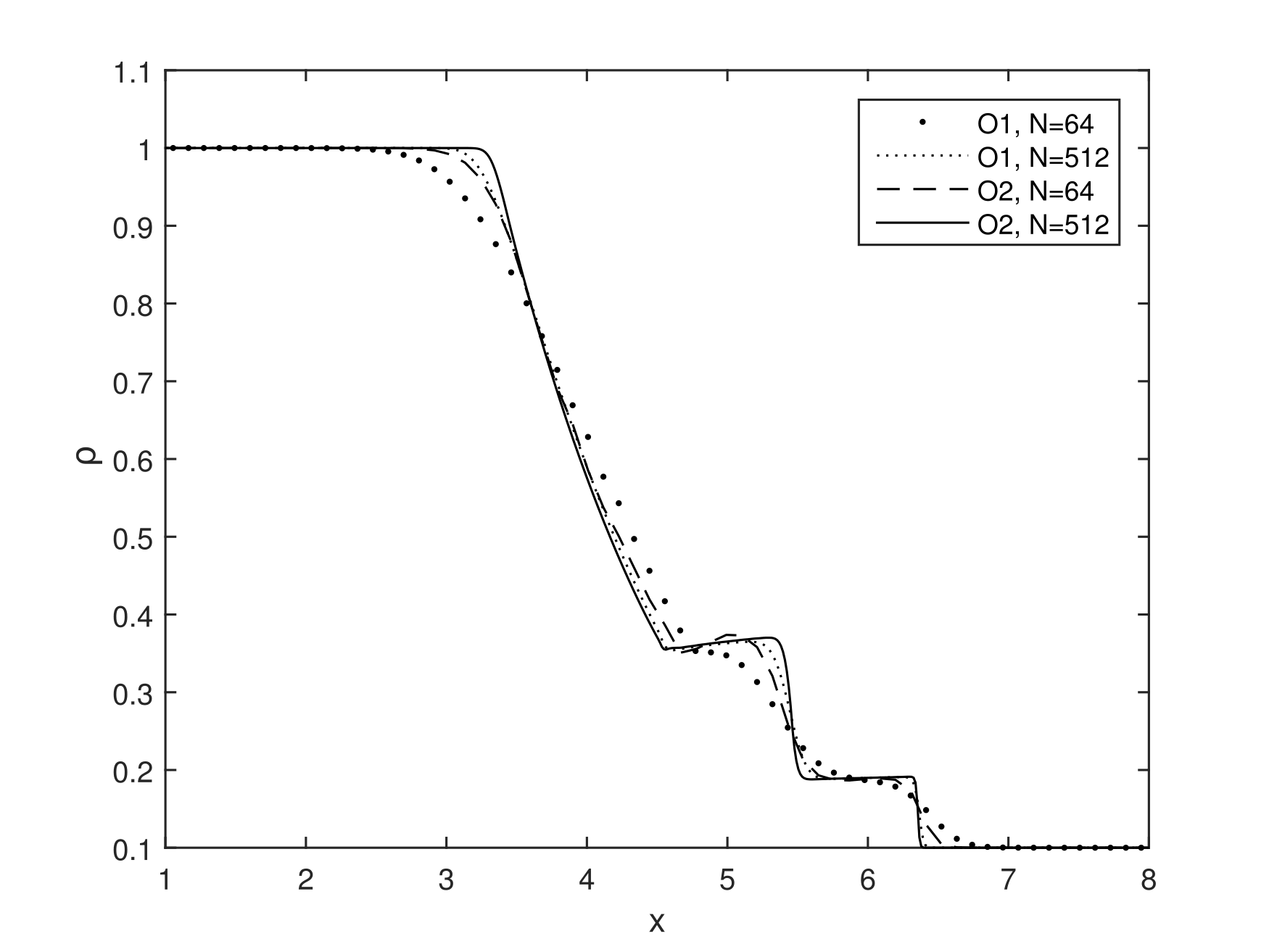} & 
	\includegraphics[width=0.4\textwidth]{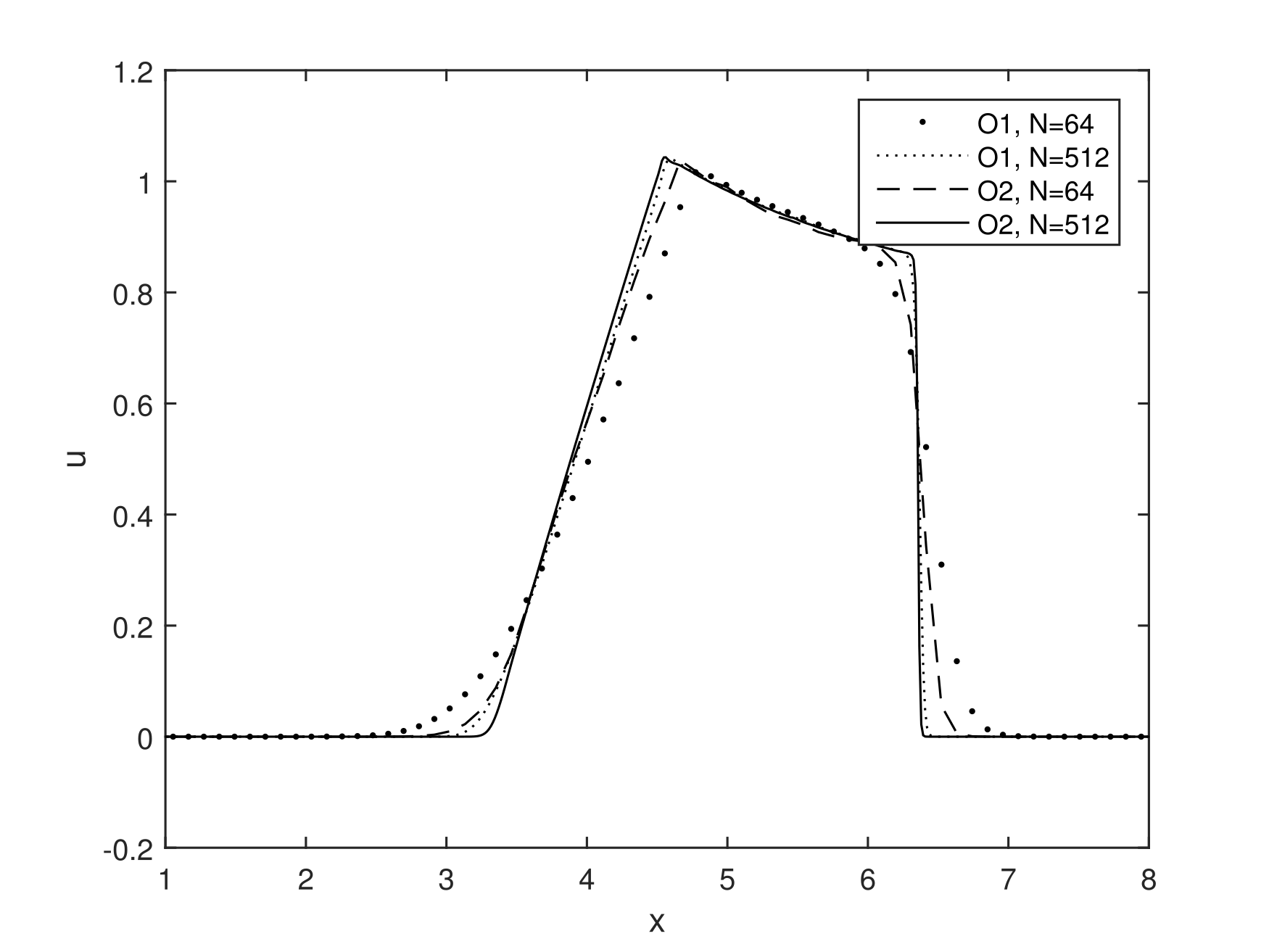} 
	\end{tabular} 
	\caption{Riemann problem at final time $t_f=1$. On the top the density and velocity profiles obtained using the HLL scheme and on the bottom the profiles obtained with the Osher-type flux. We have employed two meshes: a coarse one with $N=64$ elements and a fine one with $N=512$ elements. Moreover, we have compared the first and second order schemes.}
	\label{fig.RP1d}
\end{figure*}

To show that our method is accurate even far away from an equilibrium, we consider as initial condition a classical Riemann problem with non-vanishing angular velocity 
\begin{equation*} 
\begin{aligned} 
\rho_L = 1.0, u_L= 0, v_L=\sqrt{\frac{G m_s}{r}}, P_L = 1.0, r= r, \, \, 1 \le r \le 4.5, \\
\rho_R = 0.1, u_R= 0, v_R=\sqrt{\frac{G m_s}{r}}, P_R = 0.1, r= r, \, \, 4.5 < r \le 8,
\end{aligned} 
\end{equation*} 
and we compute the solution by employing the schemes set up to preserve the equilibrium in (\ref{eq.EqPcnst_disc}). 
We report the results obtained with the first and second order scheme and with the HLL and Osher-type flux in Figure (\ref{fig.RP1d}). Note that both schemes produce quite similar results.

\subsection{Evolution of perturbations} 

Following the idea presented in \cite{KM15_630} we have tried to study small perturbations around the equilibrium.
We have considered a density profile as in (\ref{eq.EqPcnst_cont}) and we have imposed a periodic perturbation on the velocity $u$ through the left boundary conditions, by imposing
\be
u_0 = A \sin \left ( 6 \frac{2\pi t^n}{t_f} \right ), \ \  t_f = 1.
\ee

\begin{figure}
	\includegraphics[width=0.9\columnwidth]{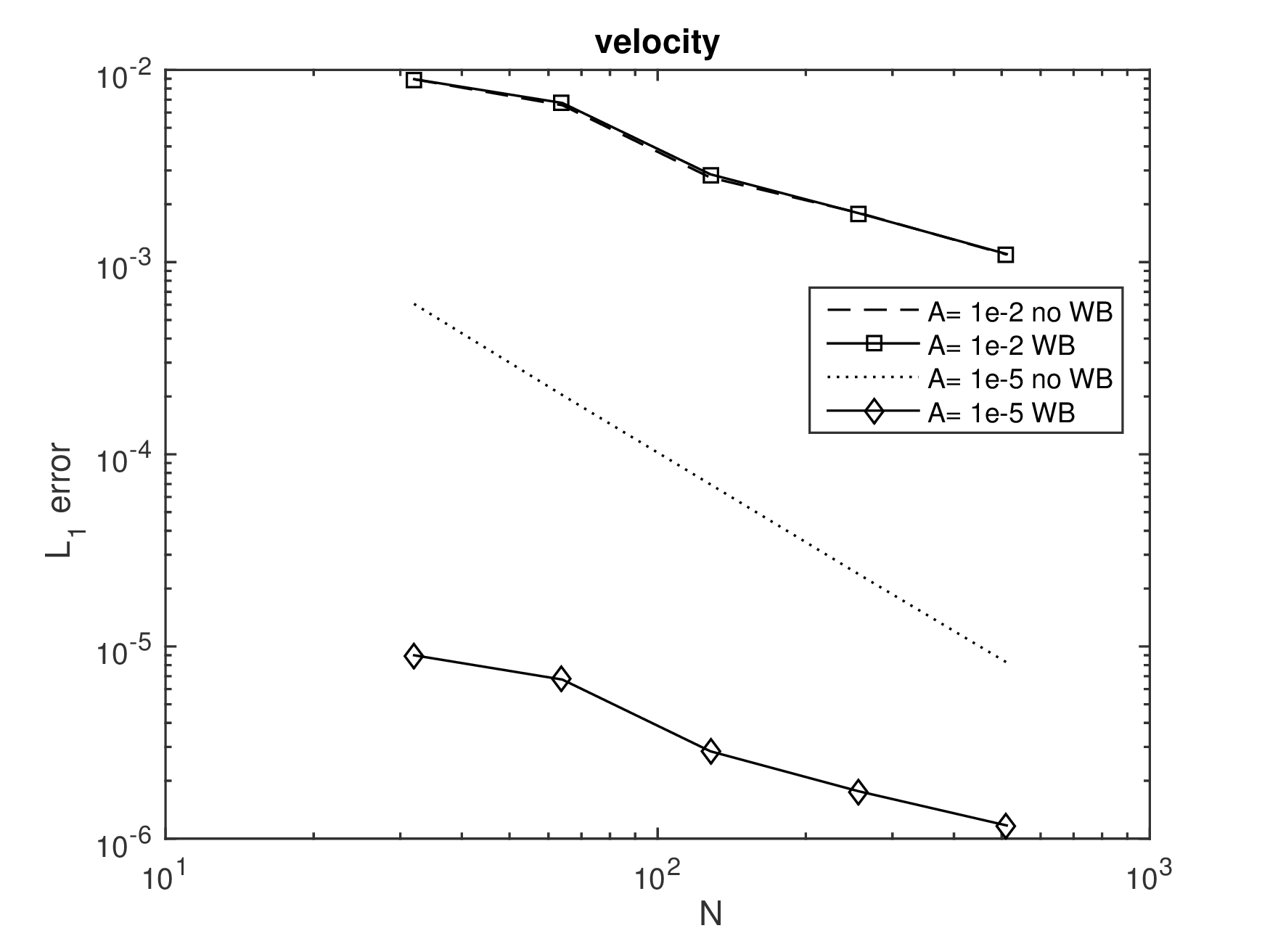} 
	\caption{Error in the $L_1$ norm between a reference solution and the numerical solutions computed with the well balanced HLL method and a second order non well balanced scheme. Well balanced and non well balanced methods perform equally well for large perturbations, while well balanced schemes perform significantly better for the small perturbation problem.}
	\label{fig.velPerturbations}
\end{figure}

Two situations are analyzed. First we consider a big perturbation, with $A=10^{-2}$ and we simulate the evolution using the second order well balanced HLL scheme and a standard second order HLL scheme using a hierarchy of grids with increasing number of cells. A reference solution computed with the second order well balanced HLL method is also considered using a fine grid composed of $N=2048$ cells. Figure (\ref{fig.velPerturbations}) shows the errors for the different meshes. Note that in this case no big differences are visible between the well balanced and not well balanced schemes as the perturbations are so large so that shocks are quickly generated and the solution is far away  from the stationary profile.  The situation changed significantly when a small perturbation is considered  
($A=10^{-5}$). In that case the well balanced method performs much better than the non well balanced scheme on the finer grid, as shown in Figure  (\ref{fig.velPerturbations}).

\section{Numerical method in two dimensions}
\label{sec.NumMethod2d}

Now, we extend our method to the two dimensional ALE context on moving nonconforming meshes. In particular we are interested in numerical schemes able to approximate accurately nontrivial 
equilibrium solutions along the radial direction given by \eqref{eq.EquilibriaConstraint}-\eqref{eq.PressureAndGravForces}. 
Hence, in general at the equilibrium $v\neq0$, which implies that $\g \neq 0$ and makes it difficult to design well balanced schemes on general meshes. 

For this reason we are going to design a numerical scheme on moving meshes that inherits the well balanced property of the previous one-dimensional scheme in the radial direction and we need to impose   some conditions on the shape of the moving cells so that the two components of the flux are not completely mixed in the computation. 
We emphasize that our numerical scheme works for completely general unstructured and nonconforming moving meshes, but it will be well balanced only if the mesh 
satisfies some special conditions. 

The rest of the section is organized as follows: first we describe the domain discretization and its time evolution due to the ALE context. 
Next, we derive the one-step path-conservative ALE scheme, and we explain where the 1D well balanced techniques are employed in order to guarantee the well balancing of the scheme even in a two dimensional moving mesh framework.

\subsection{Arbitrary-Lagrangian-Eulerian scheme}
\label{ssec.MovingDomainDiscretization}

To discretize the moving domain, we consider a nonconforming mesh $\mathcal{T}_\Omega^n$ which covers the computational domain $\Omega(\x, t^n) = \Omega^n$ at the time $t^n$ with a total number $N_E~=~N~\times~M$ of quadrilateral elements $T_i^n$, $i=1,\dots,N_E$. 
{We refer to our mesh as \textit{nonconforming} because each edge can be shared between 
\textit{more} than two elements and a node can lie on an edge not only at its extremities, i.e. we explicitly allow so-called \textit{hanging nodes}. This gives us more flexibility in the grid motion and helps to maintain a high quality mesh.}

\begin{figure}	
	\begin{center}
		\includegraphics[width=0.6\linewidth]{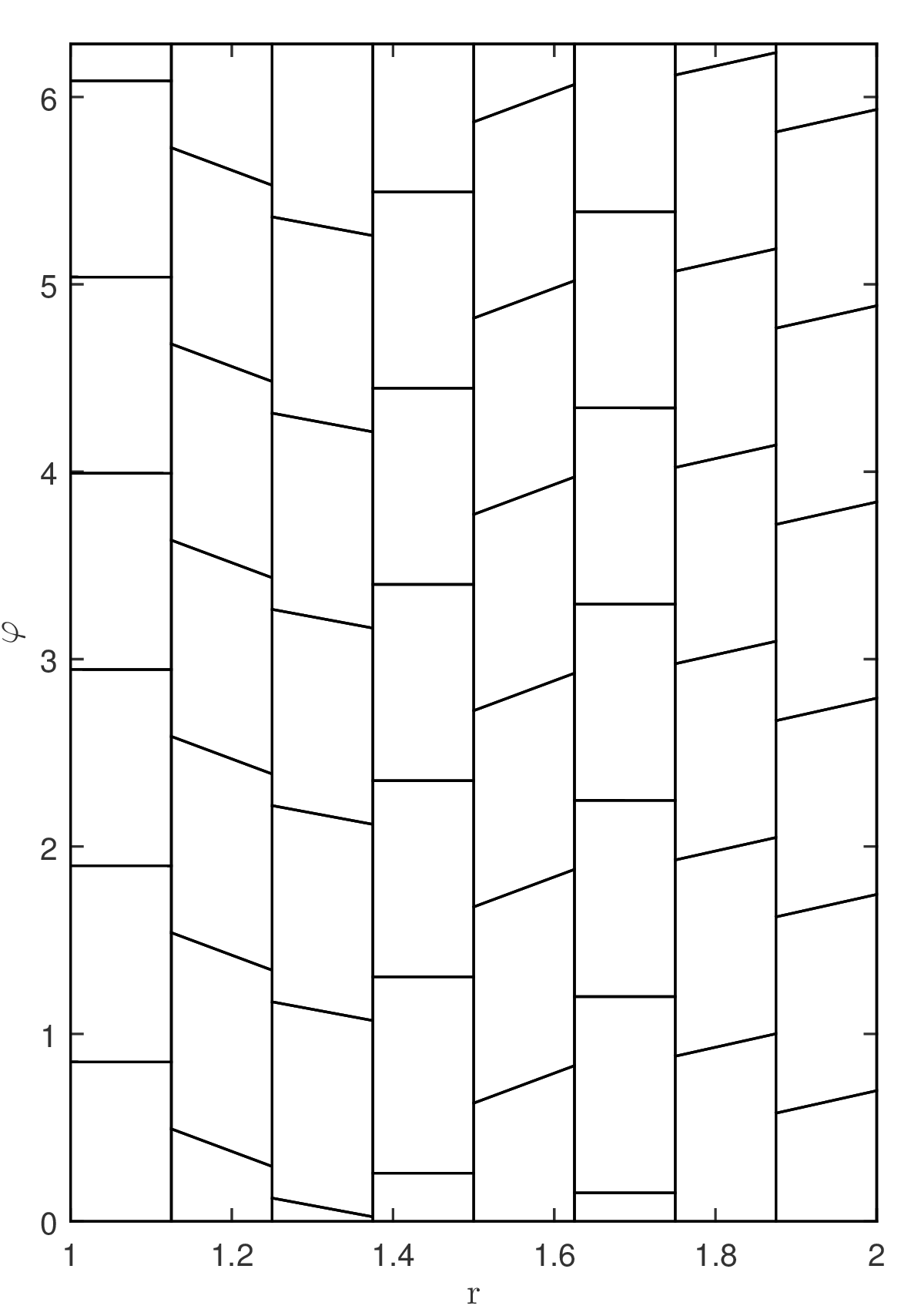}
	\end{center}	
	\caption{Example of a mesh that allows a well balanced treatment of the fluxes. Each element has two vertical edges and the other two are parallel between them. Besides the vertical edges lie on straight lines and the barycenters are aligned along $r=r_i$. Moreover the domain is periodic so that $\varphi=0$ coincides with $\varphi=2\pi$.} 
	\label{fig.acceptedMesh}
\end{figure}

The elements should satisfy the following conditions:
\begin{enumerate}
\label{item.mehsConstraints}
\item
their barycenters should be aligned along straight lines with $r = r_i,\  i\! =\! 1, \dots N$, 
\item
the two bounding edges of each element in radial direction must be aligned with $r = r_{i \pm 1/2} = const., \ i \!=\! 1, \dots N\!+\!1$,
\item the other two bounding edges must be parallel between them. 
\end{enumerate}
For example a Cartesian grid satisfies these conditions, but we could accept even something more general (which allows us to move the computational domain). See Figure~\ref{fig.acceptedMesh} for a general mesh that satisfies the above constraints. In Section \ref{ssec.2dFirstOrderMethod} these choices will be justified.

The method we are going to employ to solve \eqref{eq.EulerPolarWellBalanced} belongs to the family of the Arbitrary-Lagrangian-Eulerian (ALE) finite volume schemes. 
This kind of schemes is characterized by a moving computational mesh: at each time step the new position of all the nodes has to be recomputed according to a prescribed mesh velocity, which generally is chosen as close as possible to the local fluid velocity (as it is in the purely Lagrangian framework), but it can also be set to zero (to reproduce the 
Eulerian case), or it can be chosen arbitrarily. 

The aim of these methods is to reduce the numerical dissipation errors due to the convective terms and to capture contact 
discontinuities sharply. 
For this reason it is particularly well suited for our situation, where the gas at the equilibrium is advected with the known equilibrium velocity field 
$ \mbf{V}(\x) =\left(u^E(\x), v^E(\x)\right)$ which reads 
\be
\label{eq.exactVelocity}
u^E(\x) = 0, \quad v^E(\x) = \sqrt { r \left( \frac{G m_s}{r^2} - \zeta_r \right) }.
\ee

Note that the \apriori \ knowledge of the velocity field significantly simplifies the application of an ALE scheme: indeed, we can move the nodes following directly the exact equilibrium 
velocity, which is not affected by any physical or numerical perturbation. In general the coordinates of a node $k$ are evolved from time $t^n$ to time $t^{n+1}$ according to
\be
\x^{n+1}_k = \x^n_k + \Delta t \overline{\mbf{V}}_k^n
\ee
where $\overline{\mathbf{V}}^n_k$ is obtained using  the node solver of Cheng and Shu.  
Cheng and Shu introduced in \cite{chengshu1} and \cite{chengshu2} a very simple and general formulation to obtain the final node velocity, which is chosen to be the arithmetic average velocity among all  the contributions coming from $\mbf{V}$ evaluated at the barycenter of the Voronoi neighbors of node $k$. 

This allows us to control the movement of the mesh avoiding the violation of the above conditions: indeed the radial component of $\overline{\mbf{V}}_k^n$ will be always zero, hence nodes will slide along straight lines with $r=const.$ where the edges lie. Moreover, since the barycenters are placed on the straight lines with $r=r_i$,  all nodes lying on the same edge will move with the same velocity  maintaining the parallelism constraint between the edges. 

Moreover, the presence of known straight slide lines makes it possible to apply the algorithm described in \cite{gaburro2016direct} for a nonconforming treatment of the mesh motion: thanks to this technique we are able to preserve a high level of grid quality of the moving mesh even in the case of strong shear flows that originates in Keplerian discs due to the differential rotation. 
For all the details about the nonconforming motion of nodes along sliding lines (insertion and deletion of nodes and edges, computation of the velocity of new nodes, and flux computation in the case of more than two neighbors), we refer to \cite{gaburro2016direct}, with the only difference that in our case the sliding interfaces are prescribed \apriori \ and do not need to be automatically detected by the  algorithm.


For the sake of clarity, we briefly recall here the concept of \textit{space-time control volumes} employed in our direct ALE scheme. 

\begin{figure}
	\centering
	$\vcenter{\hbox{\includegraphics[width=0.25\textwidth]{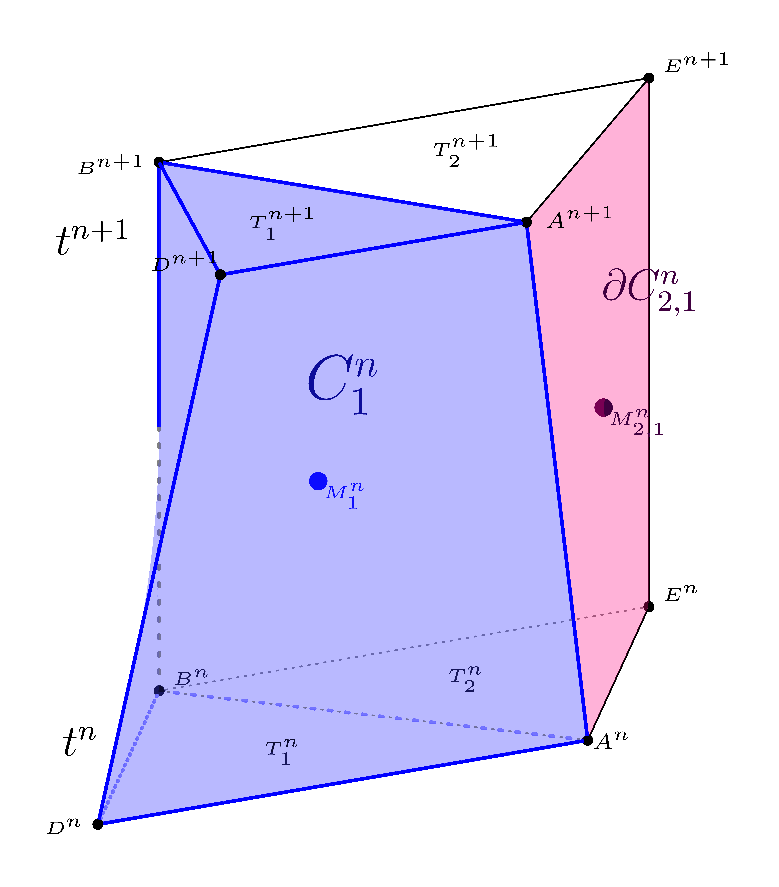}}}$ \hspace{-1cm} \hspace*{.2in}
	$\vcenter{\hbox{\includegraphics[width=0.2\textwidth]{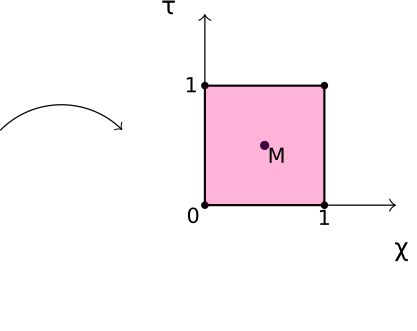}}}$
	\caption{Left. In blue we show the physical space-time control volume $C^n_1$ obtained by connecting via straight line segments each vertex of $T_1^n$ with the corresponding vertex of $T_1^{n+1}$, and its space-time midpoint $M_1^n$. In pink we show one of the lateral surfaces of $C_2^n$, $\partial C_{2,1}^n$, together with its space-time midpoint $M_{2,1}^n$. Right. The reference system $(\chi,\tau)$ adopted for the bilinear parametrization of the lateral surfaces $\partial C^n_{ij}$.}
	\label{fig.Cin}
\end{figure}

Let $T_i^n$ and $T_i^{n+1}$ denote the space control volumes respectively at time $t^n$ and $t^{n+1}$.
A space-time control volume $C_i^n$ is obtained by connecting each vertex of the element $T_i^n$ via \textit{straight} line segments with the corresponding vertex 
of $T_i^{n+1}$. For a graphical interpretation one can refer to Figure~\ref{fig.Cin}, where we have reported an example of a control volume and the parametrization of one of its the lateral space-time surfaces. A lateral space-time  surface is denoted by $\partial C_{ij}^n$ where the index $i$ refers to the element $C_i^n$ and the index $j$ refers to the neighbor $j$ of $C_i^n$.

For each control volume we have to compute the normal vectors, the areas and the space-time midpoints of all its sub--surfaces
\be
\partial C^n_i = \left( \bigcup \limits_{j} \partial C^n_{ij} \right) 
\,\, \cup \,\, T_i^{n} \,\, \cup \,\, T_i^{n+1}.
\label{dCi}
\ee
The upper space-time sub-surface $T_i^{n+1}$ and the lower space-time sub-surface $T_i^{n}$ are the simplest, since they are orthogonal to the time coordinate. The space-time unit normal vectors are respectively $\mathbf{\tilde n} = (0,0,1)$ and $\mathbf{\tilde n} = (0,0,-1)$. 
Area and barycenter can be easily computed, since $T_i^{n}$ and  $T_i^{n+1}$ are quadrilaterals. We will denote the area of $T_i^n$ with $ |T_i^n|$ and use the notation 
$\mathbf{\tilde{x}} = (r,\varphi,t)$ for the space-time coordinate vector.  

Next, the lateral space-time surfaces of $C_i^n$ are parametrized using a set of bilinear basis functions as
\be
\partial C_{ij}^n = \mathbf{\tilde{x}} \left( \chi,\tau \right) \!=\!
\sum\limits_{k=1}^{4}{\beta_k(\chi,\tau) \, \mathbf{\tilde{X}}_{ij,k}^n },	
\quad 0 \leq \chi \leq 1,  \	0 \leq \tau \leq 1, 										 
\label{eq.SurfParBeta}
\ee
where $\mathbf{\tilde{X}}_{ij,k}^n$ represent the physical space-time coordinates of the four vertices of $\partial C_{ij}^n$, and the $\beta_k(\chi,\tau)$ functions are defined as follows 
\begin{eqnarray}
\beta_1(\chi,\tau) = (1-\chi)(1-\tau), && 
\beta_2(\chi,\tau) = \chi(1-\tau), \nonumber\\
\beta_3(\chi,\tau) = \chi\tau, && 
\beta_4(\chi,\tau) = (1-\chi)\tau.
\label{BetaBaseFunc}
\end{eqnarray}
The mapping in time is given by the transformation 
\be
t = t_n + \tau \, \Delta t, \qquad  \tau = \frac{t - t^n}{\Delta t}, 
\label{eq.timeTransf}
\ee 
hence the Jacobian matrix $J_{\partial C_{ij}^n}$ of the parametrization is 
\be
J_{\partial C_{ij}^n} = \left( \begin{array}{ccc} \vec{e}_r & \vec{e}_\varphi & \vec{e}_t \\[2pt] 
	\frac{\partial r}{\partial \chi} & \frac{\partial \varphi}{\partial \chi} & \frac{\partial t}{\partial \chi} \\[4pt] 
	\frac{\partial r}{\partial \tau } & \frac{\partial \varphi}{\partial \tau } & \frac{\partial t}{\partial \tau } 
\end{array} \right) = \left( \begin{array}{c} \mathbf{\tilde{e}} \\[2pt]  \frac{\partial \mathbf{\tilde{x}}}{\partial \chi} \\[4pt]  \frac{\partial \mathbf{\tilde{x}}}{\partial \tau} \end{array} \right). 
\label{JacSTsurf}
\ee
The space-time unit normal vector $\mathbf{\tilde n}_{ij}$ can be evaluated computing the normalized cross product between the transformation vectors of the mapping \eqref{eq.SurfParBeta}, i.e.
\be
| \partial C_{ij}^n| = \left| \frac{\partial \mathbf{\tilde{x}}}{\partial \chi} \times \frac{\partial \mathbf{\tilde{x}}}{\partial \tau} \right|, 
\quad 
\mathbf{\tilde n}_{ij} = \left( \frac{\partial \mathbf{\tilde{x}}}{\partial \chi} \times \frac{\partial \mathbf{\tilde{x}}}{\partial \tau}\right) / | \partial C_{ij}^n|,
\label{n_lateral}
\ee
where $| \partial C_{ij}^n|$ is the determinant of the Jacobian matrix $J_{\partial C_{ij}^n}$ and represents also the area of the lateral surfaces.
Moreover, exploiting the parametrization in \eqref{eq.SurfParBeta}-\eqref{eq.timeTransf} and choosing 
$\chi = 0.5$ and $\tau =0.5$ we recover the coordinates $M_{i,j}^n$ of the space-time midpoint of the lateral surfaces.

\subsection{Well balanced direct ALE scheme}
\label{sec.WBALE2d}

In order to obtain a space-time formulation of a direct path-conservative ALE scheme, as 
proposed in \cite{dumbser2013high}, the governing PDE \eqref{eq.generalform} is first reformulated in a space-time divergence form as
\be
\tilde \nabla \cdot \tilde{\F}(\Q) + \tilde{\B}(\Q) \cdot \tilde{\nabla} \Q = \S(\Q), \qquad 
\tilde \nabla  = \left( \partial_r, \, \partial_\varphi, \, \partial_t \right)^T  
\label{PDEdiv}
\ee 
with
\begin{equation*} 
\tilde{\F}  = \left( \mathbf{F}, \, \Q \right)^T\! = \left( \mathbf{f}, \, \mathbf{g}, \, \Q \right)^T\!, \  \tilde{\B}  = \left (\B, \textbf{0} \right )^T\! = \left (\B_1, \textbf{0}, \textbf{0} \right)^T\!, \ \text{and } \S = \0, 
\end{equation*} 
and it is then integrated over the space-time control volume $C_i^n$ 
\be
\label{STPDE}
\int_{C_i^n} \left ( \tilde \nabla \cdot \tilde{\F}(\Q) + \tilde{\B}(\Q) \cdot \tilde{\nabla} \Q \right ) \, d\mathbf{x} dt= \0 \, .
\ee 
Now, the space-time volume integral of $\tilde \nabla \cdot \tilde{\F}(\Q)$  can be rewritten using the Gauss theorem as
\be
\label{eq.I1}
\int_{\partial C^n_i}    \tilde{\F} \cdot \mathbf{\tilde n} + 
\int_{C_i^n}    \tilde{\B} \cdot \tilde{\nabla} \Q = 
\0 , 
\ee
where $\mathbf{\tilde n} = (\tilde n_r,\tilde n_\varphi,\tilde n_t)$ is the outward pointing space-time unit normal vector on the space-time surface $\partial C^n_i$.
	
Taking into account the jump of $\tilde{\B}$ at the interfaces, the final high order ALE one-step finite volume scheme is then obtained from equation \eqref{eq.I1} as
\be
\label{eq.Scheme}
|T_i^{n+1}| \, \Q_i^{n+1} =  |T_i^n| \, \Q_i^n 
& - \sum \limits_{j} \,\,\int_0^1 \!\! \int_0^1 
| \partial C_{ij}^n| \ \tilde{\mathbf{D}}_{ij}  \cdot \mathbf{\tilde n}_{ij} d\chi d\tau \!\!\!\!\!\!\!\!\!\!\!\! \!\!\!\!\\
& - \int_{C_i^n}    \tilde{\B}(\q_i^n) \cdot \tilde{\nabla} \q_i^n
\ d\mathbf{x} dt \\
\ee 
where $\q_i^n(\x,t)$ is a well balanced second order reconstruction of the conserved variables $\Q$ inside cell $T_i$ at time $t^n$, and 
the discontinuity of the solution at the space-time sub--face $\partial C_{ij}^n$ is resolved by a well balanced path-conservative ALE flux $ \tilde{\mathbf{D}}_{ij} \cdot \mathbf{\tilde n}_{ij}$, which accounts for the jump in the discrete solution between two neighbors across the intermediate space-time lateral surface.  

In particular when the lateral surface is shared between more than two control volumes we have to compute the flux across each sub-piece and sum each contribution (see \cite{gaburro2016direct} for further details).

\subsubsection{Well balanced ALE numerical flux function}
\label{ssec.2dFirstOrderMethod}

The core of the well balanced method in \eqref{eq.Scheme} is the design of the well balanced space-time flux function. 
Its final expression will be 
\be 
\label{eq.WBspacetimeFlux_2}
 {\tilde{\mathbf{D}}_{ij}} \cdot \mathbf{\tilde n}_{ij}  =  \quad 
 & \frac{1}{2} \left  ( {\tilde{\F}(\q^+) + \tilde{\F}(\q^-)} + { \mathcal{B}_{ij} \left(  \q^+ - \q^- \right)  } \right ) \cdot \mathbf{\tilde n}_{ij} \!\!\!\!\!\!\\
  	- &\frac{1}{2} \Vn_{ij} \left (\q^+-\q^{-} \right ),
\ee
where $\q^-$ is the value of the reconstructed numerical solution inside the element $C_i^n$ evaluated at the space-time midpoint $M_{i,j}^n$ of the lateral surface $\partial C_{ij}^n$, and $\q^+$ is the evaluation at the same point of the reconstructed numerical solution inside the neighbor $C_j^n$ at $\partial C_{ij}^n$.
Besides, generalizing the notation introduced in Section \ref{sec.NumMethod1d},  $\tilde \F$ is the physical flux, the term ${ \mathcal{B}_{ij} \left(  \q^+ - \q^- \right)}$ represents a well balanced way to write the non-conservative products, and  $\Vn_{i+\frac{1}{2}} \left (\q^+-\q^- \right )$ is the viscosity term.

As already pointed out, according to \cite{Pares2006}, the numerical flux should satisfy the following properties
\be
\label{eq.dcond1}
{\tilde{\mathbf{D}}_{ij}}(\Q, \Q) \cdot \mathbf{\tilde n}_{ij} = \0 \quad \forall \Q \in \Omega, \ \text{and} 
\ee
\be
\label{eq.dcond2}
{\tilde{\mathbf{D}}_{ij}}(\q^-, \q^+) \cdot \mathbf{\tilde n}_{ij} =   \int_0^1 \mathbf{A}^{\!\! \mathbf{V}}_{\mathbf{n}}\left(\Phi(s;\q^-, \q^+)\right)\de{\Phi}{s}(s;\q^-,  \q^+)ds,
\ee
where,  due to the ALE framework, 
\be 
\label{eq.ExtendedJacobian2}
& \mathbf{A}^{\!\! \mathbf{V}}_{\mathbf{n}}(\Q) = \sqrt{\tilde n_r^2 \!+\! \tilde n_\varphi^2}  \left ( \!\left (   \de{\mathbf{{F}}}{\Q} + {\B}  \right ) \! \cdot \! \mathbf{{n}}  \! -\! 	\left(\mathbf{V} \!\cdot\! \mathbf{n}\right) \mbf{I}\, \right ), \\   
&\mathbf{n} = (n_r, n_\phi) = \frac{(\tilde n_r, \tilde n_\varphi)^T}{\sqrt{\tilde n_r^2 + \tilde n_\varphi^2}},  \!\!
\ee 
with $\mathbf{I}$ representing the identity matrix and $\mathbf{V} \cdot \mathbf{n}$ denoting the local normal mesh velocity.  

We explain now how to discretize $\mathcal{B}_{ij}$ and $ \mathcal{V}_{ij}$ in \eqref{eq.WBspacetimeFlux_2} in a well balanced way.
Here we perform our reasoning edge--by--edge and we distinguish two situations: the first one across the vertical edges, which evolving in time originate a surface orthogonal to the radial direction, easier to be treated, and the second one across the other two parallel edges (see the constraints stated at the beginning of Section \ref{ssec.MovingDomainDiscretization})).
For the sake of clarity, in Appendix \ref{app.WellBalancing_1element} we present the proof that our scheme is well balanced taking into account a single element.

First of all, it is easy to see that the flux across the sub--surfaces evolved from the vertical edges coincides with the one dimensional flux. 
Indeed, in this case,  $\mathbf{n} = (n_r, 0)$,  $\mbf{V} = (0, v)$ and so $\mbf{V}\cdot \mbf{n} = 0$.
Hence $ \mathbf{A}^{\!\! \mathbf{V}}_{\mathbf{n}}(\Q) = \mbf{J_f}(\Q) + \B_1(\Q)  $ which coincides with \eqref{eq.ExtendedJacobian}. 
So we can discretize $\mathcal{B}_{ij}$ as stated in \eqref{eq-B}-\eqref{b2}-\eqref{eq-b3}-\eqref{eq-b4}, and $\mathcal{V}_{ij}$ by using the Osher-Romberg method \eqref{eq.ViscosityOR} or the modified HLL scheme as described in Section \ref{ssec.hll1d}.
Therefore the scheme is well balanced in the radial direction and second order accurate provided that the reconstruction $\q_i^n$ and the integrals in \eqref{eq.Scheme} are computed in a well balanced manner and with second order of accuracy (see Section \ref{ssec.2d}).

For what concerns the flux through the other two surfaces (see Point (\textit{iii}) of Section \ref{ssec.MovingDomainDiscretization}) let us first state the following remark.
\begin{remark}
\label{remark.geometry}
Given an element $T_i^n$ consider its two edges which are parallel between them but not vertical. Their evolution in time originates two parallel surfaces with the same areas and opposite normal vectors. Moreover call $T^n_{j_1}$ and $T^n_{j_2}$ the two neighbors of $T^n_i$ through these edges. Since the barycenters of $T^n_i$, $T^n_{j_1}$ and $T^n_{j_2}$ are aligned on the same vertical line, i.e. their $r$-coordinate is the same, the equilibrium values $\Q_i^E$, $\Q_{j_1}^E$ and $\Q_{j_2}^E$ coincide.
\end{remark}

Now let us rewrite \eqref{eq.ExtendedJacobian2} as
\be 
\label{eq.ExtendedJacobian3}
& \mathbf{A}^{\!\! \mathbf{V}}_{\mathbf{n}}(\Q) =  \sqrt{\tilde n_r^2 \!+\! \tilde n_\varphi^2}  \left ( \left (  \mbf{J}_{\f} + \B_1 \right )n_r + \mbf{J}_{\g}  n_\varphi - \left( \mbf{V} \!\cdot\!\mbf{n} \right ) \mbf{I}\, \right ).
\ee 
and \eqref{eq.dcond2} as
\be
\label{eq.dcond3}
{\tilde{\mathbf{D}}_{ij}}(\q^-, \q^+) \!\cdot\! \mathbf{\tilde n}_{ij} \! =  \sqrt{\tilde n_r^2 \!+\! \tilde n_\varphi^2}  \! \int_0^1  \Bigl (  &  \left (  \mbf{J}_{\f} + \B_1 \right )n_r  + \mbf{J}_{\g}  n_\varphi \\ & - \left( \mbf{V} \!\cdot\!\mbf{n} \right) \mbf{I} \  \Bigr) \,  \de{\Phi}{s}(s)ds.
\ee 
Thus, by exploiting the linearity of the integral, we can give the discretization of $ \tilde{\mathbf{D}}_{ij} \cdot \mathbf{\tilde n}_{ij} $ in \eqref{eq.WBspacetimeFlux_2} as the sum of the following contributions
\be 
\label{eq.WBspacetimeFlux_3}
{\tilde{\mathbf{D}}_{ij}} \cdot \mathbf{\tilde n}_{ij}  =  \quad
& \frac{1}{2} \left  ( {\f(\q^+) + \f(\q^-)} + { \mathcal{B}_{ij} \left(  \q^+ - \q^- \right)  } \right ) \tilde{n}_r \\
+ & \frac{1}{2} \left  ( {\g(\q^+) + \g(\q^-)} \right ) \tilde{n}_\varphi \\
+ &\frac{1}{2} \left (\q^+ + \q^-\right )\, \tilde{n}_t - \frac{1}{2} \Vn_{ij} \left (\q^+-\q^{-} \right ).
\ee
Note that, whereas the discretization of $\tilde{\F}$ and of $\mathcal{B}_{ij}$ can be splitted, the same cannot be done automatically for the viscosity $\mathcal{V}_{ij}$, whose expression depends on the chosen method (Osher-Romberg, HLL or others).

The expression in \eqref{eq.WBspacetimeFlux_3} results to be well balanced, provided that a well balanced expression for $\mathcal{V}_{ij}$ is given. 
Indeed the first row coincides with the one dimensional flux along the radial direction for which $\mathcal{B}_{ij}$ is given by \eqref{eq-B}-\eqref{b2}-\eqref{eq-b3}-\eqref{eq-b4} that are well balanced.
With regards to the second line we know that in general it is not zero evaluated at the equilibrium because, as already pointed out at the beginning of the section, $\g$ is not zero evaluated at the equilibrium. 
But, if we consider, together with the flux between $T^n_i$ and $T^n_{j_1}$, also the flux between $T^n_1$ and $T^n_{j_2}$ and we sum them up, we can see that all the values at the equilibrium cancel exactly, thanks to the properties stated in Remark \ref{remark.geometry}, that follows from the geometrical constraints we have imposed on our mesh.
Finally, the same argument shows that also the third line goes to zero when $\q^- = \Q_i^E$ and $\q^+ = \Q_{j_1,j_2}^E$.


\noindent \textbf{Viscosity term}


\noindent To end with the formulation of the well balanced ALE flux \eqref{eq.WBspacetimeFlux_2} across this second kind of surfaces, we have to provide an expression for the viscosity $\mathcal{V}_{ij}(\q^+-\q^-)$ which vanishes on stationary solutions \eqref{eq.EquilibriaConstraint}-\eqref{eq.PressureAndGravForces}.

First of all, it is easy to generalize the Osher-Romberg scheme introduced in Section \ref{ssec.OsherRomberg}. 
Indeed in the two dimensional ALE context the viscosity matrix introduced in \eqref{eq.OsherMatrixFormal} can be written as
\be
\Vn_{ij}(\q^+-\q^-) = \int_0^1  \left |\, \mathbf{A}^{\!\! \mathbf{V}}_{\mathbf{n}}(\Q)  \left ( \Phi(s)  \right) \,  \right | \partial_s \Phi(s), \quad   0 \le s \le 1. 
\ee 
Following the same reasoning of Section \ref{ssec.OsherRomberg} we get the following expression
\be
\label{viscosity_n_osher2}
\Vn_{ij}(\q^+-\q^-) =\sum_{j=1}^l \omega_j \text{sign}\left( \mathbf{A}^{\!\! \mathbf{V}}_{\mathbf{n}}(\Phi(s_j)\right) \frac{\mathcal{R}_j}{2\epsilon_j},
\ee
where 
\be
\label{eq.OsherRomberg2d_Rterm}
\mathcal{R}_j= \tilde{\F}(\Phi(s_j+\epsilon_j))-\tilde{\F}(\Phi(s_j-\epsilon_j))  \\
 + \tilde{\Bn}_{\Phi_j}\left(\Phi(s_j+\epsilon_j)-\Phi(s_j-\epsilon_j)\right)
\ee
is discretized as explained in the 1D case above and the Romberg quadrature formula with $l=3$ is still used. 
Hence, if $\q_{i}^n$ and $\q_{i+1}^n$ lie on the same stationary solution $\Phi(s)=\Phi^E(s)$ and $\mathcal{R}_j=\0$, $j=1, \dots, l$.

Thus, the extension to two dimensions of the Osher-Romberg scheme results to be straightforward. The only drawback is that the complete eigenstructure of the extended Jacobian matrix  $\mathbf{A}^{\!\! \mathbf{V}}_{\mathbf{n}}$ should be computed, which could be costly in particular when edges are not parallel to the axis (we underline that $\mathbf{A}^{\!\! \mathbf{V}}_{\mathbf{n}}$ does not enjoy the property of rotational invariance that characterizes the Euler equations in Cartesian coordinates).
As counter part, the method is very little dissipative and allows us to obtain very good results in convective transport problems. 

The generalization of the HLL scheme is simpler.
Equation \eqref{eq.hll_matrixform} can be rewritten in two dimensions as
\be
\Vn_{ij}(\q^+-\q^-)= \alpha^0_{ij} \mathbf{I}_{ij} (\q^+-\q^-) + \alpha^1_{ij} \mathcal{R}_{ij},
\ee
where $\mathbf{I}_{ij}$ is the identity matrix, 
\be 
\mathcal{R}_{ij}=\F(\q_{i+1})-\F(\q_i)+\Bn_{ij}(\q^+-\q^-)
\ee
(which can be discretized as described in Section \ref{ssec.2dFirstOrderMethod} to maintain the well balanced properties), and 
$\alpha^{0,1}_{ij}$ can be computed as in \eqref{eq.alphaHll} being $S^L$ and $S^R$ the minimum and the maximum eigenvalues of $\mathbf{A}^{\!\! \mathbf{V}}_{\mathbf{n}}(\q^n_{i,i+1})$.

For the same reasons stated in Section \ref{ssec.hll1d}, $\mathbf{I}_{ij}$ must be replaced by a matrix that vanishes when a stationary solution is considered.
In particular we choose the following identity modification
\be
\tilde{\mathbf{I}}_{ij} = \tilde{\mathbf{I}}_{i+1/2}\, n_r + \mathbf{I} n_\varphi,
\ee
where $\tilde{\mathbf{I}}_{i+1/2}$ is given by \eqref{viscosity_HLL}, which we already know to be well balanced for stationary solutions. 
Moreover it follows from Remark \ref{remark.geometry} that when $n_\varphi \ne 0$ the term $\mathbf{I} n_\varphi$ cancels at the equilibrium (by considering the two contributions of the neighbors $T^n_{j_1}$ and $T^n_{j_2}$ of $T^n_i$).

\subsubsection{2nd order well balanced reconstruction}
\label{ssec.2d}
The missing ingredient for \eqref{eq.Scheme} to be well balanced up to second order is the definition of a second order well balanced reconstruction operator.
As in the one dimensional case we are going to employ a combination of a smooth stationary solution together with the standard MUSCL method, hence our reconstruction will be of the form
\be
\label{eq.FinalReconstruction2}
\q_i^n(\x,t) = \Q^E_i(\x,t) + \mathcal{P}^{f}_i(\x,t), \quad \x \in C_i^n,
\ee 
where, as in Section \ref{ssec.1d_2ndorder}, $\mathcal{P}^{f}_i(\x,t)$ is the standard MUSCL method applied in order to reconstruct the fluctuations with respect to the given stationary solution computed for all the neighbors $T_j^n$ of $T_i^n$ as 
\be
\label{eq.FluctuationForRec2}
\Q^{f\!,n}_j \!= \Q_j^{n} \!- \Q^{E\!,n}_j.
\ee

The expression of the reconstruction operator is
\be
\label{reconstruction_fluctuation2}
\mathcal{P}^{f,n}_i (\x,t) = \Q_i^{f,n} + \Phi_i \nabla \Q_i^{f,n} (\x-\x_i) + \partial_t\Q_i^n (t -t^n),
\ee 
where $\x_i$ is the barycenter of cell $T_i^n$.

To compute $\nabla Q_i^{f\!,n}$ we use the standard MUSCL method (see \cite{leer5})  together with the Barth and Jespersen limiter (see \cite{BarthJespersen}). We would like to remark that the employed methods are standard, the novelty is in the fact that both are applied only to the fluctuations.
  
  
Finally, the term $\partial_t \Q_i^n$ indicates the time derivative of $\Q$ and  it can be computed using a discrete version of the governing equation
\be
\partial_t \Q_i^n = \left( \mbf{J}_{\f} + \B_1 \right )\!|_{\x_i} \partial_r \Q (\x_i) + \mbf{J}_\g |_{\x_i} \, \partial_\varphi \Q(\x_i),
\ee
evaluated at the barycenter $\x_i$ of $T_i^n$. 
In particular the gradient of the conserved variables must be expressed as the gradient of the equilibrium plus the previously computed gradient of the fluctuation, i.e. 
\be
\label{eq.gradientEq+fluct}
\nabla Q = \nabla \Q^E_i + \nabla \Q^f_i = \nabla \Q^E_i + \nabla \Q_i^{f,n}, 
\ee 
in order to preserve the well balancing.
The same idea of \eqref{eq.gradientEq+fluct} can be exploited in order to rewrite
\be
\int_{C_i^n}    \tilde{\B}(\q_i^n) \cdot \tilde{\nabla} \q_i^n
\ d\mathbf{x} dt,
\ee
where, as in Section \ref{ssec.1d_2ndorder}, the equilibrium terms cancel and the remaining terms all contain fluctuations.
So the integral can be computed through the mid-point quadrature rule which is second order accurate on the fluctuations without affecting the equilibrium.

\section{Numerical results in two dimensions}
\label{sec.Results2d}

\subsection{Equilibrium preservation }
\label{sec.2dEquilibriumPreservation}

\begin{figure}
	\begin{center}
	\includegraphics[width=0.7\linewidth]{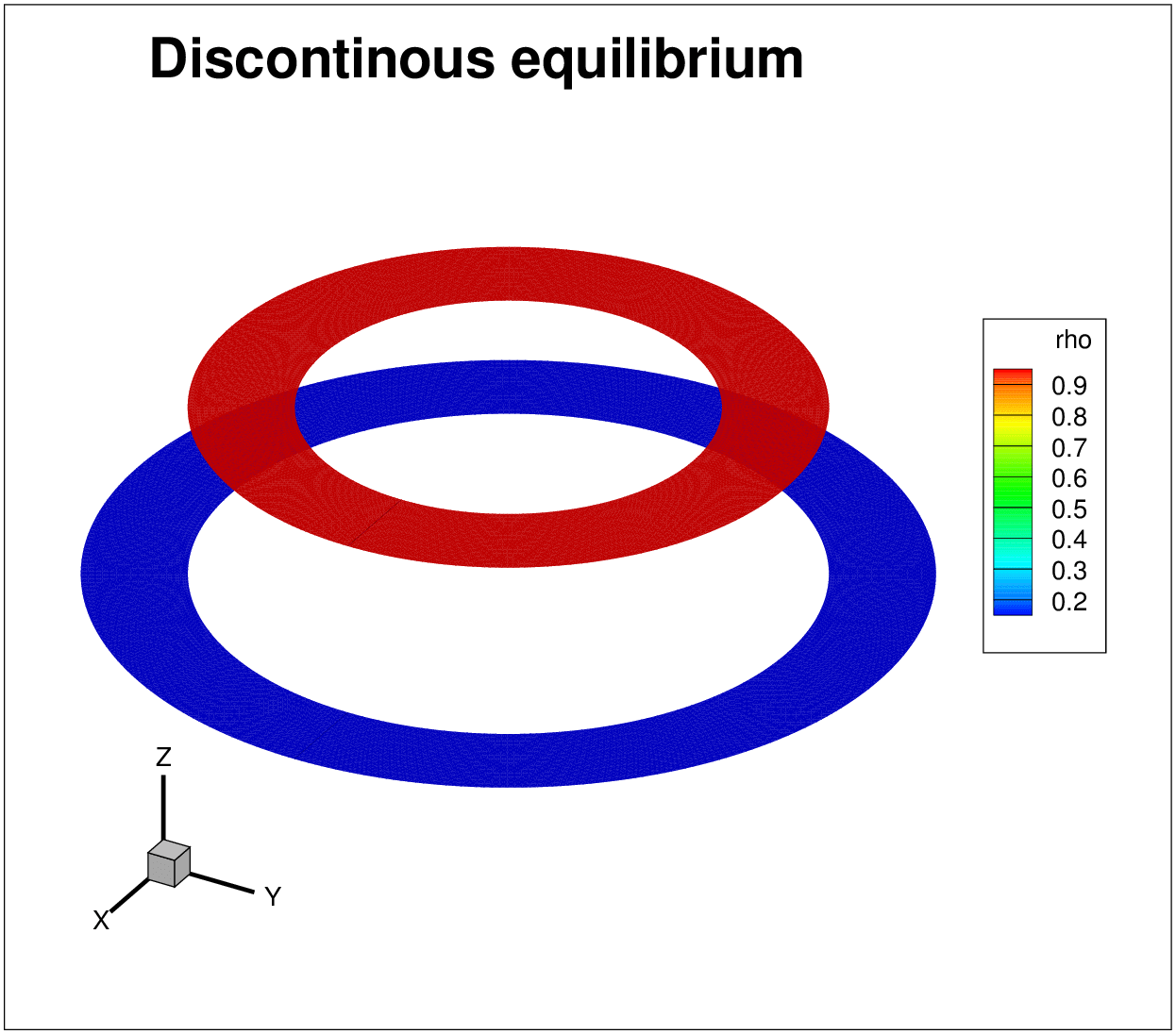}
	\end{center}	
	\caption{Discontinous density profile for the equilibrium solution considered in the test case of Section \ref{sec.2dEquilibriumPreservation}. }
	\label{fig.Equilibrium preservation 2d}	
\end{figure}

First of all we want to show the accuracy of our scheme in preserving some equilibrium of interest. We consider a discontinuous equilibrium 
\begin{equation} 
\label{eq.equilibiurm2d_rhoconst}
\begin{aligned}   
&	\rho = 1, \  \text{if } r < r_m, \quad \rho = 0.1, \  \text{if } r \ge r_m,   \\ 
&	u = 0, \quad 	v = \sqrt{ \frac{G m_s}{r} }, \quad 	P = 1,  
\end{aligned} 
\end{equation} 
with $r_m = 1.5$, $G =1$, $m_s = 1$, over the computational domain $ [r, \varphi] \in [1,2] \times [0, 2\pi]$.
In Figure~\ref{fig.Equilibrium preservation 2d} we depict the density profile at the equilibrium and in Table~\ref{tab.EquilibriumPreservation_2} we report the maximum error, committed using the HLL flux, with respect to the exact solution after long computational times over a coarse mesh, both for order $1$ and $2$. The equilibrium results to be perfectly preserved.

{Then we consider a hydrostatic equilibrium without tangential velocity, so that the gravity force is perfectly balanced with the pressure gradient. The initial data reads 
\be
& \rho = 1, \quad \u=(u,v) = \0, \quad P = 1/r, \quad G = m_s = 1.
\ee 	
We consider a computational domain $ [r, \varphi] \in [1,2] \times  [0, 2\pi]$ covered by a coarse mesh of $20\times40$ elements.
In Table \ref{tab.EquilibriumPreservation_v0} we show the error between the analytical solution and our numerical solution obtained with the second order Osher-Romberg scheme. Since the scheme is 
exactly well balanced the errors are maintained at the order of machine precision for very long computational times. Similar results are also achieved with our well balanced HLL-type flux. }

\begin{table}	
	\caption{Maximum error between the exact and the numerical density obtained with the first and the second order well balanced methods (using the HLL flux). We underline that similar results have been obtained using the Osher-Romberg flux and that the same precision is achieved for the velocities.} 
	\begin{center}
	\begin{tabular}{ccc} 
		\hline 
		\multicolumn{3}{c}{ points  $20\times40$}        \\[-2pt] 
		\hline
		\quad time\qquad  &  $\mathcal{O}1$  & $\mathcal{O}2$   \\[-2pt] 
		\hline
		\quad 10    & 7.32E-13  &   4.20E-13   \\
		\quad 40    & 2.83E-12  &   8.18E-12   \\
		\quad 80    & 3.92E-12  &   1.72E-11   \\
		\quad 100   & 2.25E-12  &   1.99E-11   \\[-2pt] 			
		\hline 
	\end{tabular}	
\end{center}
	\label{tab.EquilibriumPreservation_2}
\end{table}

\begin{table}	
	\caption{{Hydrostatic equilibrium. Maximum error in $L_\infty$ norm between the exact solution and the numerical results for density, velocity and pressure at different times. The values refer to the second order Osher-Romberg ALE scheme, but similar results have been obtained at first order and with the HLL-type flux.}}  
	\begin{center}
		\begin{tabular}{c||cccc} 
			\hline  time &  $E_\rho$ & $E_u$ & $E_v$ & $E_P$ \\ 			
			\hline  1  &  7.77E-15 & 3.29E-16 & 3.95E-16  & 3.33E-16  \\ 
			10  &  1.60E-14 & 3.16E-16 & 1.05E-15  & 3.33E-16  \\ 
			40  &  2.66E-14 & 3.58E-16 & 1.37E-15  & 3.33E-16 \\ 
			80  &  3.02E-13 & 1.30E-13 & 4.98E-14  & 3.87E-14 \\ 
			\hline 
		\end{tabular}	
	\end{center}
	\label{tab.EquilibriumPreservation_v0}
\end{table}

\subsection{Order of convergence}

\begin{figure}
	\begin{center}
		\includegraphics[width=0.9\linewidth]{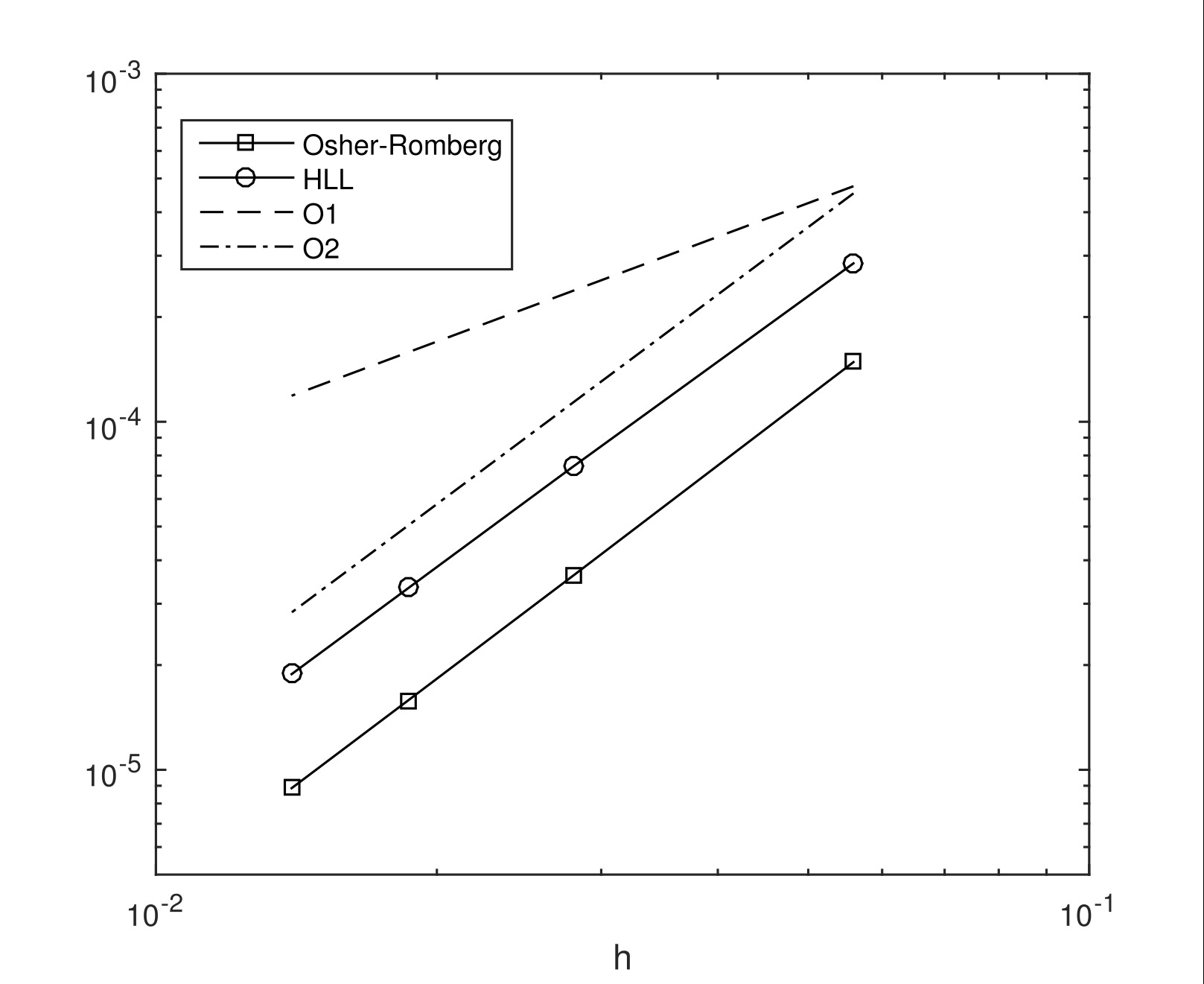}
	\end{center}	
	\caption{{Order of convergence, isentropic vortex for imposed equilibrium (eq.) given by \eqref{eq.equilibiurm2d_rhoconst}, i.e. different from the initial data of the isentropic vortex \eqref{eqn.shuvortex}. We report the $L_1$ error norm of the density obtained with our well-balanced Osher-Romberg and HLL ALE schemes. The dashed lines represent the theoretical slopes of order one and two, respectively.}}
	\label{fig.order_isentropicVortex}	
\end{figure}

\begin{figure}
	\begin{center}
		\includegraphics[width=0.9\linewidth]{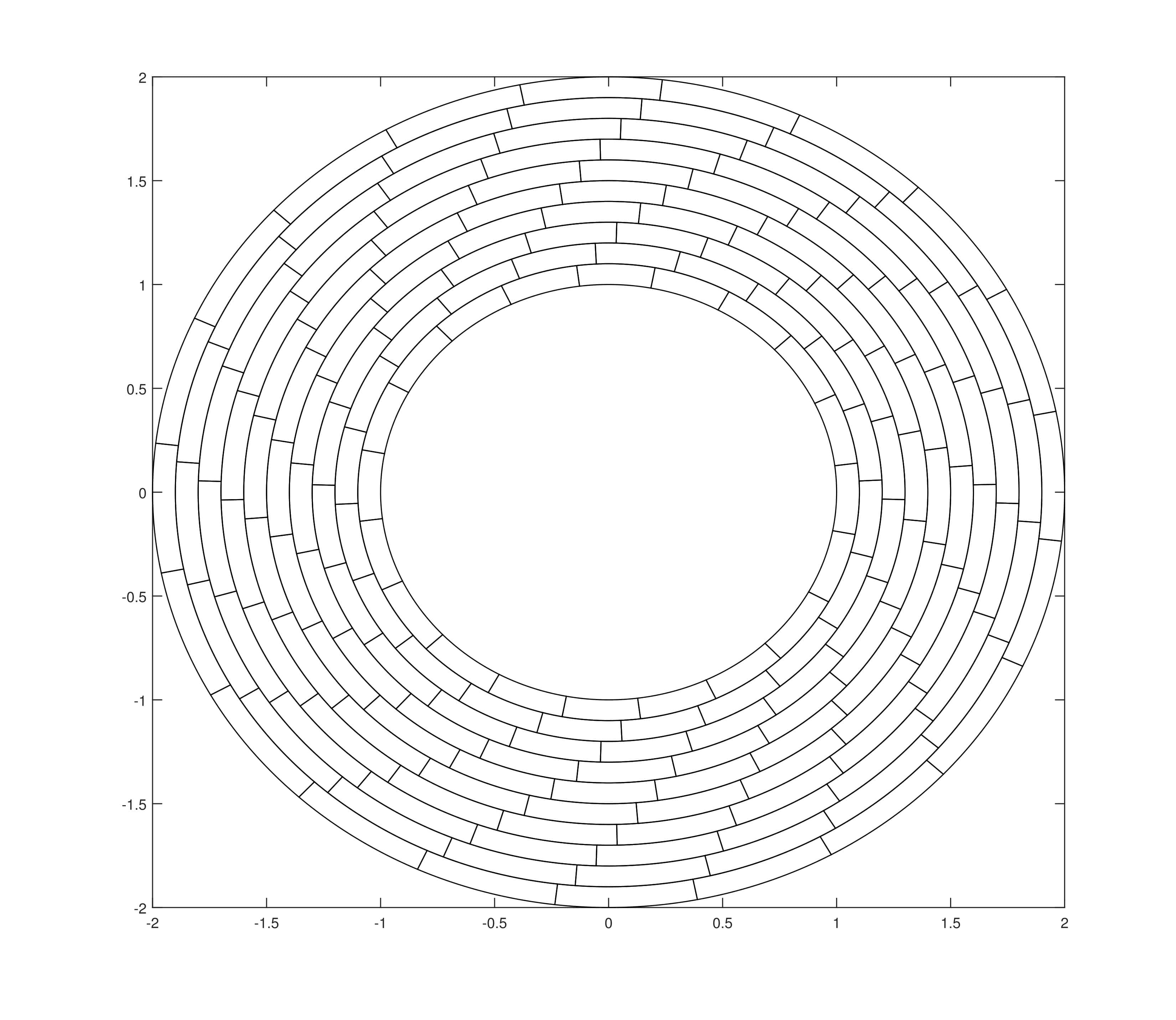} 
	\end{center}	
	\caption{{Isentropic vortex, final mesh. We report the final mesh configuration at time $t_f=1$ obtained with our Osher-Romberg scheme in the case of a very coarse mesh of $10\times20$ elements so that the nonconforming motion is clearly visible.}}
			\label{fig.mesh_isentropicVortex}	
\end{figure}

{To study numerically the order of convergence of our method we consider a smooth isentropic vortex, similar to the one proposed in \cite{hu1999weighted}. 
The initial condition in polar coordinates is given by 
\be 
& 
\rho = 1+ \delta\rho, \quad u = 0, \quad v= \delta v, \quad P = 1 + \delta P, \\ 
& 
\delta v = r \frac{\epsilon}{2\pi} e^{\frac{1-r^2}{2}}, \quad 
\delta T = - \frac{(\gamma-1)^{\epsilon^2}}{8\gamma\pi} e^{1-r^2}, \\
& \delta P = (1+\delta T)^{\frac{1}{\gamma-1}} -1, \quad 
\delta \rho = (1+\delta T)^{\frac{\gamma}{\gamma-1}} -1, 
\label{eqn.shuvortex}
\ee
with $\epsilon=5$, $G=0$, $m_s=0$ and $\gamma =1.4$ and the computational domain defined as $[r,\phi] = [1,2]\times[0,2\pi]$. The final time is $t_f=1$.}
{Our new scheme is able to preserve this equilibrium up to machine precision if we impose the above initial data \eqref{eqn.shuvortex} also as the equilibrium 
profile to be preserved. 
However, it is also possible to impose a different equilibrium profile to be maintained, e.g. the one given by \eqref{eq.equilibiurm2d_rhoconst}. In this way, equilibrium and 
initial condition are not close one to the other so the method comes back to its standard order of convergence, i.e. second order. Refer to Table \ref{tab.order.2d.shu} and 
Figure \ref{fig.order_isentropicVortex} for the numerical results, which confirm that our scheme is indeed second order accurate away from the prescribed equilibrium profile.}
{Finally, we would like to remark that we are working with a moving nonconforming grid. In Figure \ref{fig.mesh_isentropicVortex} we report an example of the 
final mesh configuration obtained with our Osher-Romberg scheme.}
  
\begin{table}	
	\caption{{Order of convergence, isentropic vortex. We report the results obtained with our second order accurate well-balanced Osher-Romberg ALE scheme. The mesh size $h$ is computed as the maximum incircle diameter of the elements of the final mesh. The errors refer to the $L_1$ norm of the difference between our numerical solution and the exact one. 
	The last column refers to the setting where the initial data \eqref{eqn.shuvortex} are also imposed as the smooth known equilibrium to be maintained, hence in this case the scheme is accurate up to machine precision. 
	The other results {are for  the setting where the code is used to evolve a different equilibrium profile  \eqref{eq.equilibiurm2d_rhoconst} that does not coincide with the initial data \eqref{eqn.shuvortex},} so that we can show its formal order of accuracy.}}  
	\begin{center}
		\begin{tabular}{c||cc||c} 
			\hline  mesh size $h$ &  $E_\rho$, eq. \eqref{eq.equilibiurm2d_rhoconst} & $\mathcal{O}(L_1)$ & $E_\rho$, eq. \eqref{eqn.shuvortex}   \\ 			
			\hline  5.59E-2  &  1.48E-4 & - &  1.86E-14    \\ 
			2.80E-2  &  3.60E-5 & 2.04 & 1.45E-13     \\ 
			1.86E-2  &  1.58E-5 & 2.03 & 4.78E-13    \\ 
			1.40E-2  &  8.85E-6 & 2.02 & 5.36E-13    \\ 
			\hline 
		\end{tabular}	
	\end{center}
	\label{tab.order.2d.shu}
\end{table}

\subsection{Riemann problem}
To show the correctness of our method also in the presence of shock waves we solve a classical Riemann problem with non-vanishing angular velocity using both the well balanced HLL and Osher-Romberg ALE schemes. 

\begin{figure*}
		\includegraphics[width=0.45\linewidth]{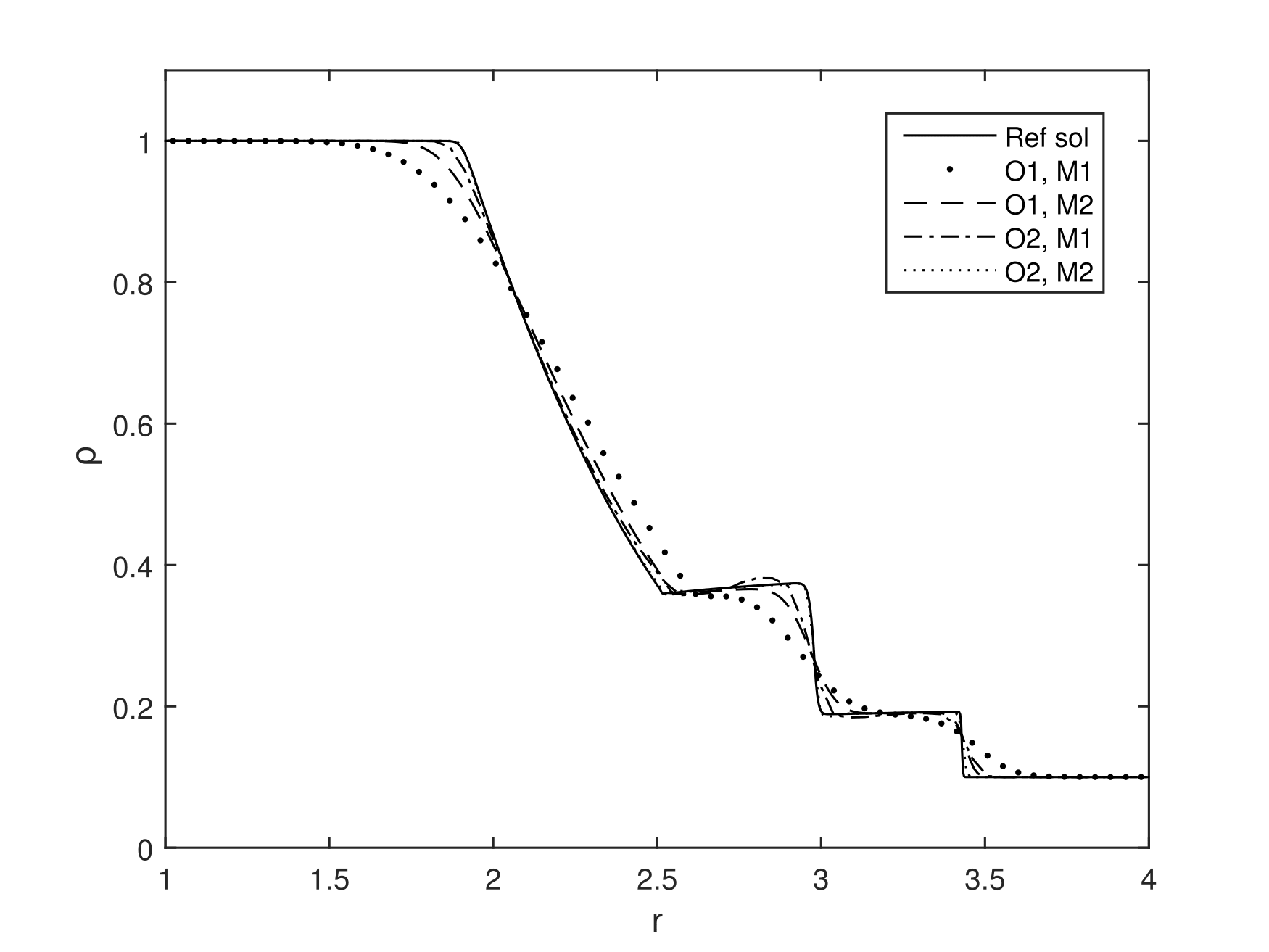}
   	\includegraphics[width=0.45\linewidth]{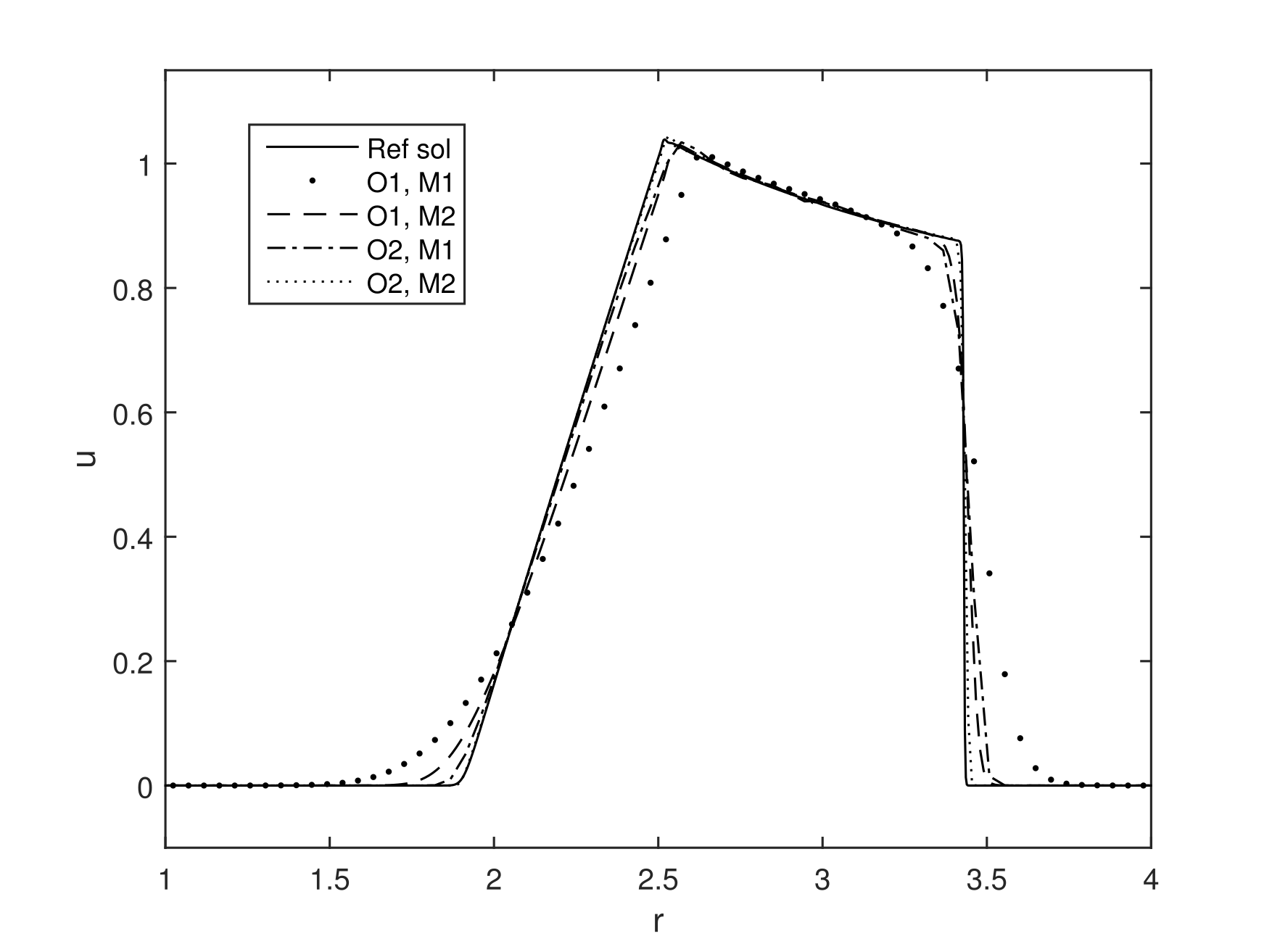}
	\caption{Riemann problem in a 2D domain. The test heave been carried out over two meshes: the first one, $M_1$, with $64 \times 20 $ control volumes and the second one, $M_2$, with $256 \times 40$ control volumes. The reported results have been obtained using the well balanced HLL scheme with first and second order of accuracy. On the left we report the results for the density and on the right for the velocity at the final time $t_f = 0.5$. The graphs have been obtained as a 1D cut along $\varphi = \pi/2$.  One can observe that the second order scheme captures the discontinuities sharply. The results are compared against a reference solution obtained with our second order well balanced HLL scheme in one space dimension with $N=1024$. } 
	\label{fig.2dRiemannProblem}
\end{figure*}

\begin{figure*}
	\includegraphics[width=0.45\linewidth]{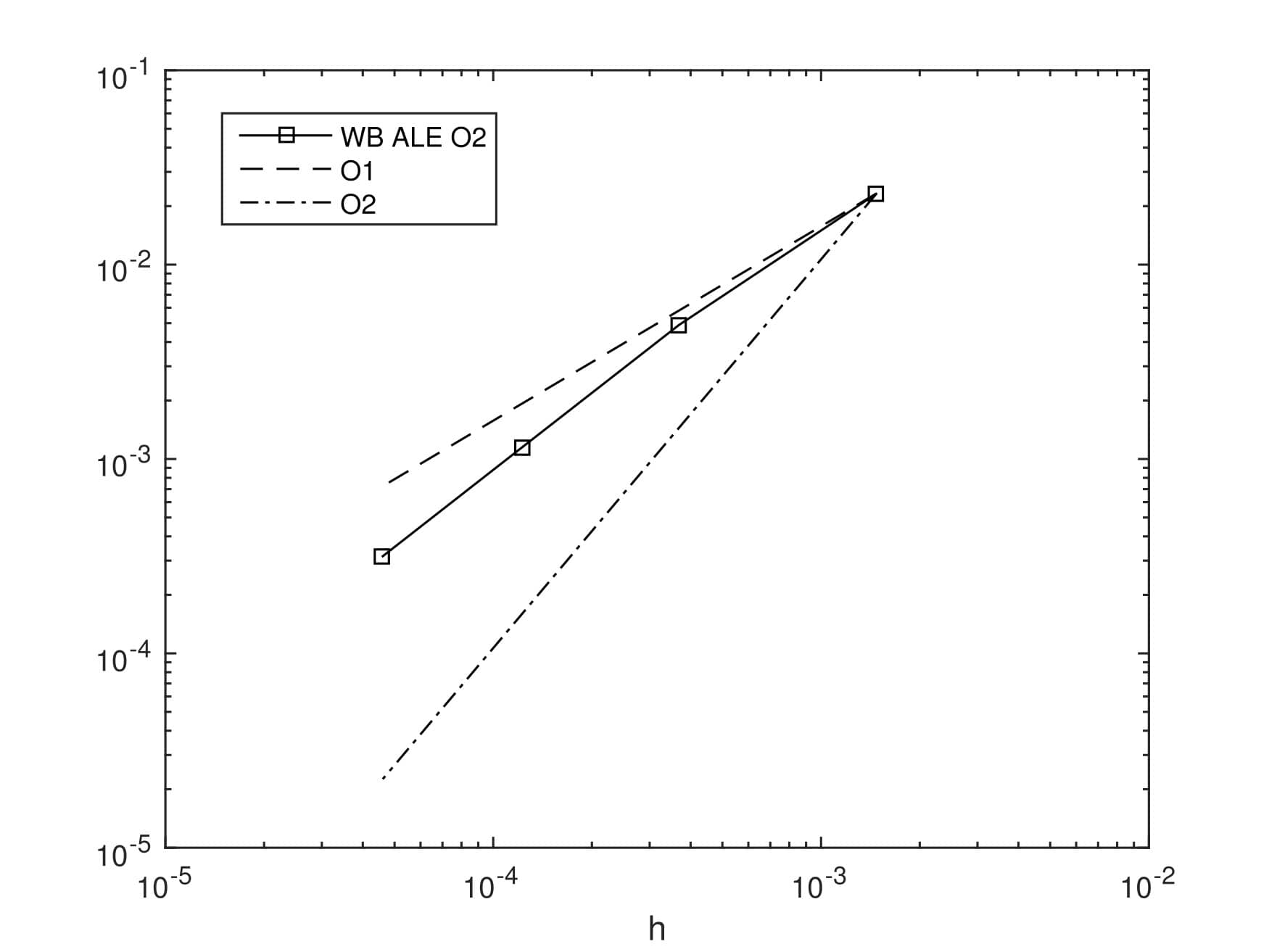}
		\includegraphics[width=0.47\linewidth]{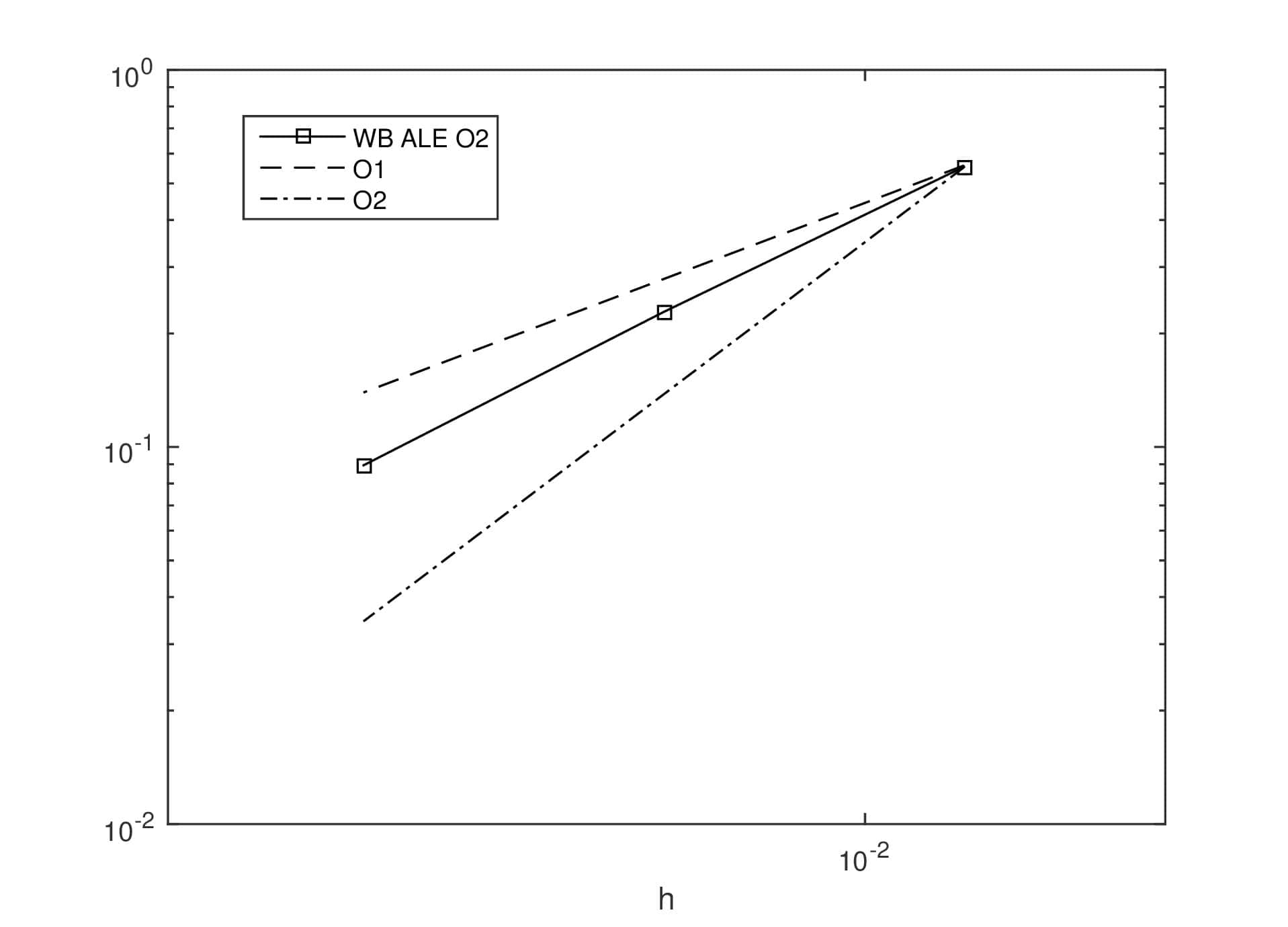}
	\caption{{Convergence results. Left: we refer to the Riemann problem \eqref{eq.RP_2d_init} and compare the results obtained with our WB ALE HLL code with a fine grid reference solution. Right: we refer to the Noh shock test of Section \ref{ssec.Noh}  and we compare our results with the exact solution. Note that the $L_1$ norm of our numerical errors are depicted with squares and is compared with the theoretical slopes of order one and two (dotted lines), respectively. It is evident that the method is better than first order accurate even in presence of shocks.} }
	\label{fig.convergenceResults}
\end{figure*}

We consider the computational domain $[r, \varphi] =  [1,4] \times [0, 2 \pi] $ and we impose the following initial conditions
\be		
\label{eq.RP_2d_init}
&\rho = 1, \  \text{if } r < r_m, \quad \rho = 0.1, \  \text{if } r \ge r_m, \\ 
& u = 0, \quad v = \sqrt{G\, m_s/r}, \\
& P = 1, \  \text{if }  r < r_m, \quad P = 0.1, \  \text{if } r \ge r_m, \\
\ee
with $r_m = 2.5$. The results at the final computational time $t_f = 0.5$ are shown in Figure~\ref{fig.2dRiemannProblem} where we report a cut along $\varphi = \pi/2$ and a comparison
with a one-dimensional reference solution computed on a fine grid using 1024 elements. We note a good agreement between the numerical solution obtained with the well balanced ALE
scheme on moving non-conforming meshes and the reference solution also in this case where the solution is far from any equilibrium. 
{Moreover we show the order of convergence of our method with respect to the reference solution in Figure \ref{fig.convergenceResults}: obviously it cannot reach order two 
because of the presence of shocks. However, the observed convergence order is higher than one.}

\subsection{{Noh shock test}}
\label{ssec.Noh}

{
The Noh shock test consists of a circular infinite strength shock propagating out from the origin. We have chosen this test case to prove that our method can deal also with highly supersonic flows, low
pressure atmospheres and shocks of infinite strength. 
Consider a gas with $\gamma=5/3$ initialized with density $\rho=1$, radial velocity $u = -1$, tangential velocity $v=0$, and  pressure $P=10^{-6}$ as an approximation to zero pressure. The shock wave propagates with speed $1/3$. The exact solution inside the shock region, i.e. $r  \le \frac{\text{t}}{3}$, is given by the following relations
\be
\rho = 16, \quad P=16/3, \quad u=0, \quad v=0, \\
\ee 
and outside the shock region, i.e. $r  > \frac{\text{t}}{3}$, by
\be
\rho = 1 + \frac{t}{r}, \quad P=0, \quad u=-1, \quad v=0. \\
\ee 
We consider an initial domain $ [r, \varphi]  \in [0, 1] \times [0,\pi/2]$.
We impose periodic boundary conditions on $\phi=0=\pi/2$, and we exploit the exact solution to impose the boundary conditions at $r=0$ and on the moving outer boundary. }

\begin{figure}
	\centering
	\includegraphics[width=0.75\linewidth]{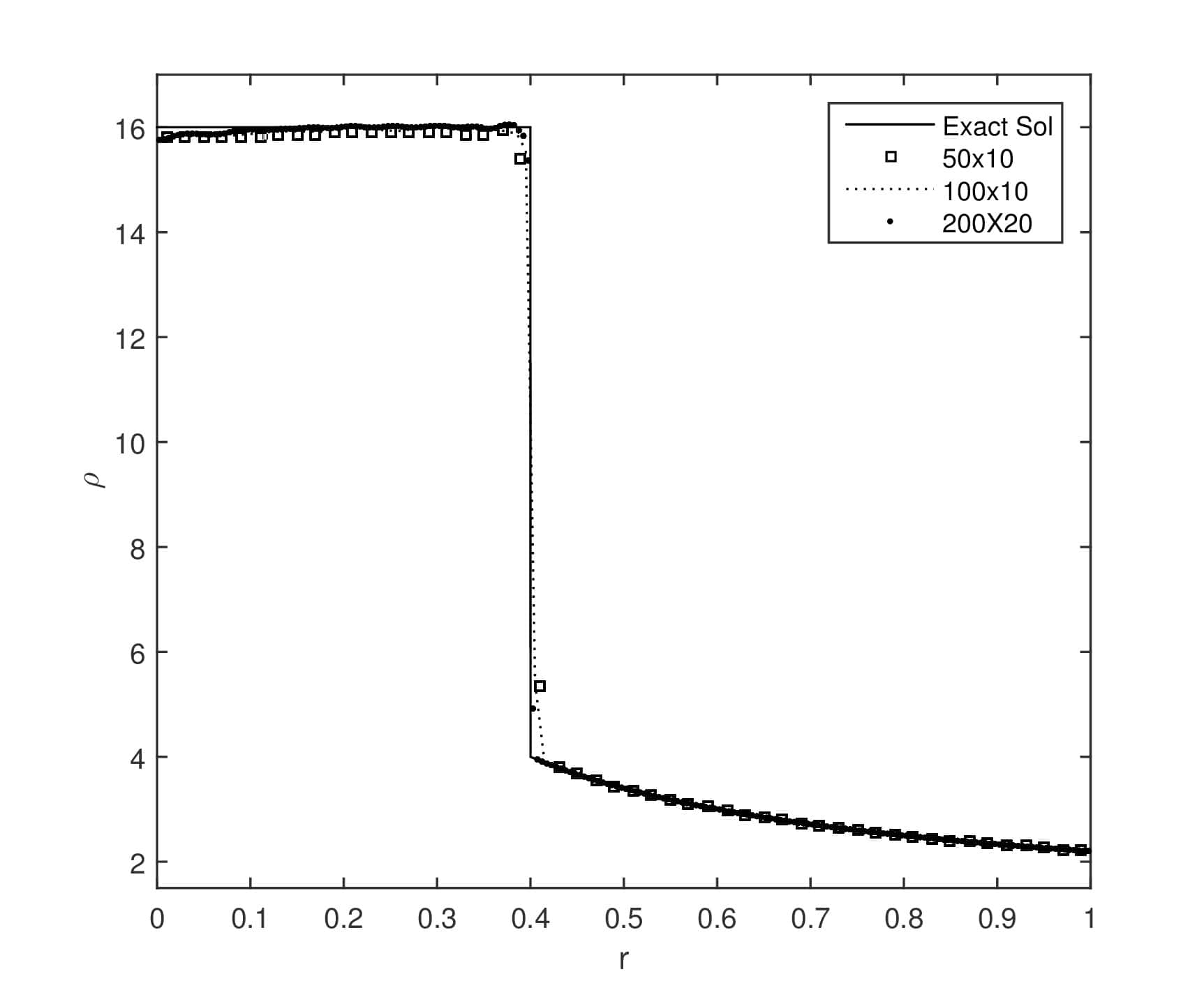}
	\caption{{Noh shock test. We show the numerical results obtained with our second order HLL-type flux at time $t_f = 1.2$ on three fixed grids with respectively $50\times10$, $100\times10$ and $200\times20$ elements. In the figure the density profile $\rho$ has been depicted along the radial direction $r$, compared with the exact solution.}} 
	\label{fig:Noh_confronto_O1O2_tf_1-2}
\end{figure}

\begin{figure*}
	\centering
	\includegraphics[width=0.275\linewidth]{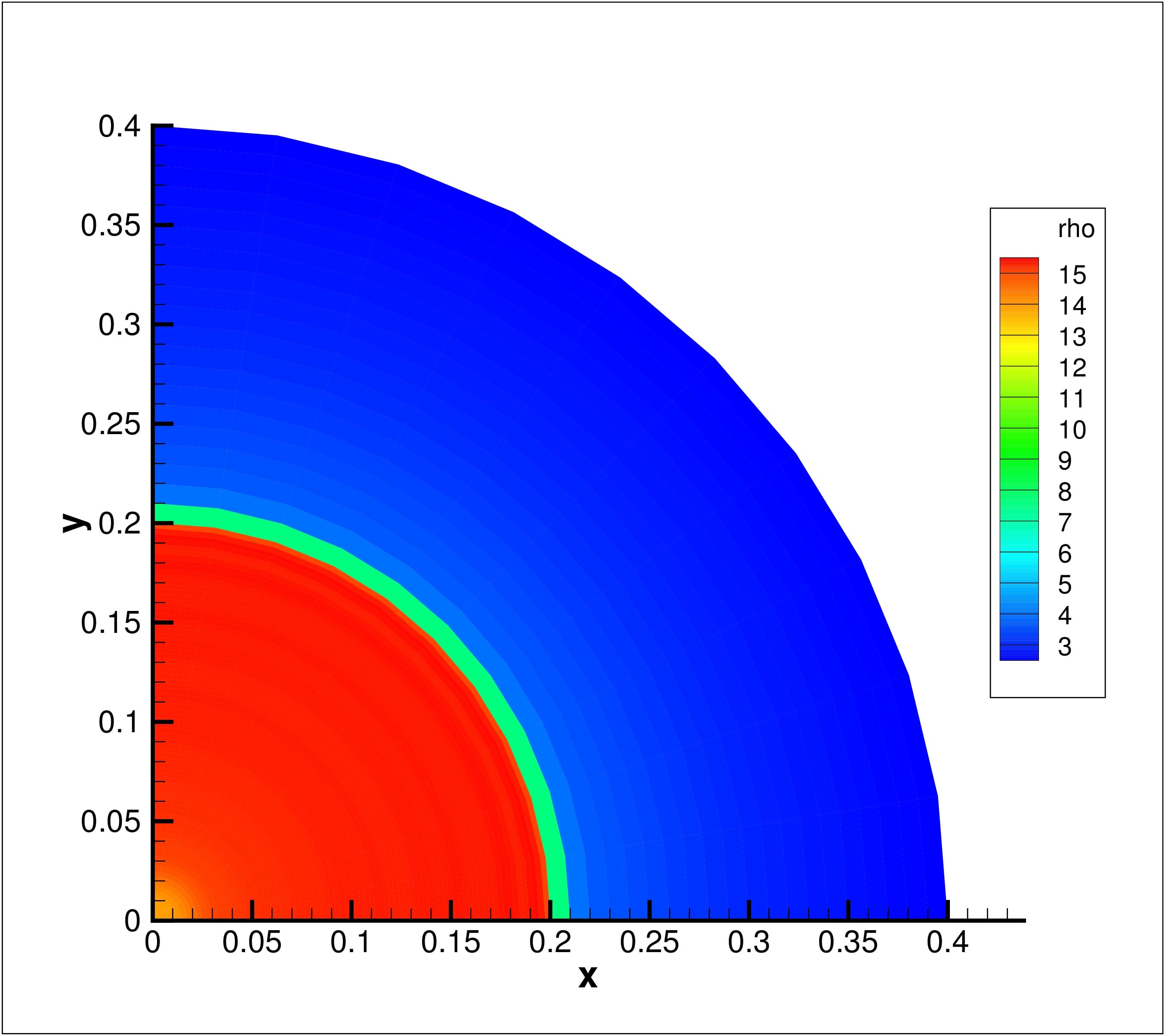} \qquad 
	\includegraphics[width=0.275\linewidth]{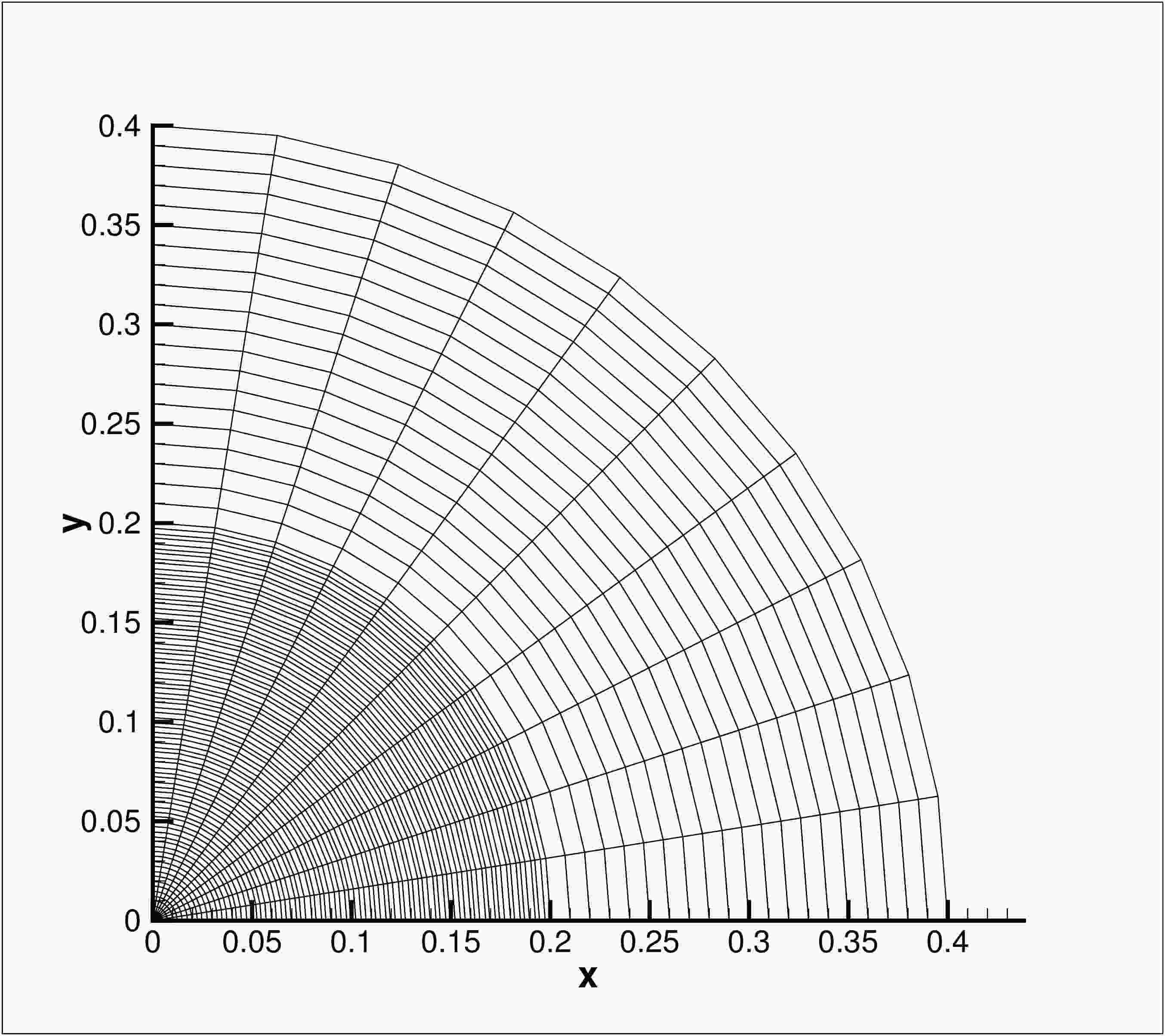} \quad
	\includegraphics[width=0.305\linewidth]{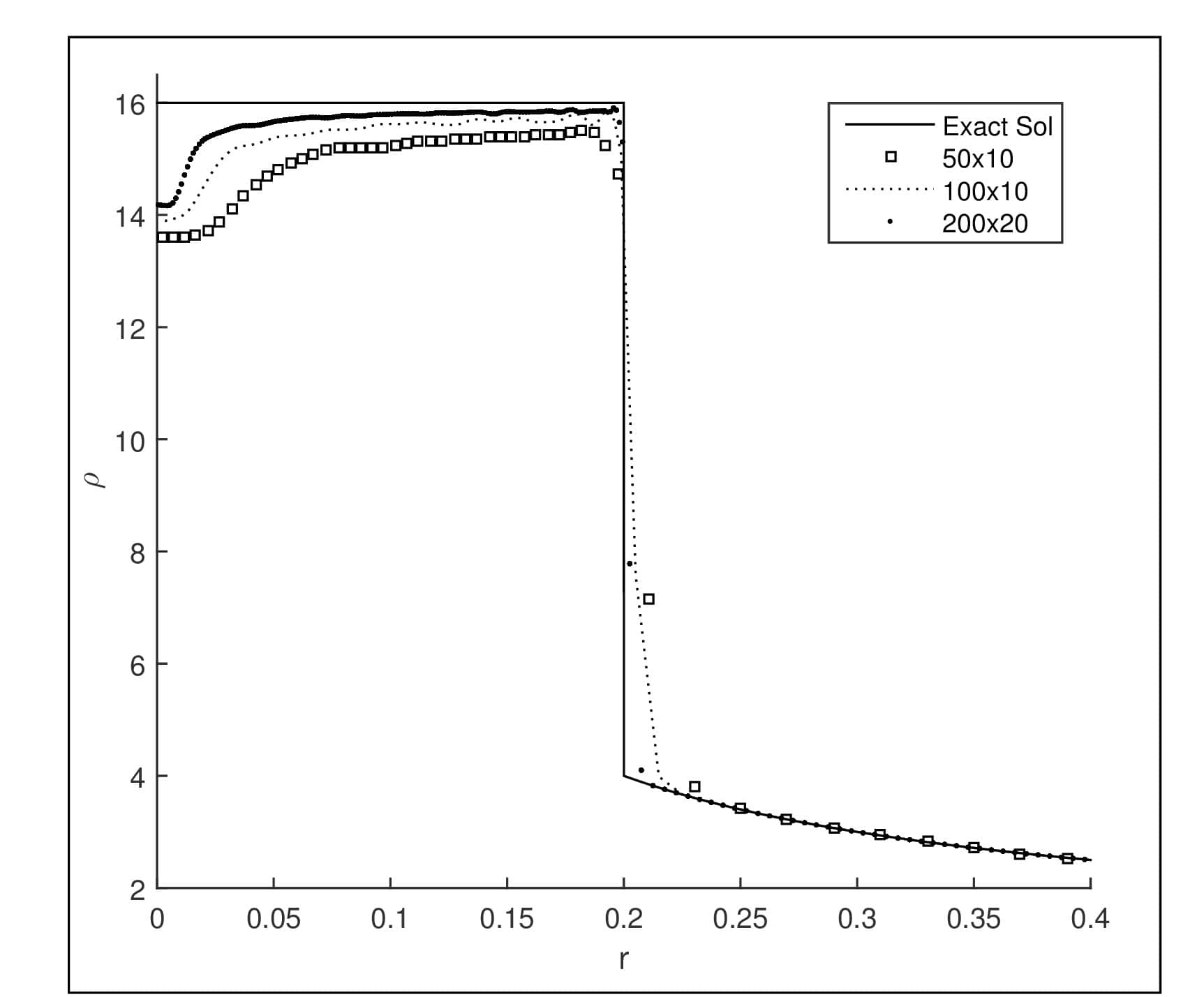}
	\caption{{Noh shock test. We show the density profile (left) and the final mesh (center) obtained with the second order ALE HLL-type scheme at time $t_f=0.6$, using a moving grid of $100\times10$ elements. On the right we compare the density profile along the radial direction $r$ with the exact solution for three different meshes with respectively $50\times10$, $100\times10$ and $200\times20$ elements.}}
	\label{fig:Noh_ALE}
\end{figure*}

{The presented results have been obtained with the HLL-type  scheme. First we have considered the Eulerian case, hence we have imposed a zero mesh velocity. The results at time $t_f=1.2$ obtained with the second order scheme are shown in Figure~\ref{fig:Noh_confronto_O1O2_tf_1-2}.
Then we have employed the ALE framework moving the mesh with the local fluid velocity. Due to the absence of shear flow, the mesh remains conforming. The results obtained with the moving mesh 
are shown in Figure \ref{fig:Noh_ALE}, where the well-known wall heating problem is visible. Apart from the wall heating, in both the cases the method shows a good agreement with the exact solution. } {For what concerns the observed convergence rate of our code in this test we refer to Figure \ref{fig.convergenceResults}.}

\subsection{{Comparison with the PLUTO code}}
\label{ssec.PLUTO_setting}

{For the following test cases that concern Keplerian discs, we compare the results obtained with our new second order well balanced Osher-Romberg scheme with the 
results one can obtain with the PLUTO code. 
PLUTO is a freely-distributed software for the numerical solution of mixed hyperbolic/parabolic systems of partial differential equations (conservation laws) targeting high Mach number flows in astrophysical fluid dynamics.
The code has been systematically checked against several benchmarks available in the literature in the papers \cite{mignone2007pluto} and \cite{mignone2011pluto}, using fixed uniform and AMR grids. 
It provides a multi-physics and multi-algorithm modular environment, where one can choose the Newtonian description for the fluid motion (\texttt{HD} option) and add a potential 
$\Phi = - \frac{Gm}{r} $  to the right hand side by setting the option $\texttt{BODY\_FORCE}$ equal to \texttt{POTENTIAL}. In this way one can study \eqref{eq.EulerCartCons} within this code. 
Then we select \texttt{POLAR} \texttt{GEOMETRY} and we do not activate any other options.} 

{The modular structure allows to choose between different numerical fluxes, limiters, spatial reconstructions and time integrators. In particular, we have selected a little 
dissipative setting by imposing \texttt{LIMITER} equal to \texttt{MC\_LIM} (the monotonized central difference limiter), or sometimes equal to \texttt{MINMOD\_LIM} (the classical minmod limiter), 
and using the \texttt{Roe solver} as numerical flux. 
Then we have compared our second order scheme with both a \textit{second order} configuration of PLUTO (with \texttt{LINEAR} reconstruction in space and \texttt{RK2} in time) and a \textit{third order} configuration (with \texttt{WENO3} reconstruction in space and \texttt{RK3} in time). 
Finally, for the comparison we set the number of elements in PLUTO either equal to the number of elements used for our scheme, or we double it in each dimension.} 

{We remark that within PLUTO special care is taken for the treatment of source and pressure terms when a polar (cylindrical or spherical) geometry is chosen, because in those cases 
the equations are discretized in angular momentum conserving form and pressure terms are treated separately. For this reason the results are more accurate than those obtained with standard 
finite volume techniques.} 


\subsection{Mass transport in a Keplerian disc}  
\label{sec.KeplDisk_transport1}

Let us consider a steady state solution of the Euler equations with gravity which satisfies the constraints in \eqref{eq.EquilibriaConstraint}-\eqref{eq.PressureAndGravForces} and with a constant density profile,
\be
\label{eq.constantEquilibrium}
\rho_E = 1, \quad 
u_E = 0, \quad 
v_E = \sqrt{ \frac{G m_s}{ r }  },  \quad 
P = 1,
\ee 
over the computational domain $[r, \varphi] \in  [1,2] \times [0, 2 \pi]  $.
At the initial time, we perturb this equilibrium solution by imposing a higher density $ \rho = 2 $ within the disc defined in Cartesian coordinates  as $ (x+1.5)^2 + y^2 \le (0.15)^2 $. 

\begin{figure*}
	\includegraphics[width=0.45\linewidth]{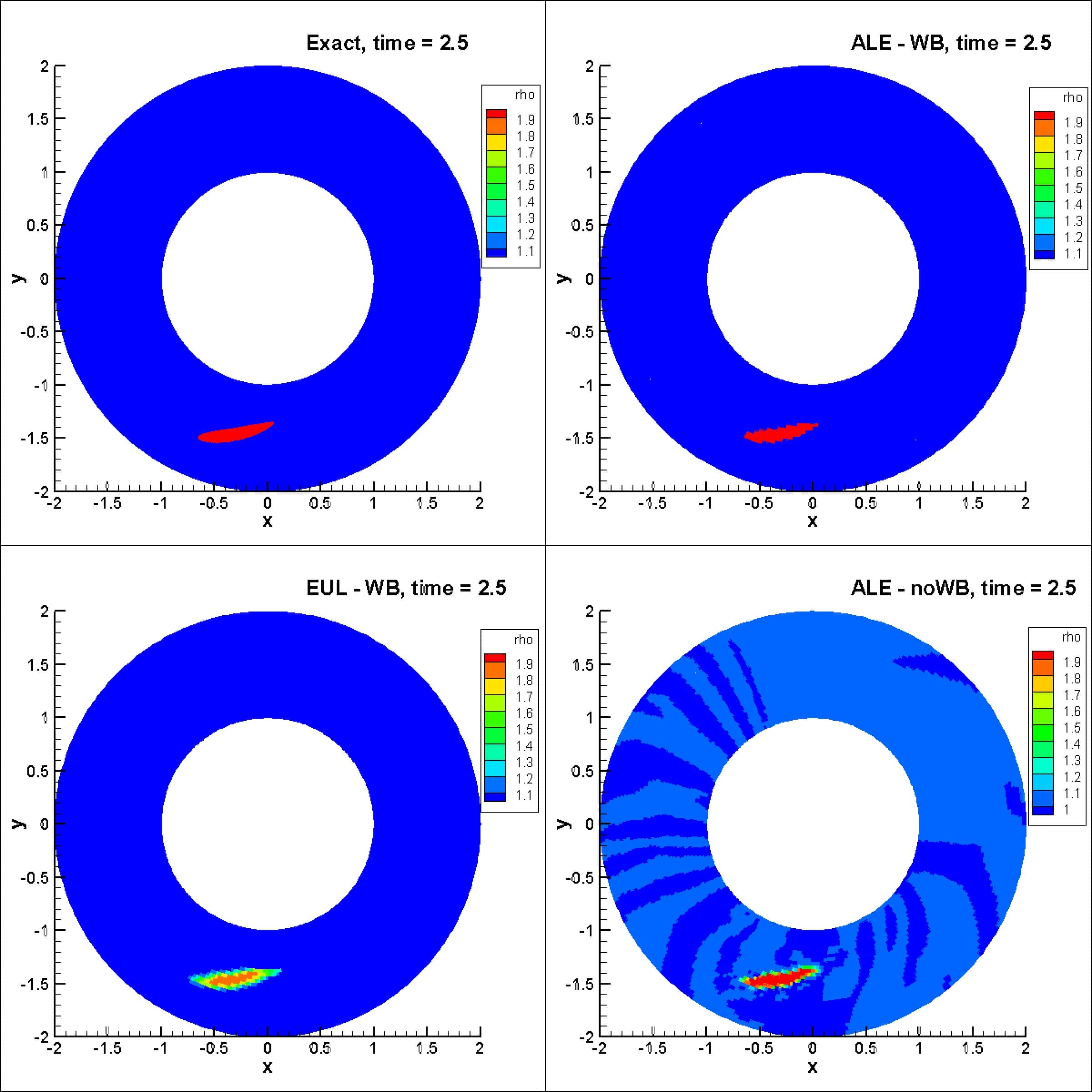} \quad 
	\includegraphics[width=0.45\linewidth]{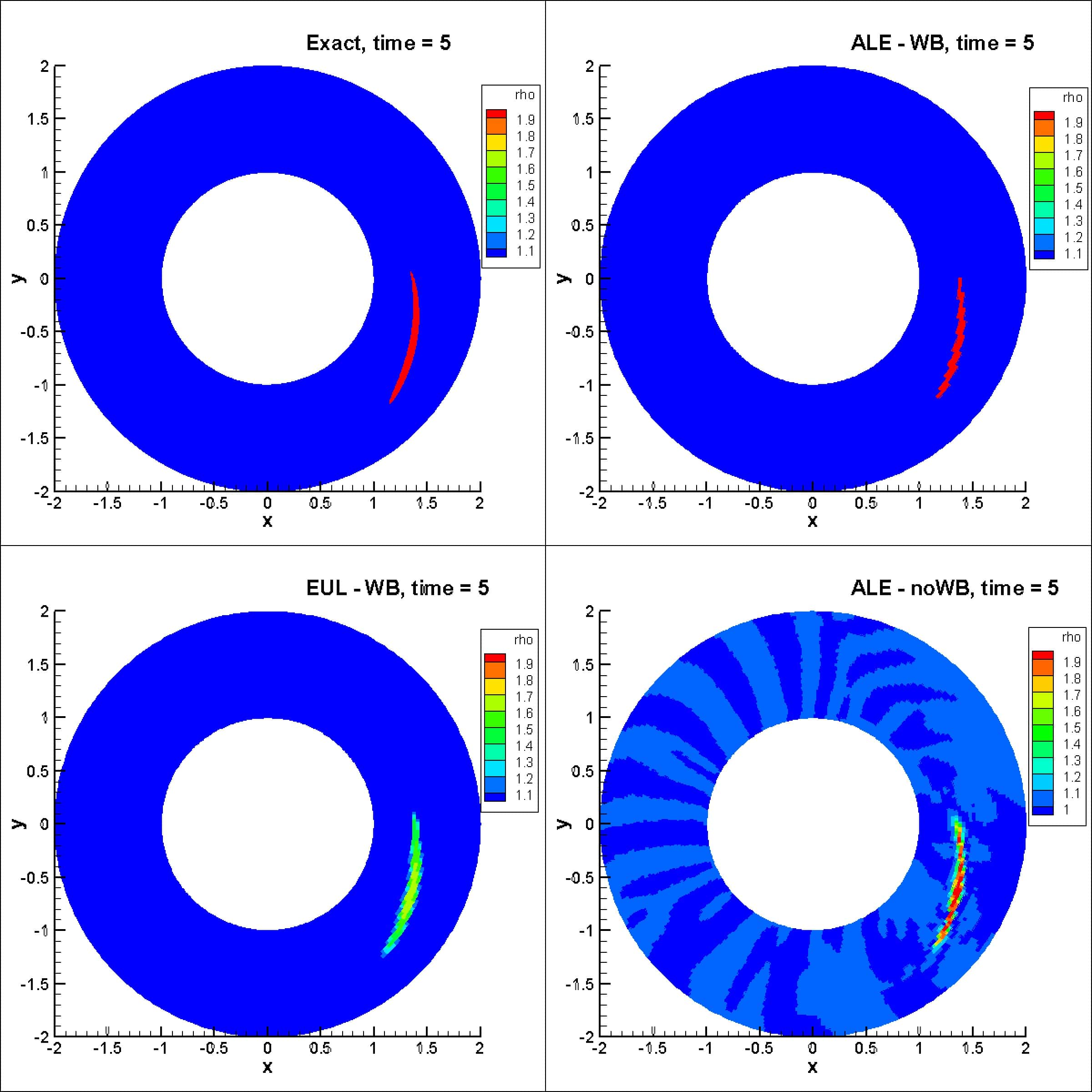}	 \\[8pt]
		\includegraphics[width=0.45\linewidth]{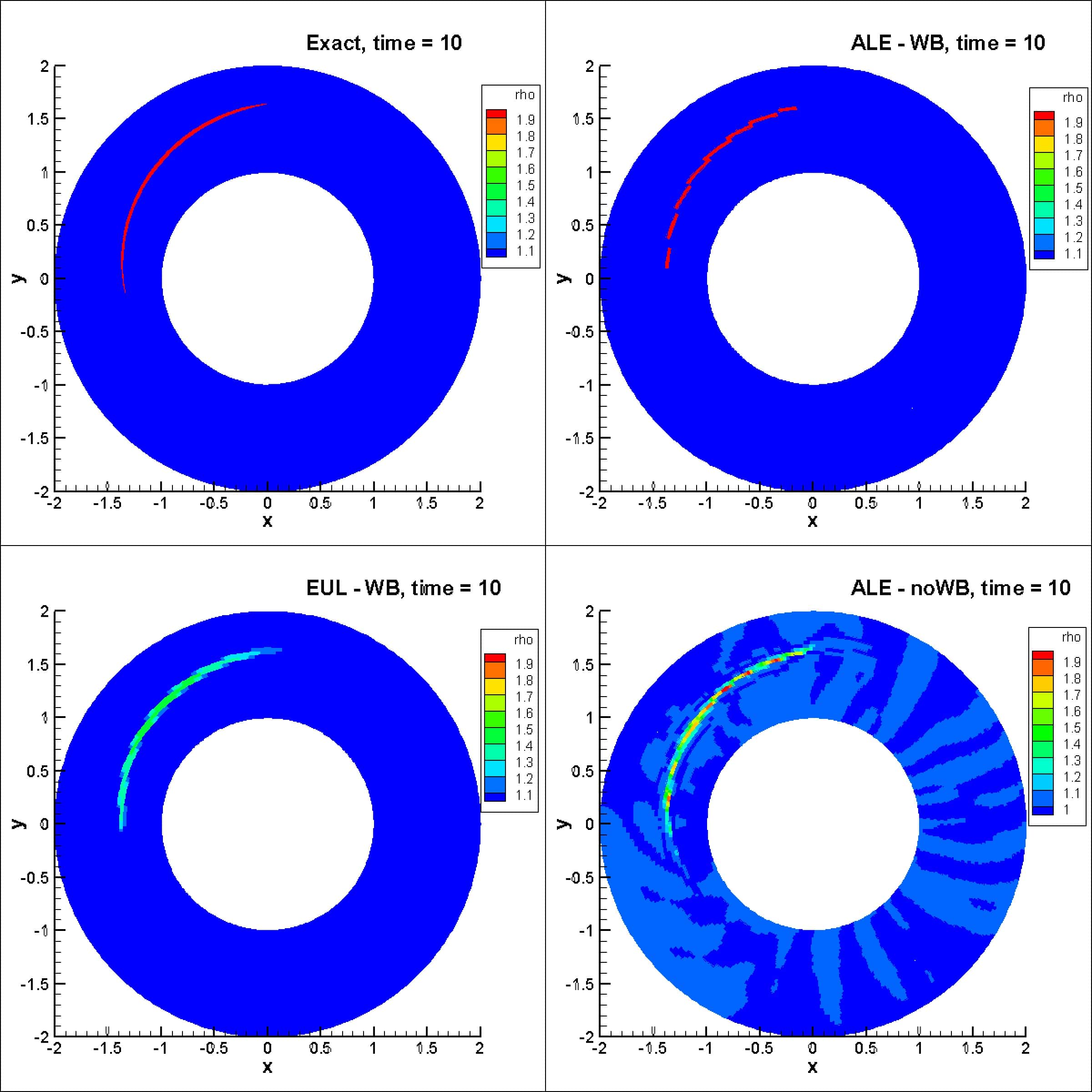} \quad		
			\includegraphics[width=0.45\linewidth]{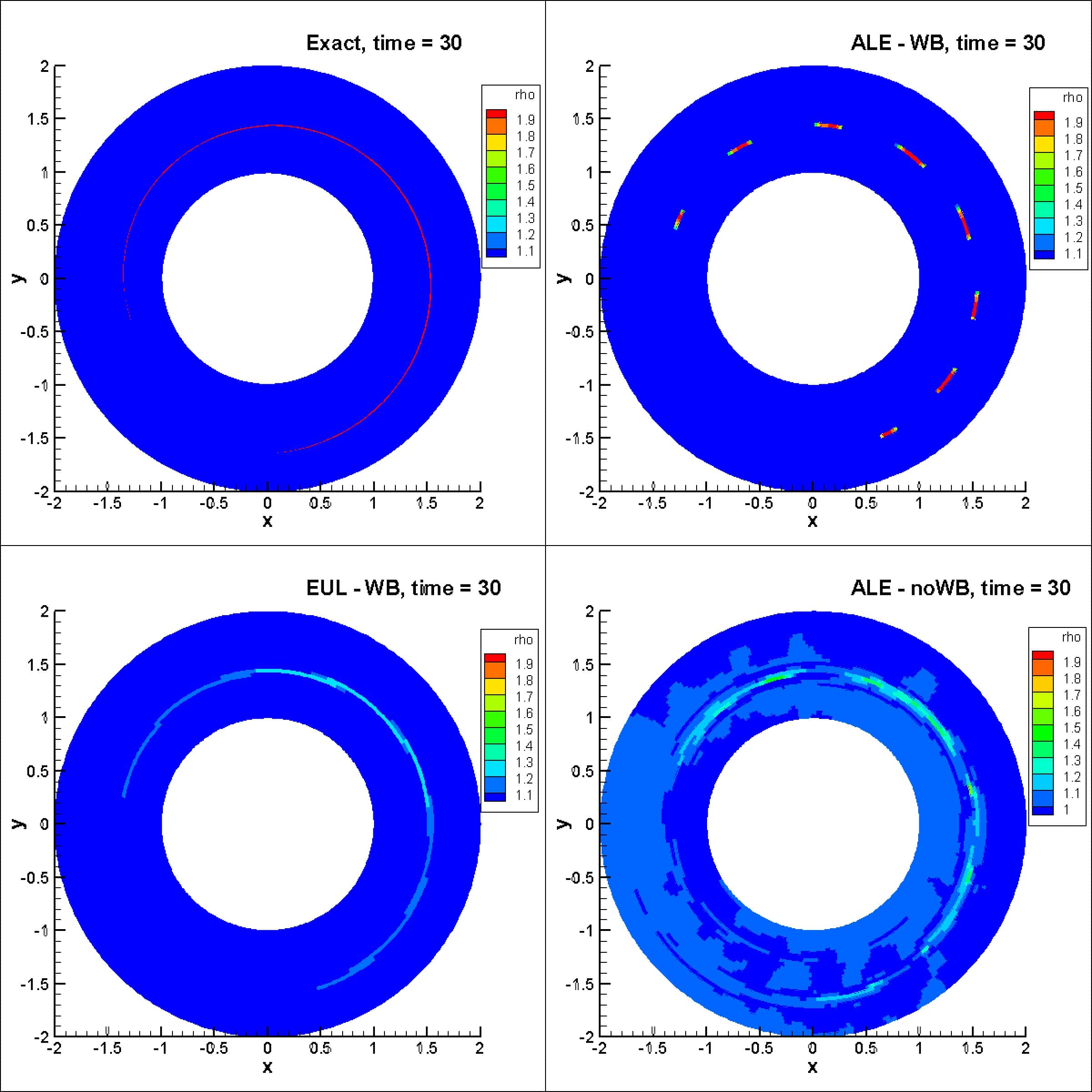} \\	
			\caption{We compare the exact solution with the numerical solutions obtained with different methods at times $t = 2.5$ (top-left), $t=5$ (top-right), $t=10$ (bottom-left), $t=30$ (bottom-right). For all the cases the employed numerical flux is an Osher-type flux. The Lagrangian algorithms show their ability in reducing the viscosity along the angular direction. The well balanced methods do not diffuse the quantities in the radial direction. When coupled together (top-right of each square) we obtain a result very close to the exact solution (top-left of each square). We want to remark that in the well balanced ALE case (top-right of each square), the quantity with higher density remains in the same cells in which it is confined at the initial time since the method is very little diffusive in any direction and the differential rotation is treated in a nonconforming way. So, after long times, the cells containing the higher density gas are no more close to each others, and this explains the figure at time $t=30$. Moreover, only the well balanced ALE scheme is able to maintain the concentration of the higher density gas. }
				\label{fig.KeplerianDisk_transport_results}
\end{figure*}

The expected result is the transport of this density fluctuation (contact discontinuity) at different velocities which are bigger at the interior and smaller at the exterior, without any dissipation. 
The velocity and the pressure field should remain constant in time, according to the equilibrium solution. 
In Figure~\ref{fig.KeplerianDisk_transport_results} we compare the results obtained with different numerical methods with the exact solution: Eulerian and ALE schemes coupled or not with the well balanced Osher Romberg scheme. 
As expected, the Eulerian scheme is very dissipative, even when coupled with our new well balanced technique. The dissipation is evident in the angular direction, since the radial velocity in this problem is zero and the Osher scheme is a complete Riemann solver that is able to resolve steady contact waves exactly. 
The ALE scheme, without well balancing does not dissipate too much in the angular direction, but if it is not coupled with a proper well balanced technique, some spurious velocity oscillations appear which lead to unphysical dissipation in the radial direction and which also produce some oscillations on the density profiles, which are evident even for short computational times. 

The coupling between the two techniques reduces the dissipation both in the radial and in the angular directions. In the computations performed with our well-balanced schemes we 
have observed that for this test problem the error in the pressure and in the velocity field was always of the order of \textit{machine precision}, since the advection of a contact discontinuity 
does \textit{not} affect the equilibrium of pressure and velocity. We emphasize that this property of conserving even non stationary equilibria (density is not constant in time here) is anything else
than trivial to achieve and to the best knowledge of the authors, the scheme presented in this paper is the first finite volume method to achieve it. 
{Referring to Table \ref{tab.EqProfile_v_and_p}, one can notice that indeed the precision achieved by our code on angular velocity and pressure is of the order of machine precision (even at time $t=30$), where instead this is not the case for various PLUTO configurations.}

\begin{figure*}
 \includegraphics[width=0.24\linewidth]{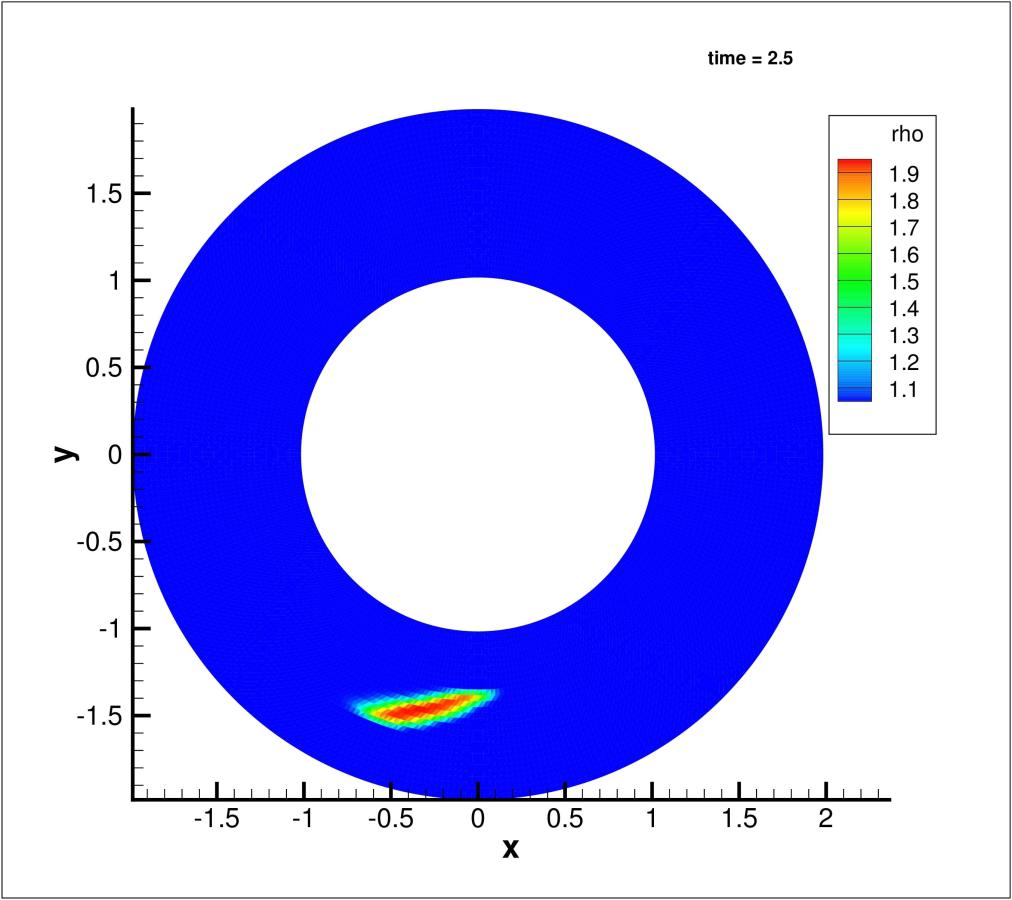} 
	\includegraphics[width=0.24\linewidth]{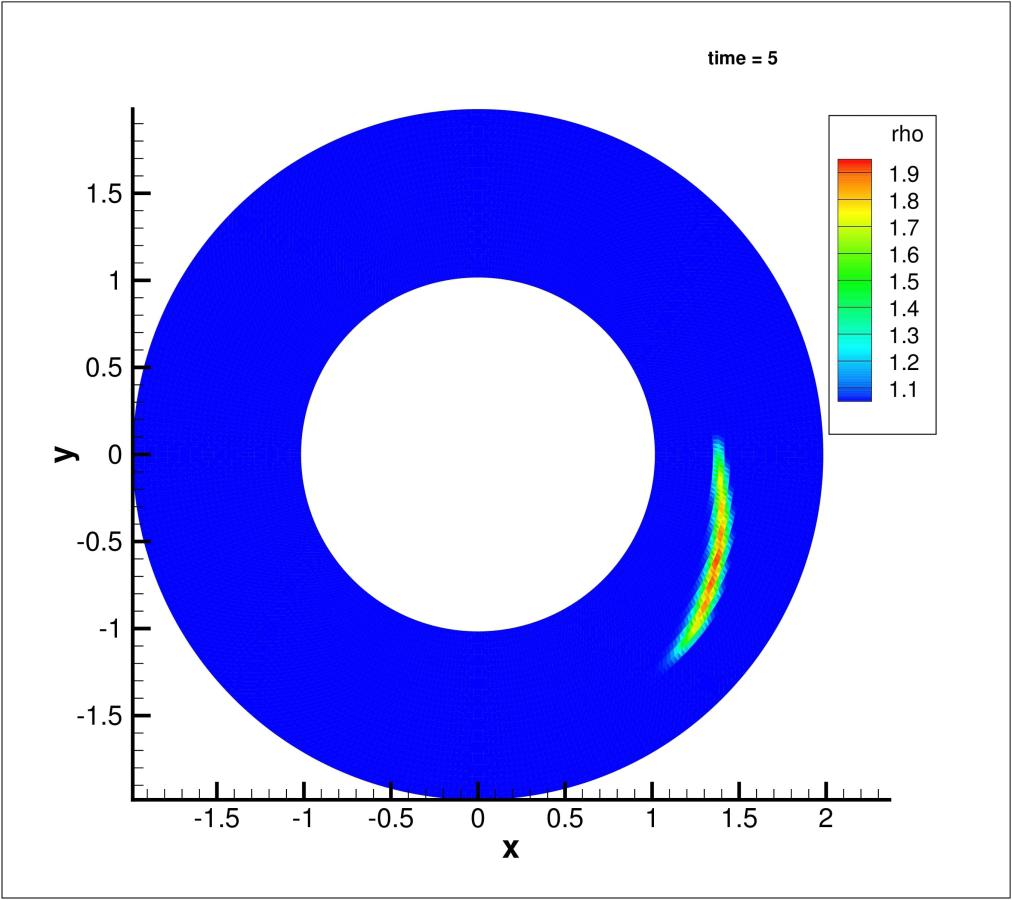}
	\includegraphics[width=0.24\linewidth]{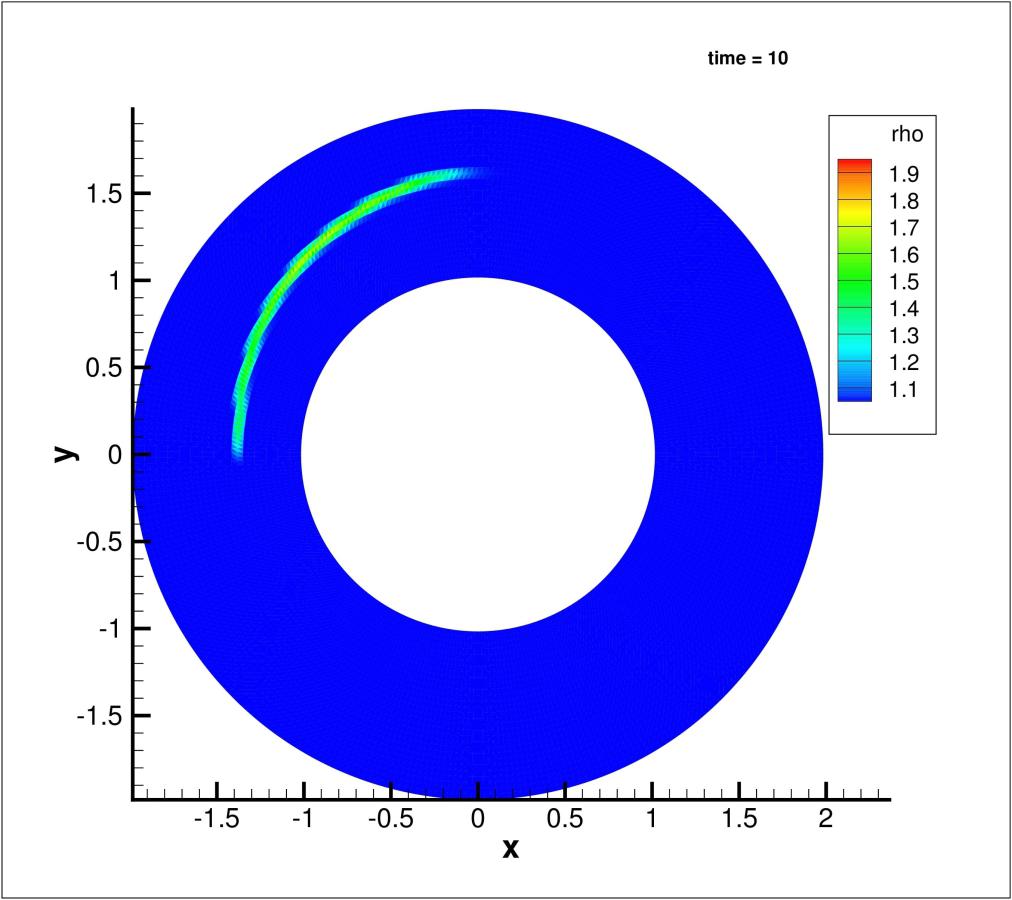} 
	\includegraphics[width=0.24\linewidth]{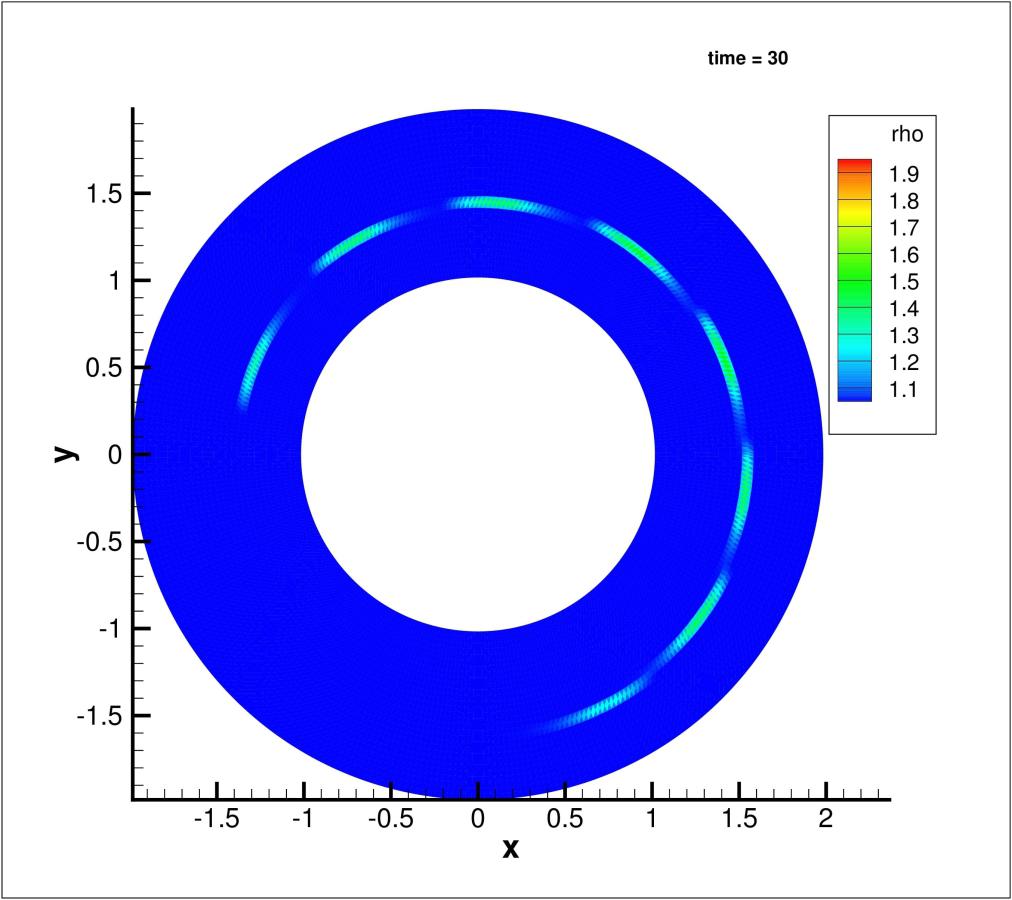} 
	\caption{{Results obtained with PLUTO, using the Roe solver combined with the mc\_lim limiter, linear reconstruction in space and RK2 in time on a grid of $30\times 350$ elements. One can observe that the results are more dissipative compared to those shown in Fig. \ref{fig.KeplerianDisk_transport_results}.}}
	\label{fig.KeplerianDisk_transport_resultsPLUTO_a}
\end{figure*}

\begin{figure*}
	\includegraphics[width=0.24\linewidth]{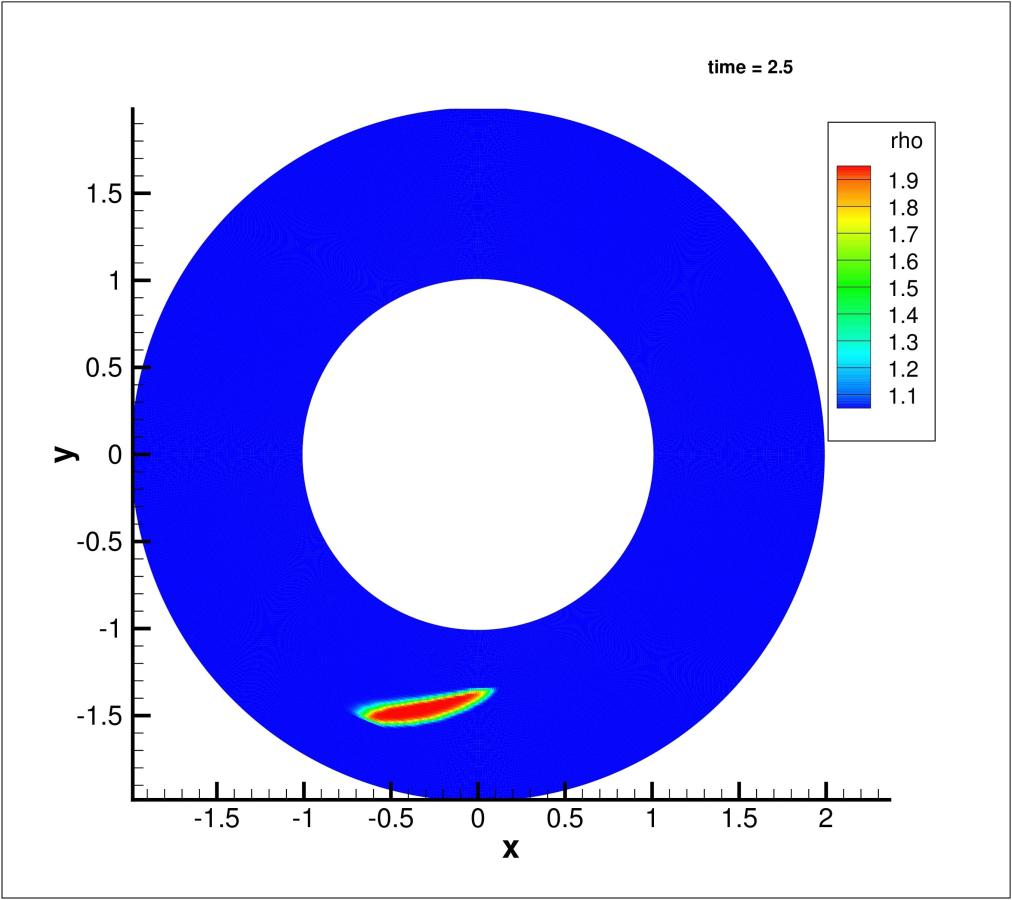} 
	\includegraphics[width=0.24\linewidth]{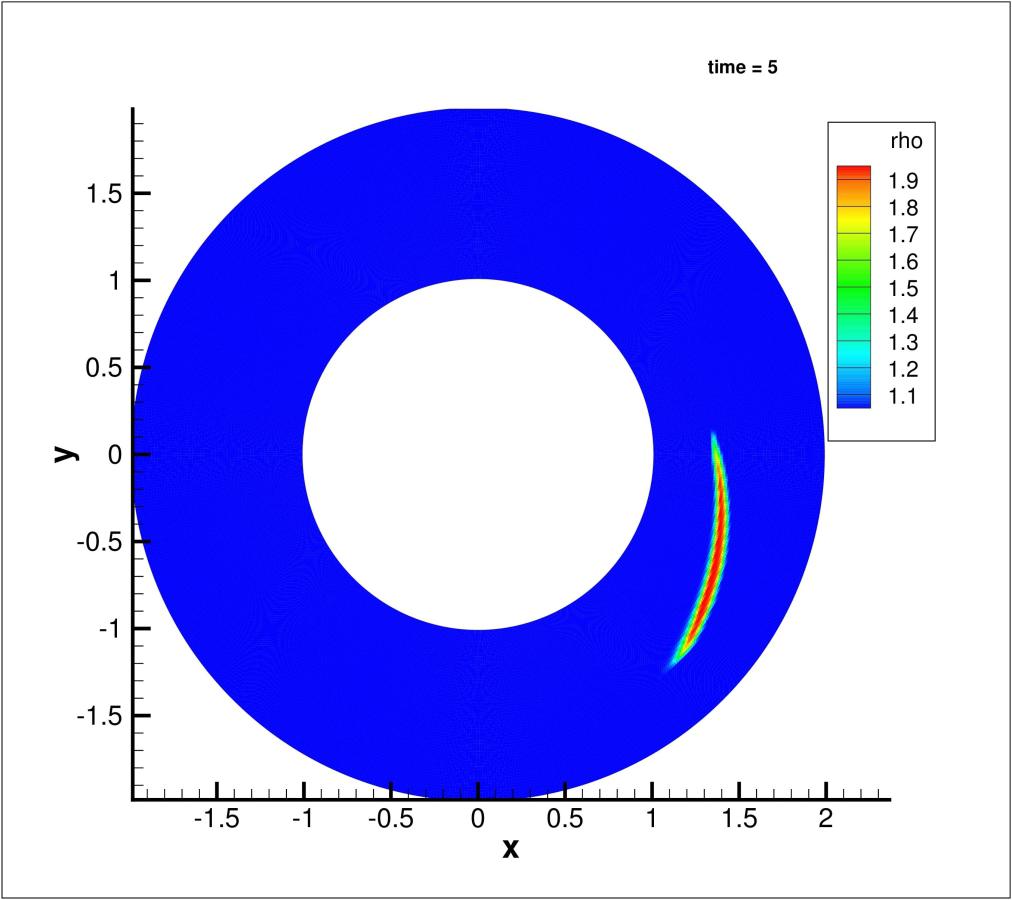}
	\includegraphics[width=0.24\linewidth]{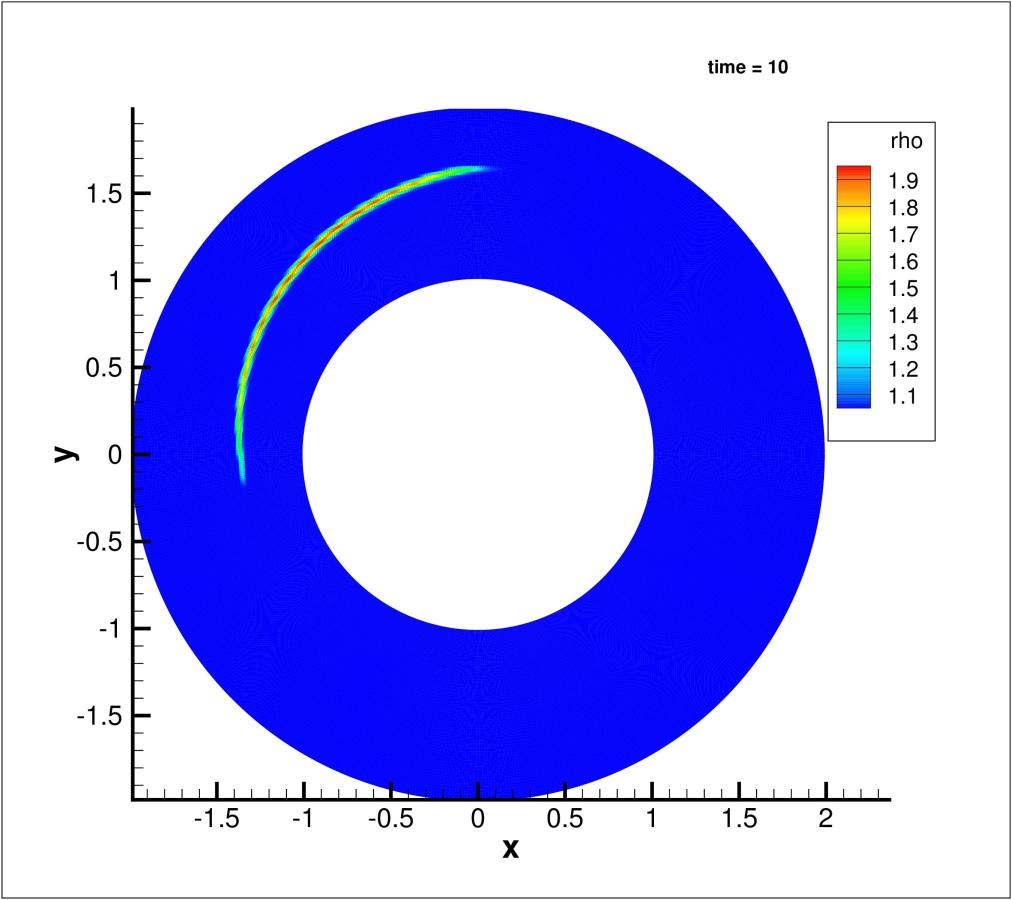} 
	\includegraphics[width=0.24\linewidth]{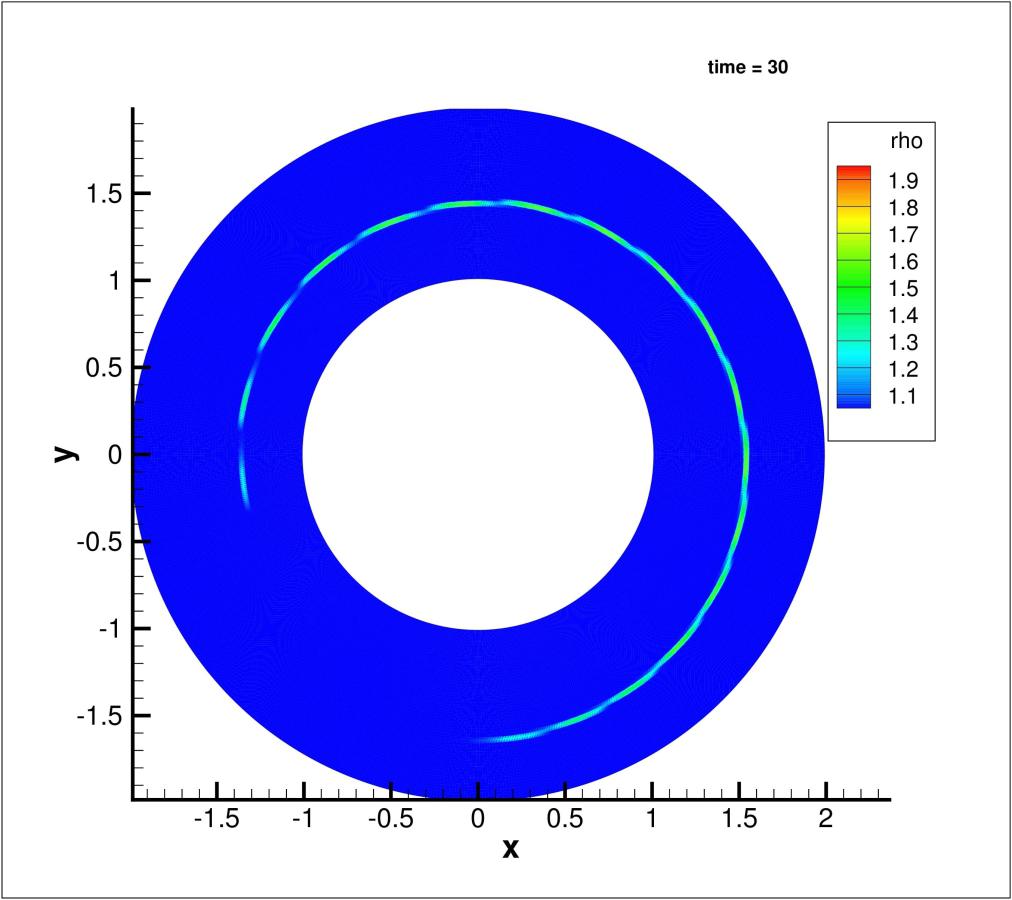} 
		\caption{{Results obtained with PLUTO, using the Roe solver, a third order WENO reconstruction in space combined with the mc\_lim limiter and a third order RK3 time integrator on a grid of $60\times 700$ elements.}}
	\label{fig.KeplerianDisk_transport_resultsPLUTO_b}
\end{figure*}

{Finally, we report the results obtained with PLUTO by selecting the configuration setting described in Section \ref{ssec.PLUTO_setting} with the \texttt{MC\_LIM}. 
First, in Figure \ref{fig.KeplerianDisk_transport_resultsPLUTO_a} we use the described second order method and $30\times350$ elements. Then, in Figure \ref{fig.KeplerianDisk_transport_resultsPLUTO_b} we use the third order method and $60\times700$ elements. In both the cases the density is dissipated faster than with our method: this shows that it is not a finer grid or a higher order of accuracy that can solve this type of problem, but a very specific treatment of the equilibrium together with the Lagrangian framework proposed in this paper. 
}


\begin{table*}	
	\caption{{The results shown in this table testify that our code is able to maintain up to machine precision even non stationary equilibria. Indeed for  the test cases presented both in Section \ref{sec.KeplDisk_transport1}  and Section \ref{sec.KeplDisk_transport2} the $L_1$ norm of the difference between the numerical solution computed with our WB ALE Osher Romberg scheme and the exact stationary profiles of angular velocity $v$ and pressure $P$, at the respective final times ($t=30$ and $t=15$), is of the order of machine precision. The other two lines refer to the results obtained with PLUTO both with second and third order of accuracy.}}  
	\begin{center}
		\begin{tabular}{ll||cc|cc} 
			\hline  \multicolumn{2}{l||}{} & \multicolumn{2}{c|}{Test Section \ref{sec.KeplDisk_transport1}} & \multicolumn{2}{c}{Test Section \ref{sec.KeplDisk_transport2}} \\
			\hline  Method &  Elements &  $|| v - v_E ||_{L_1}$ & $|| P - P_E ||_{L_1}$ &  $|| v - v_E ||_{L_1}$ & $|| P - P_E ||_{L_1}$ 		\\				
			\hline WB ALE Osher-Romberg  & $100 \times 200$ & 2.17E-12   & 7.19E-14  & 2.13E-12  &  6.36E-14 \\ 
			PLUTO O2 minmod  &  $100 \times 200$  &   5.56E-7 & 2.36E-6  & 5.44E-7 & 9.89E-6  \\ 
			PLUTO O3 mc\_lim  &  $200 \times 400$   & 1.30E-7   & 5.28E-7  & 1.49E-7  & 2.44E-6  \\ 
			\hline 
		\end{tabular}	
	\end{center}
	\label{tab.EqProfile_v_and_p}
\end{table*}

\subsection{Keplerian disc with density perturbations} 
\label{sec.KeplDisk_transport2}

\begin{figure*}
	\includegraphics[width=0.24\linewidth]{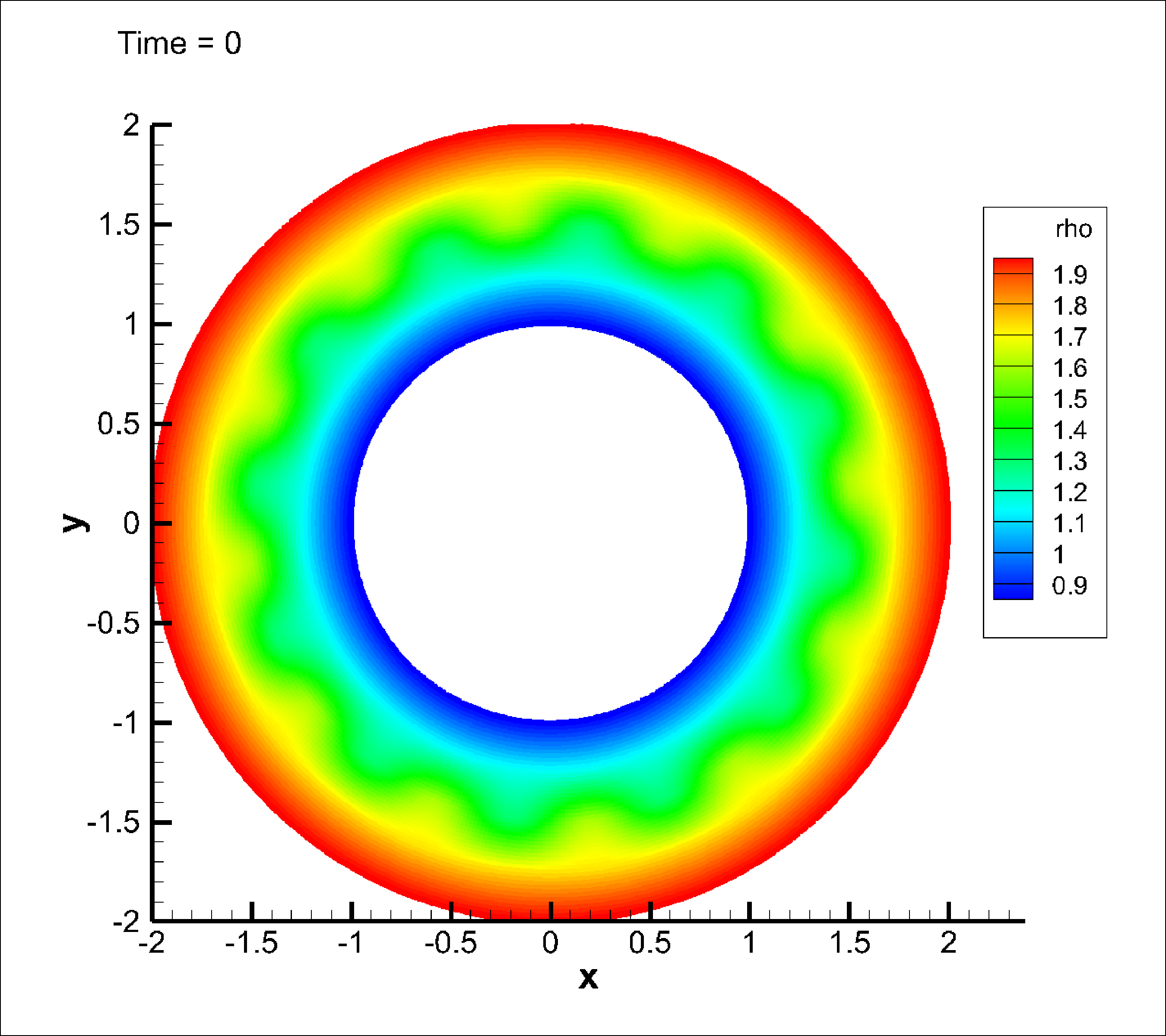} 		
	\includegraphics[width=0.24\linewidth]{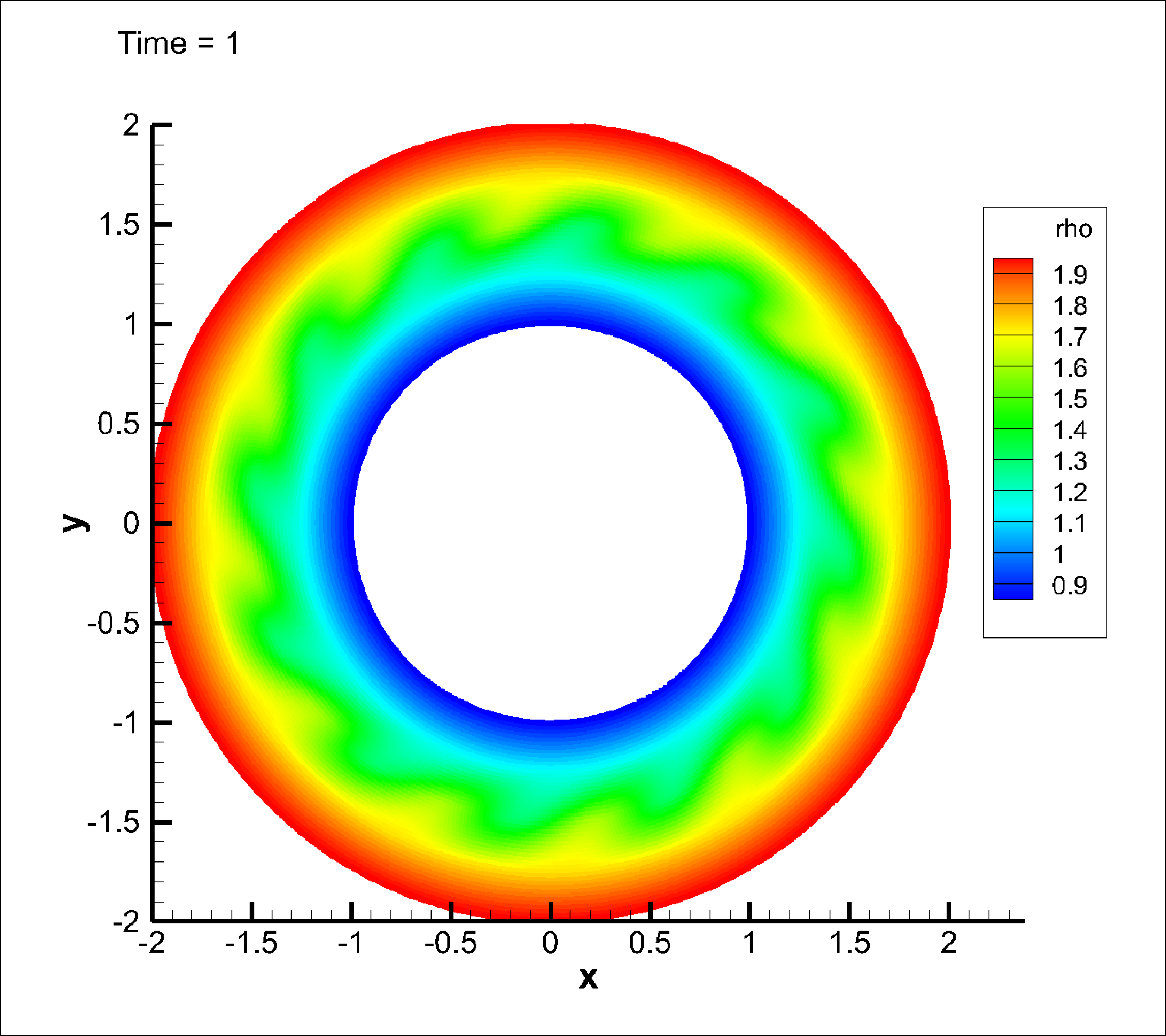} 		
	\includegraphics[width=0.24\linewidth]{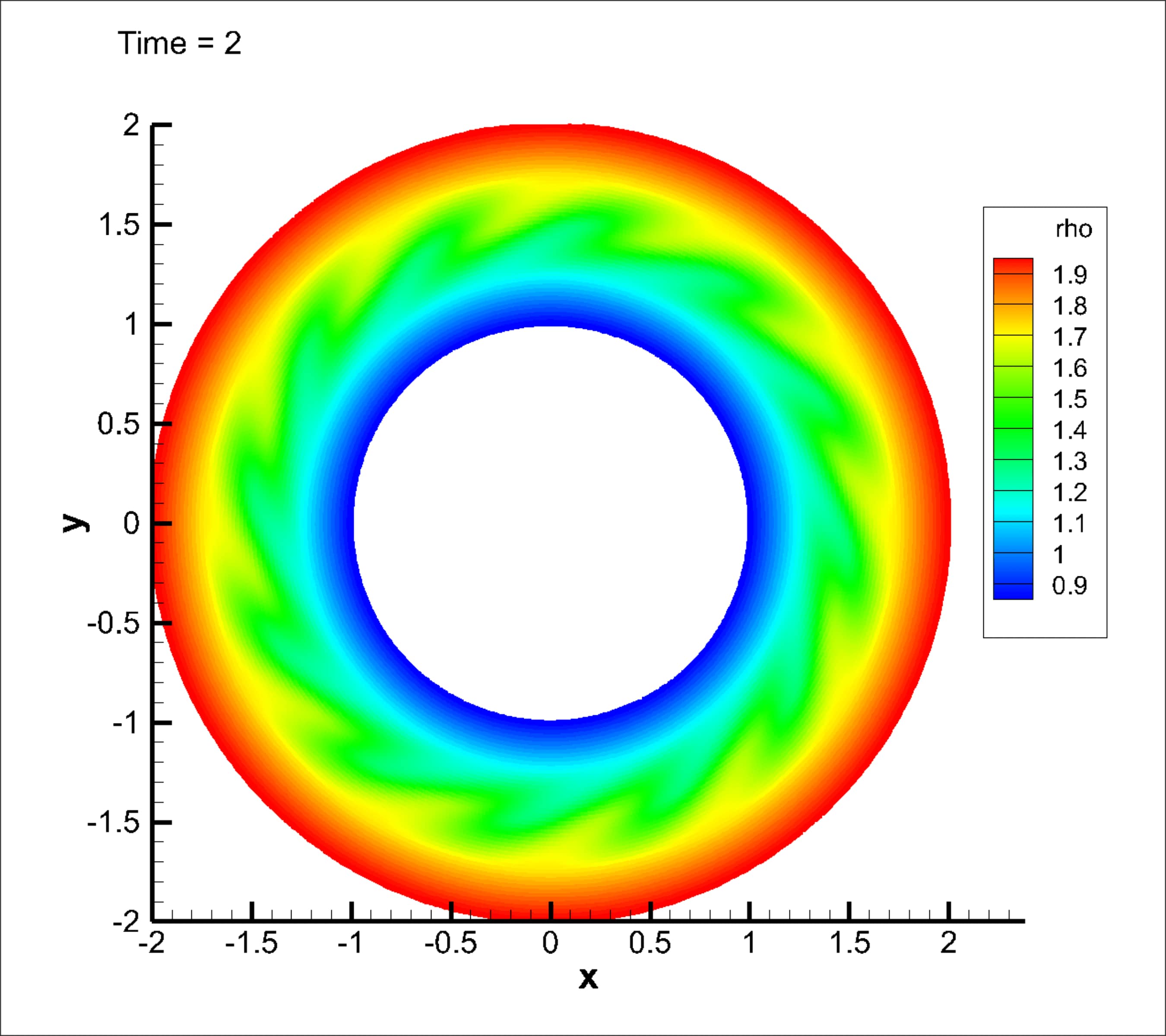} 		
	\includegraphics[width=0.24\linewidth]{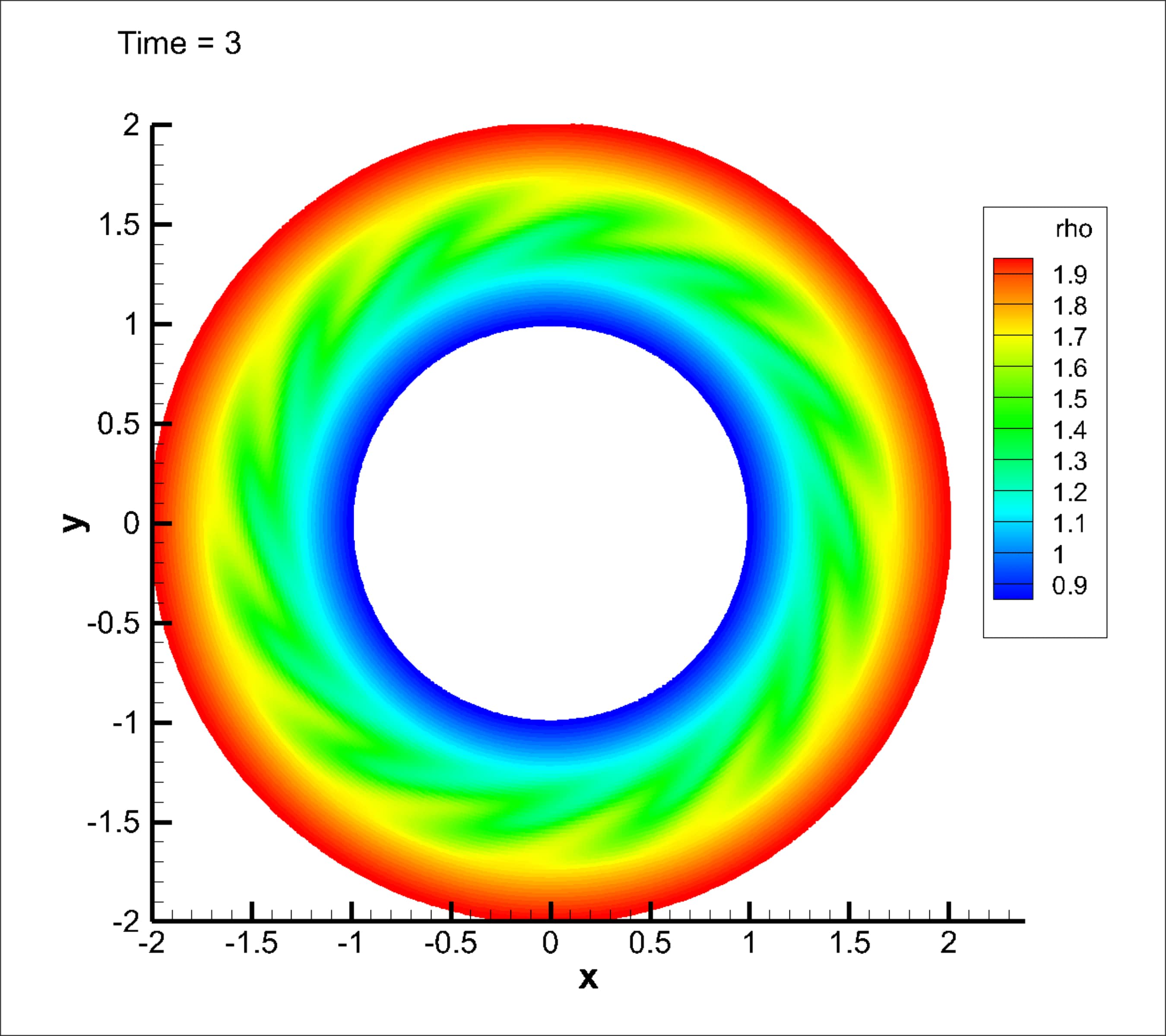} 	 \\
\	\includegraphics[width=0.24\linewidth]{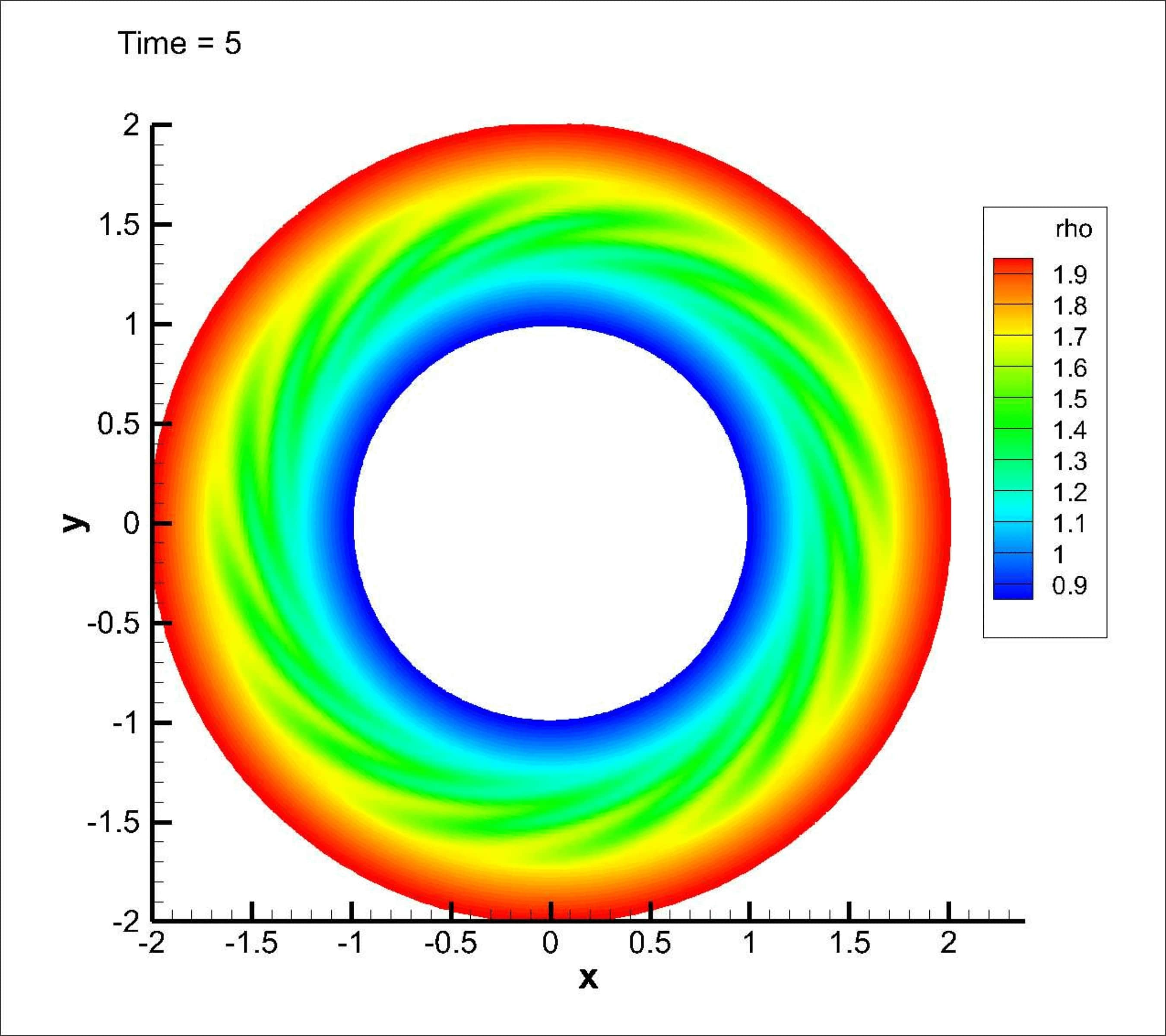}
	\includegraphics[width=0.24\linewidth]{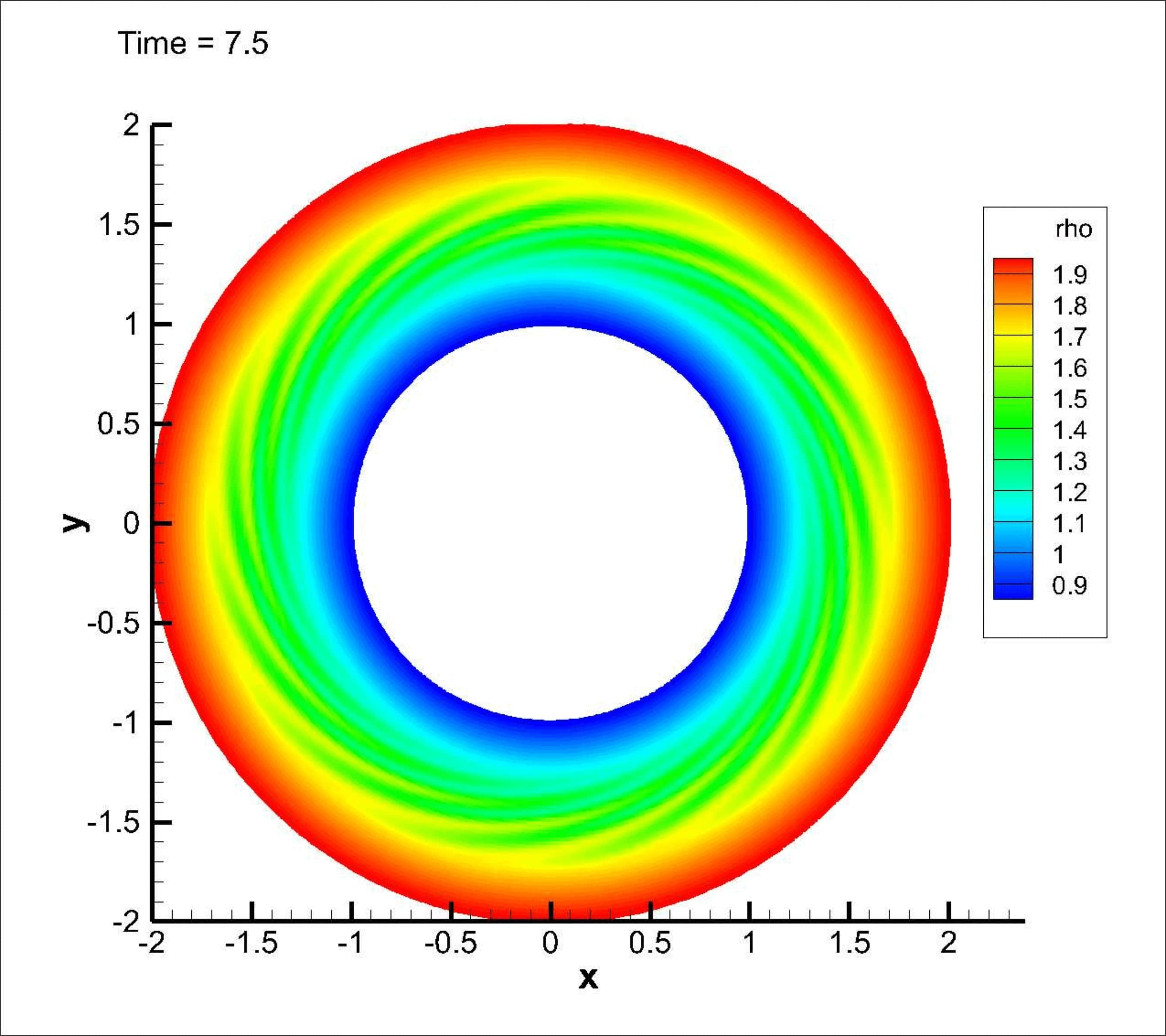}
	\includegraphics[width=0.24\linewidth]{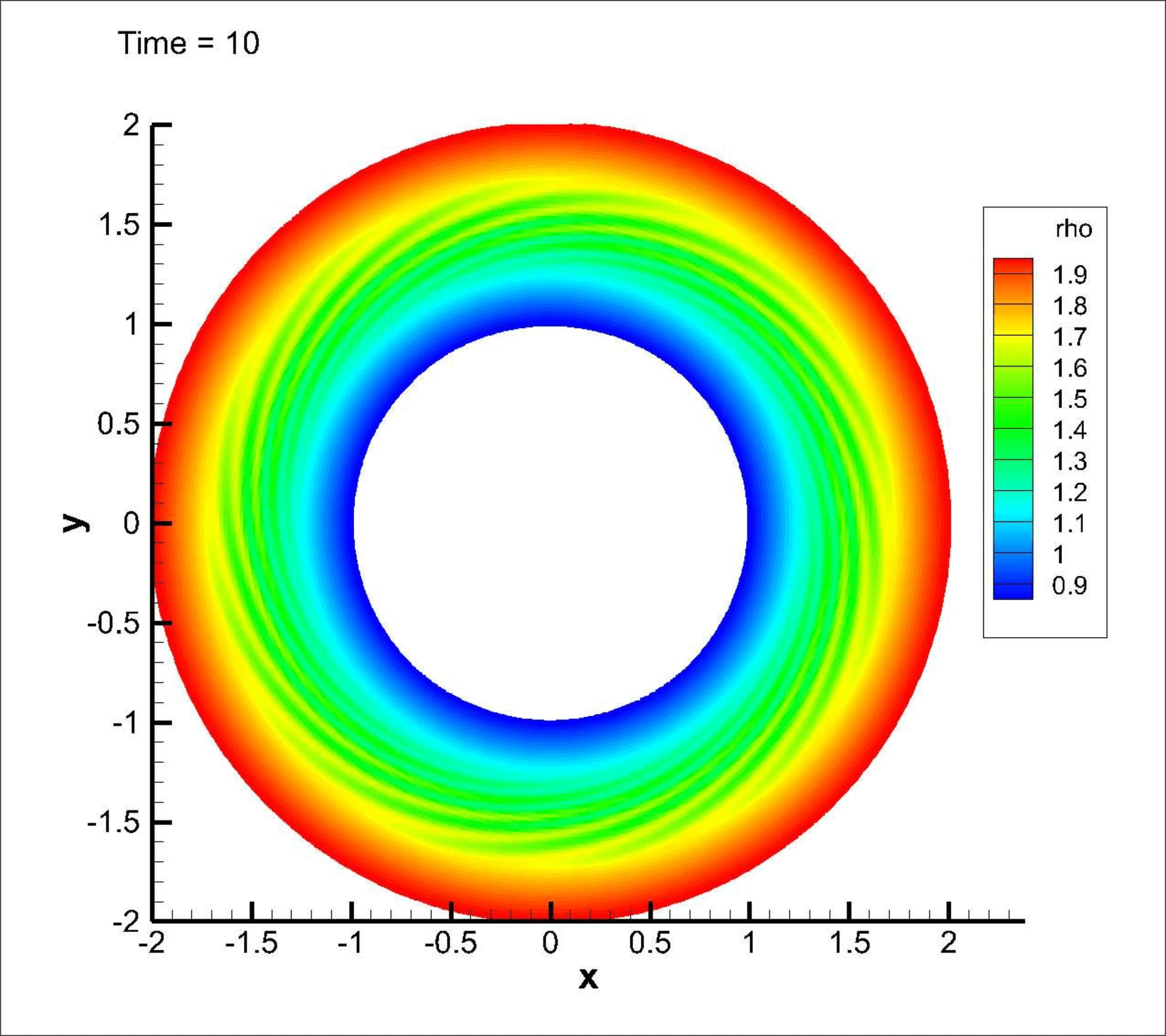} 		
	\includegraphics[width=0.24\linewidth]{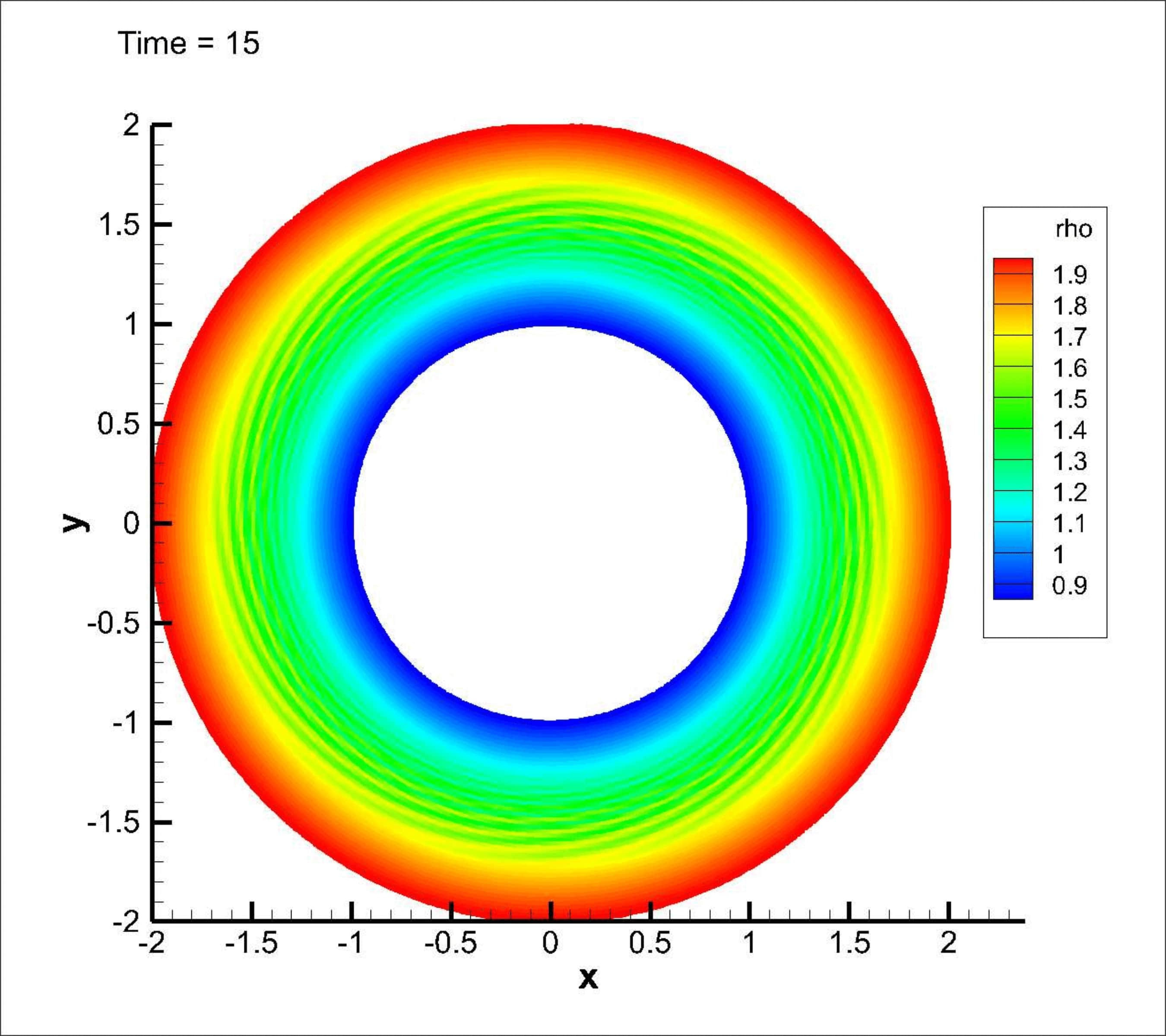}
	\caption{Evolution of periodic density perturbations in an equilibrium disc obtained with the well balanced ALE scheme with Osher-Romberg flux. The perturbations are perfectly convected (with an inner velocity bigger than the outer one), and no spurious Kelvin-Helmholtz instabilities are generated, even after long computational times. }
		\label{fig.KeplerianDisk_transportOfPerturbations}
\end{figure*}

For this test we have considered the equilibrium profile 
\be
\label{eq.constantEquilibrium2}
\rho_E = r, \quad 
u_E = 0, \quad 
v_E = \sqrt{ \frac{G m_s}{ r }  },  \quad 
P = 1,
\ee 
and we have added a periodic perturbation to the density profile as follows
\be
\rho = \rho_E + A \sin (k_1\varphi) ( 0.25 - |r_m - r| ), \ r \in [ r_1, r_2 ]
\ee
with $ A = 0.5$, $k_1 = 12$, $r_1 = 1.25$, $r_2 = 1.75$, $ r_m = 1.5$.
The goal of the this test is to show that our well balanced ALE scheme is able to maintain the equilibrium pressure and velocity exactly and that the numerical method does \textit{not} generate 
any spurious numerical perturbations of pressure and velocity that would usually lead to Kelvin-Helmholtz type flow instabilities for density fluctuations combined with shear flow as in the
above setup.  
In Figure~\ref{fig.KeplerianDisk_transportOfPerturbations} we show the evolution of the perturbations at different times. They are properly transported with different velocities with only very 
little numerical dissipation. As in the previous case we stress that the velocity and pressure remain at the \textit{equilibrium solution} up to \textit{machine accuracy} throughout the entire 
simulation. No spurious Kelvin-Helmhotz instabilities are generated, since the equilibrium pressure and velocity are exactly maintained for arbitrary long simulation times. 

\begin{figure*}
		\includegraphics[width=0.24\linewidth]{TransportOfFluctuations_t15_2}
	\includegraphics[width=0.24\linewidth]{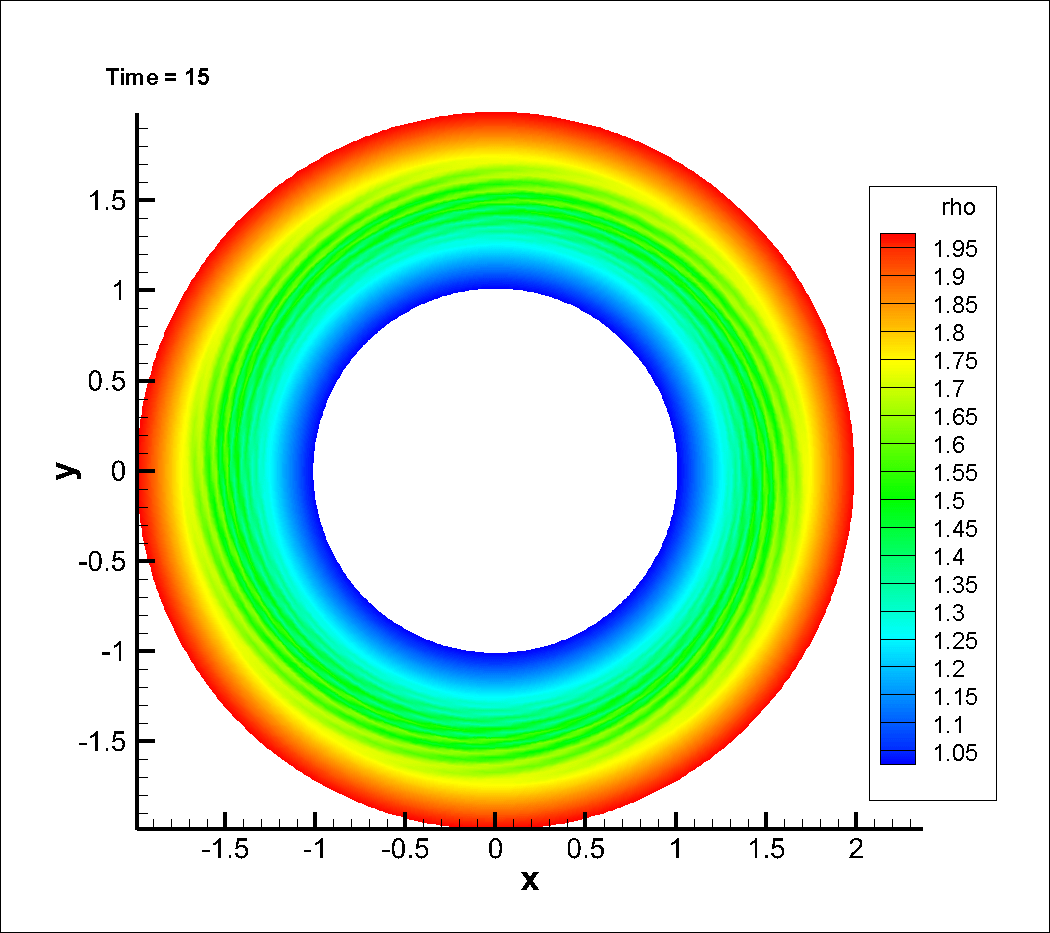}
				\includegraphics[width=0.24\linewidth]{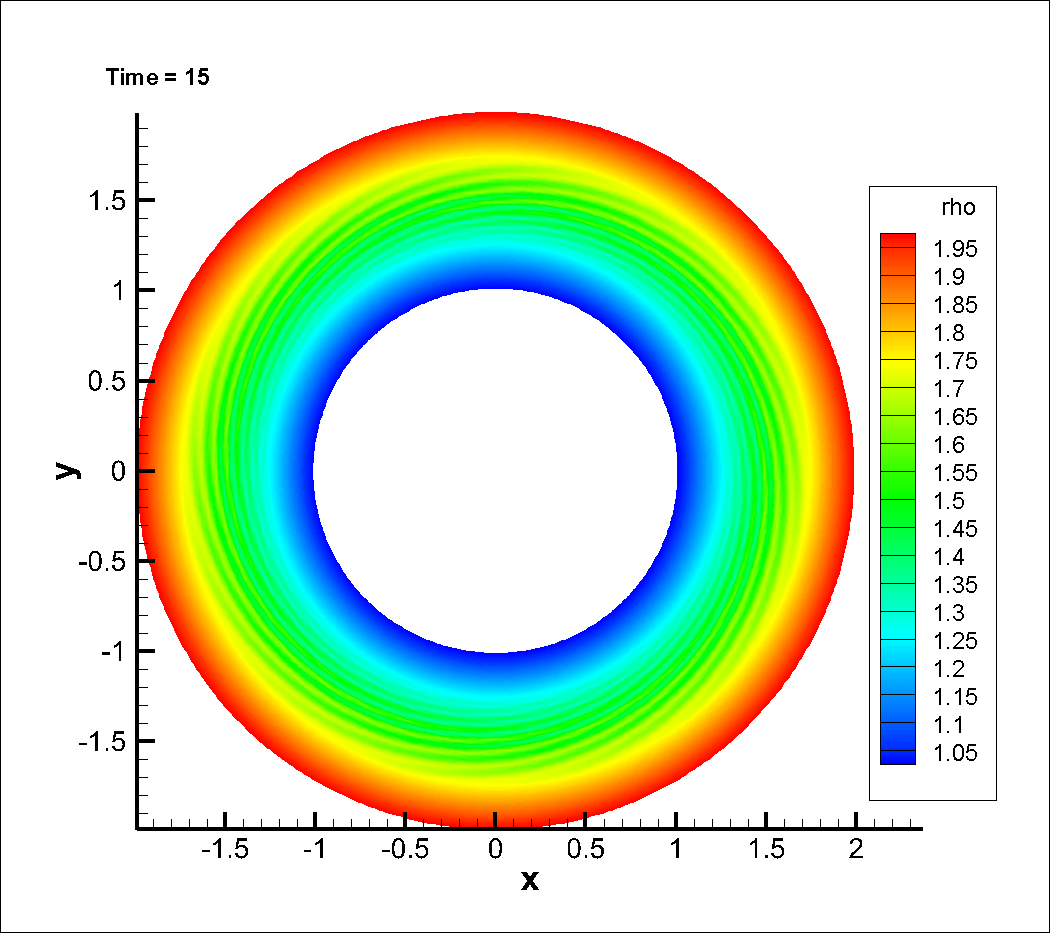}
		\includegraphics[width=0.24\linewidth]{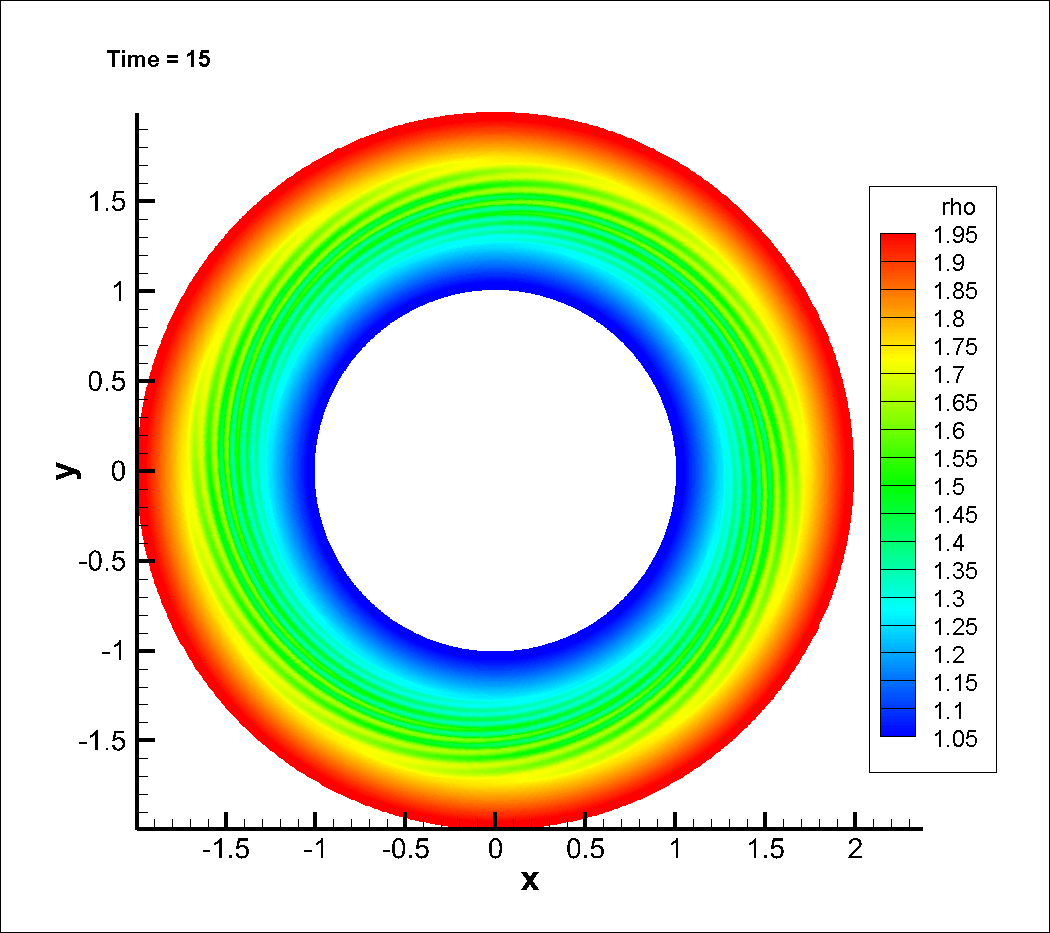}
	\caption{{Method comparison at time $t=15$. The first image is obtained with our code and $50x500$ elements. The second and the third one with PLUTO using $50\times500$ elements and respectively a second order scheme with mc\_lim limiter and a third order scheme with minmod\_lim limiter. The last image is obtained with the third order version of PLUTO using mc\_lim and 
	$100\times1000$ elements. All images are drawn with the same color map. Even if the results are similar, one can notice that to obtain the same resolution of our code we need the third 
	order version of PLUTO and a finer mesh.}}
	\label{fig.KeplerianDisk_PLUTO_transportOfPerturbations}
\end{figure*}

{Finally, we compare our result at time $t=15$ with the results obtained with PLUTO, refer to Figure \ref{fig.KeplerianDisk_PLUTO_transportOfPerturbations}. For the visualization we have always used the software Tecplot and the same colormap; even if the results look similar, one can notice that to obtain the same resolution of our code (left image of the panel) we need the third order version of PLUTO and a finer mesh (last image of the panel). We stress that our code maintains $u$ and $P$ up to machine precision, whereas PLUTO produces standard numerical errors, see Table 
\ref{tab.EqProfile_v_and_p}}

\subsection{Keplerian disc with Kelvin-Helmholtz instabilities I}

Let us consider an equilibrium solution which satisfies the equilibrium constraints in \eqref{eq.EquilibriaConstraint}-\eqref{eq.PressureAndGravForces} so that 
\begin{equation*} 
\label{eq.equilibriumForKH1}
\rho_E = \rho_0 + \rho_1 \text{tanh}\left( \frac{r-r_m}{\sigma} \right),  \, 
u_E = 0, \, 
v_E = \sqrt{ \frac{Gm_s}{r} }, \, 
P_E = 1, 
\end{equation*} 
with $G = 1$, $m_s = 1$, $\rho_0 = 1$, $\rho_1 = 0.25$, $r_m= 1.5$ and $ \sigma = 0.01$. It shows a steep gradient in the density for $ r \rightarrow 1.5$.
We consider as computational domain a ring sector with radius $r \in [1,2]$ and  $\varphi \in [0, \pi/2]$. For the boundary conditions we exploit the exact solution when $r=1,2$, and we impose periodic boundary conditions for $\varphi = 0, \pi/2$.

\begin{figure*}	
	\includegraphics[width=0.24\linewidth]{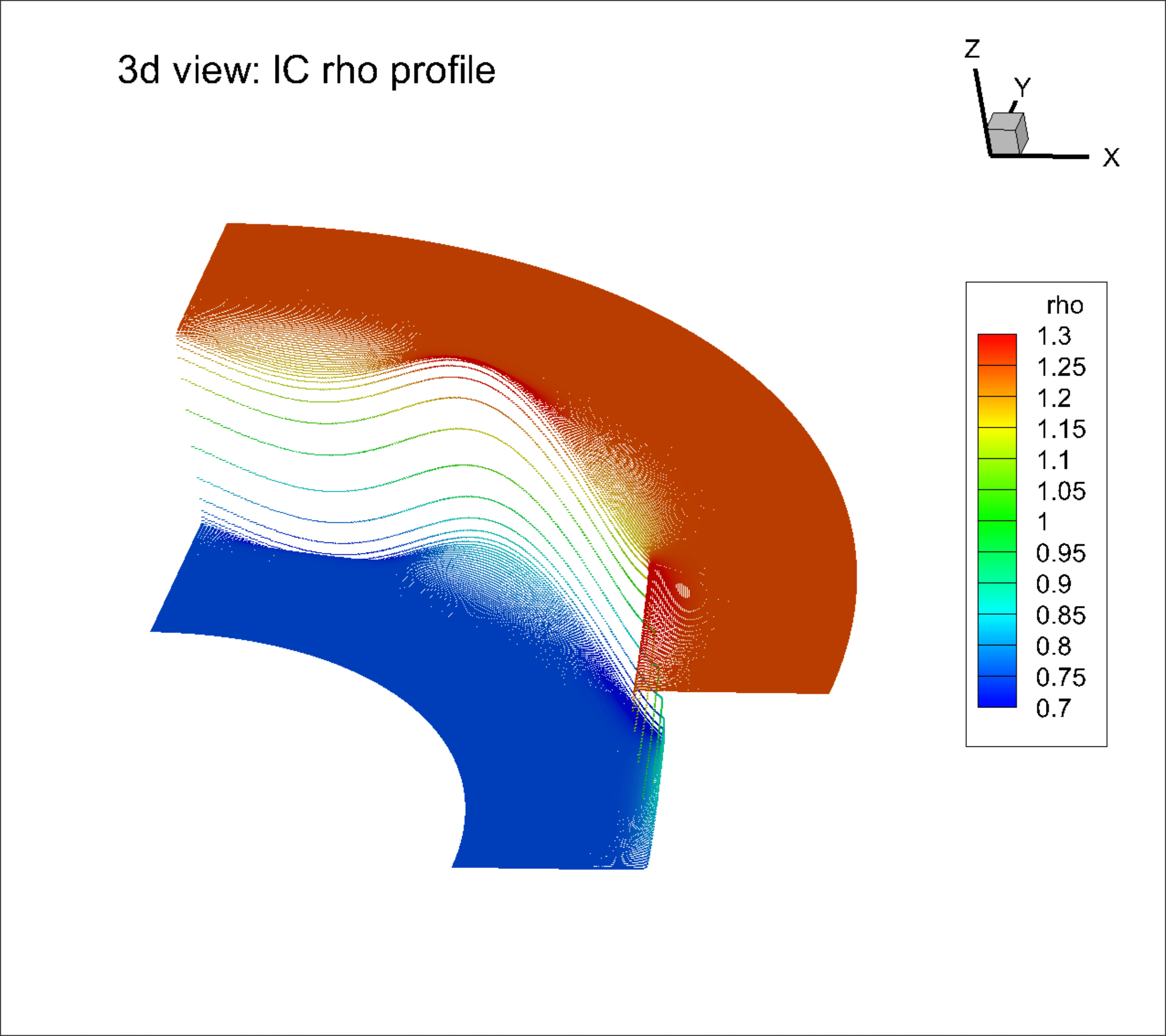} 		
	\includegraphics[width=0.24\linewidth]{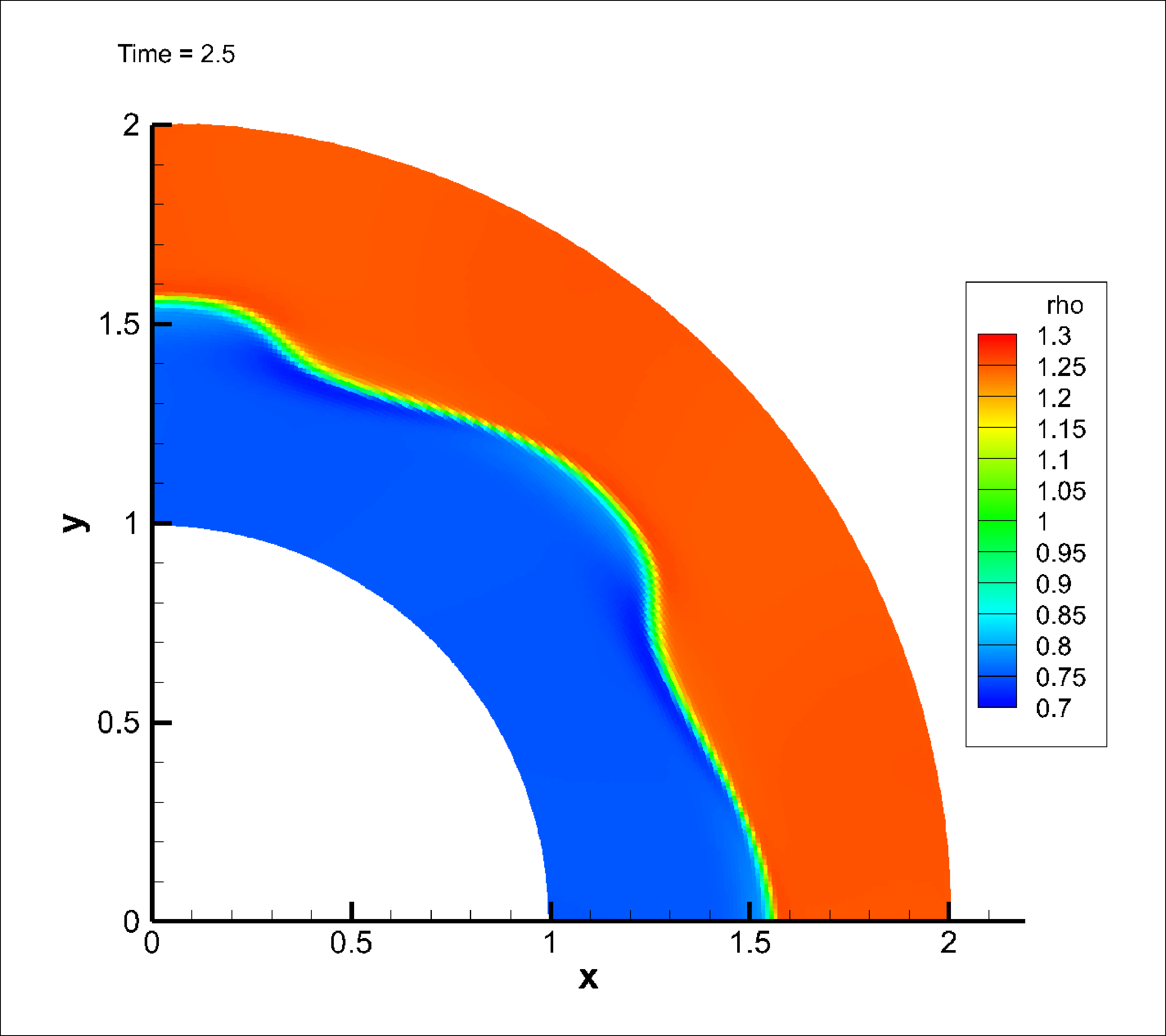} 		
	\includegraphics[width=0.24\linewidth]{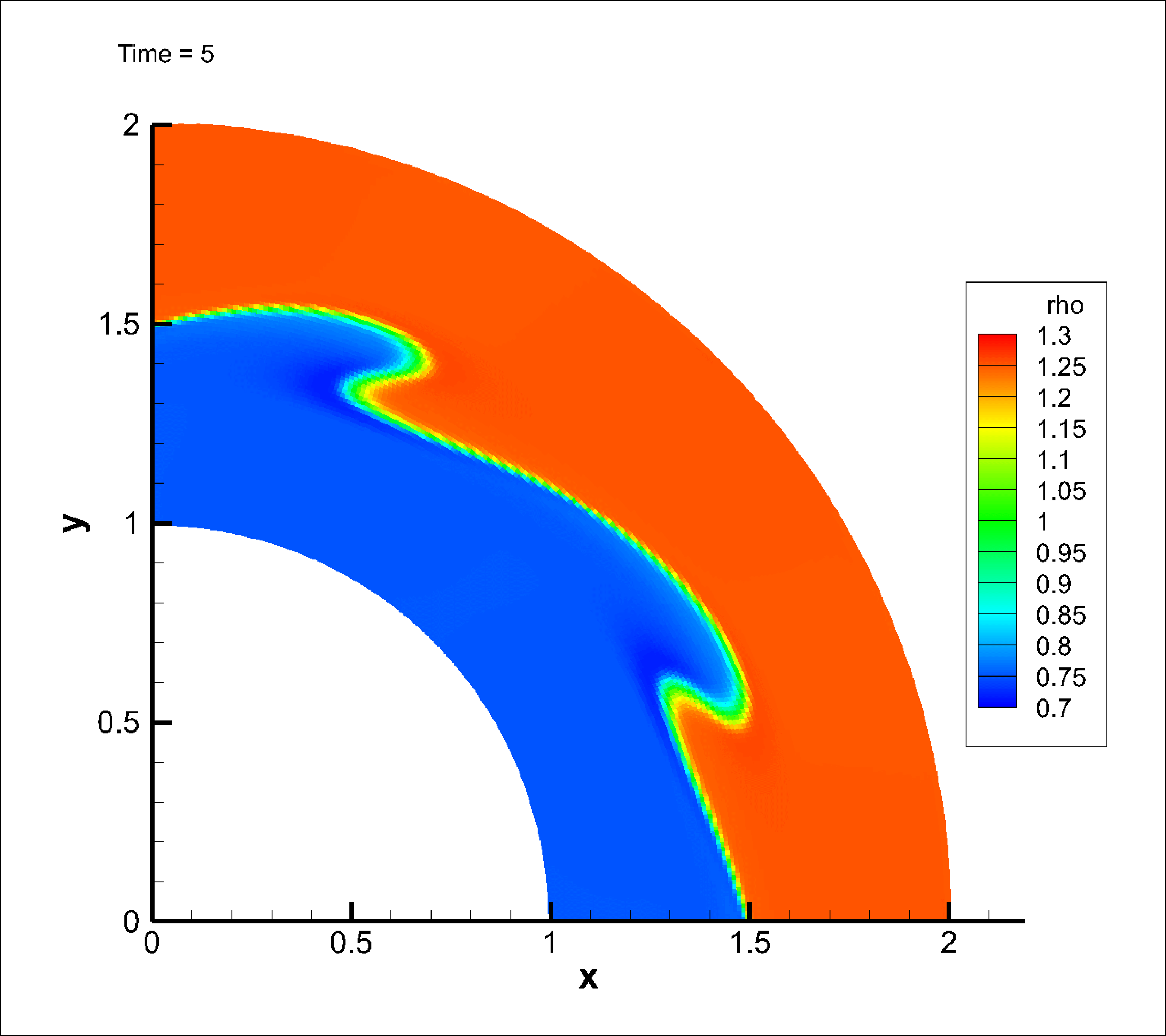} 		
	\includegraphics[width=0.24\linewidth]{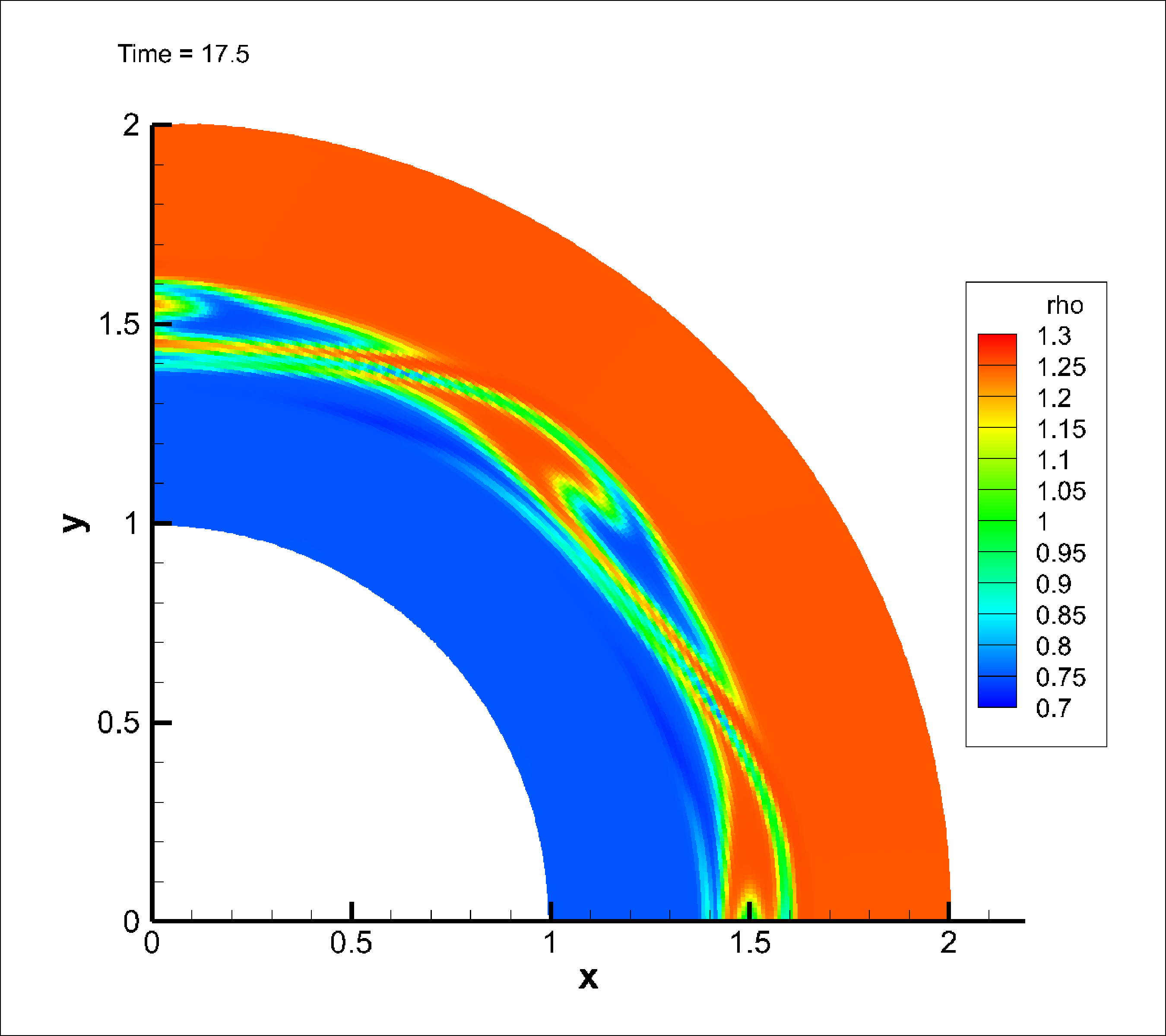} 	 \\
	\ \includegraphics[width=0.24\linewidth]{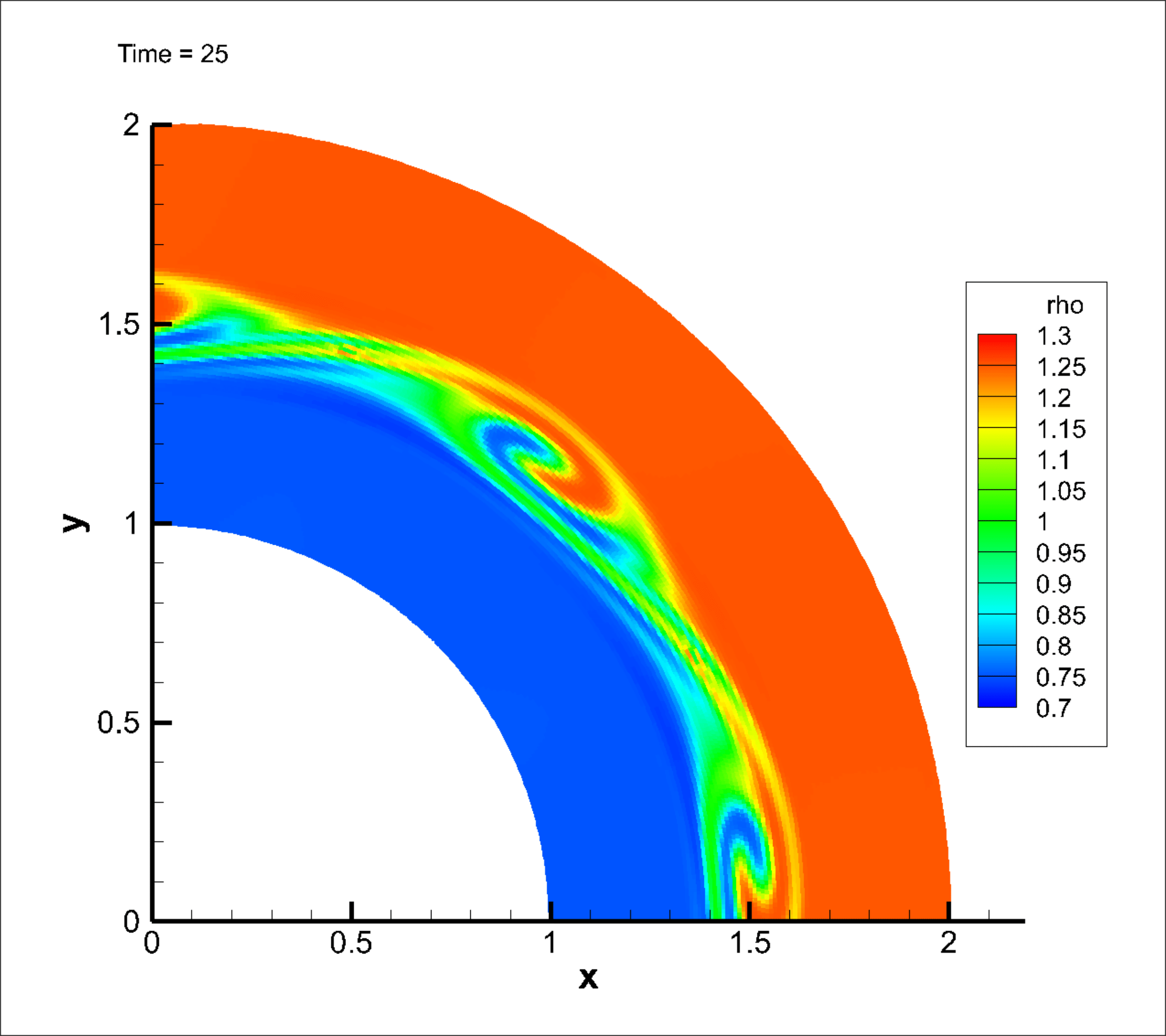}
	\includegraphics[width=0.24\linewidth]{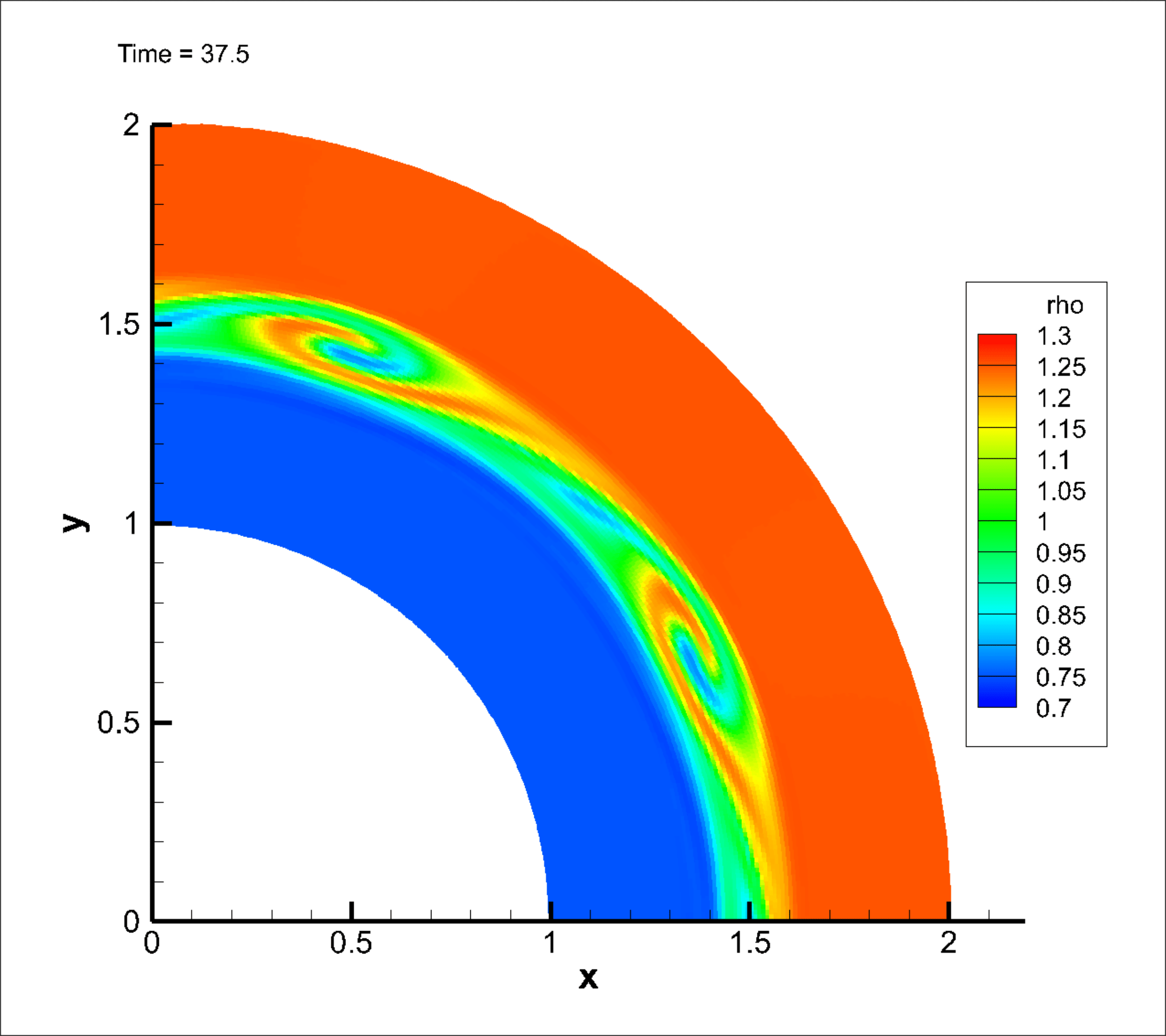}
	\includegraphics[width=0.24\linewidth]{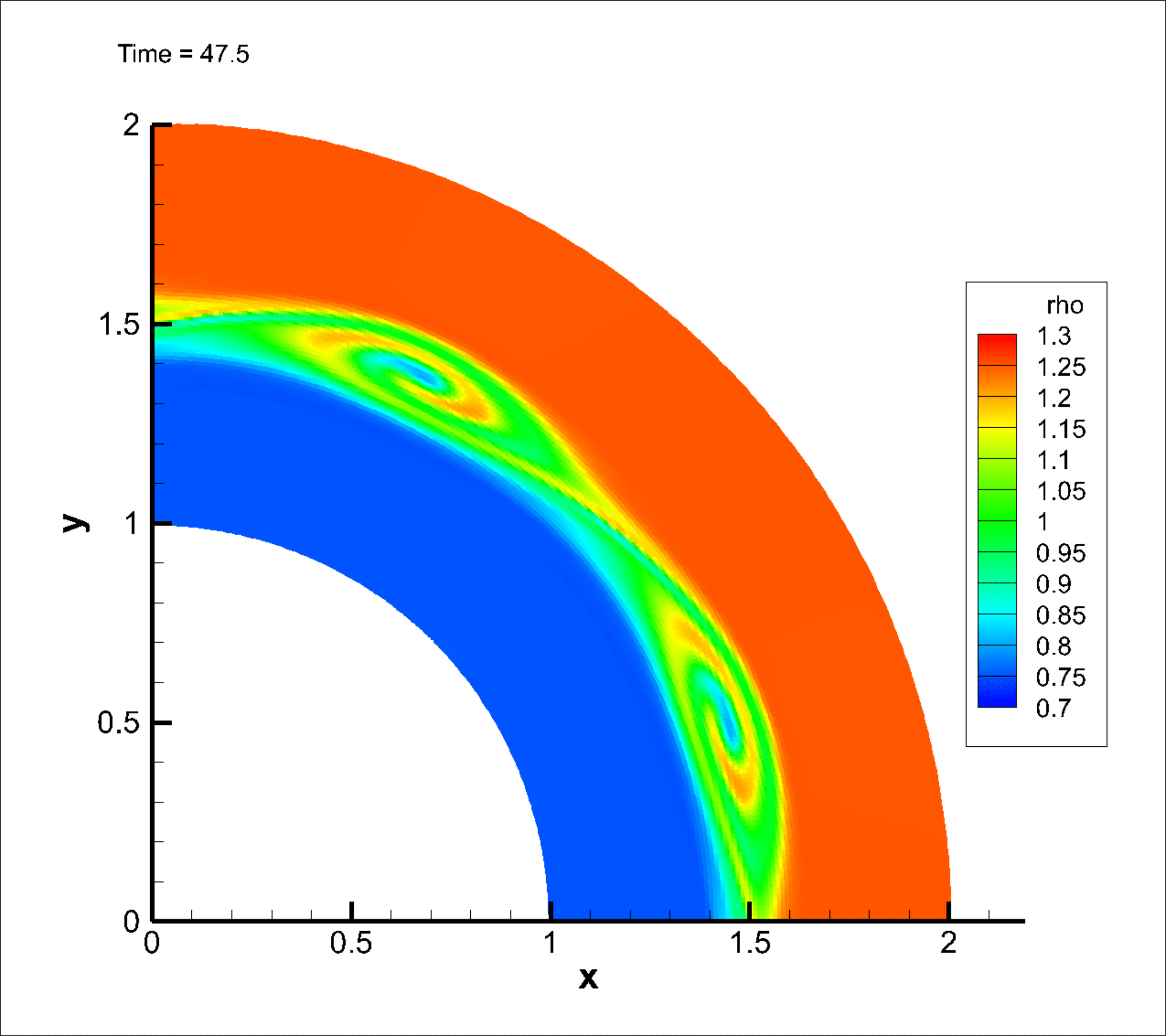} 		
	\includegraphics[width=0.24\linewidth]{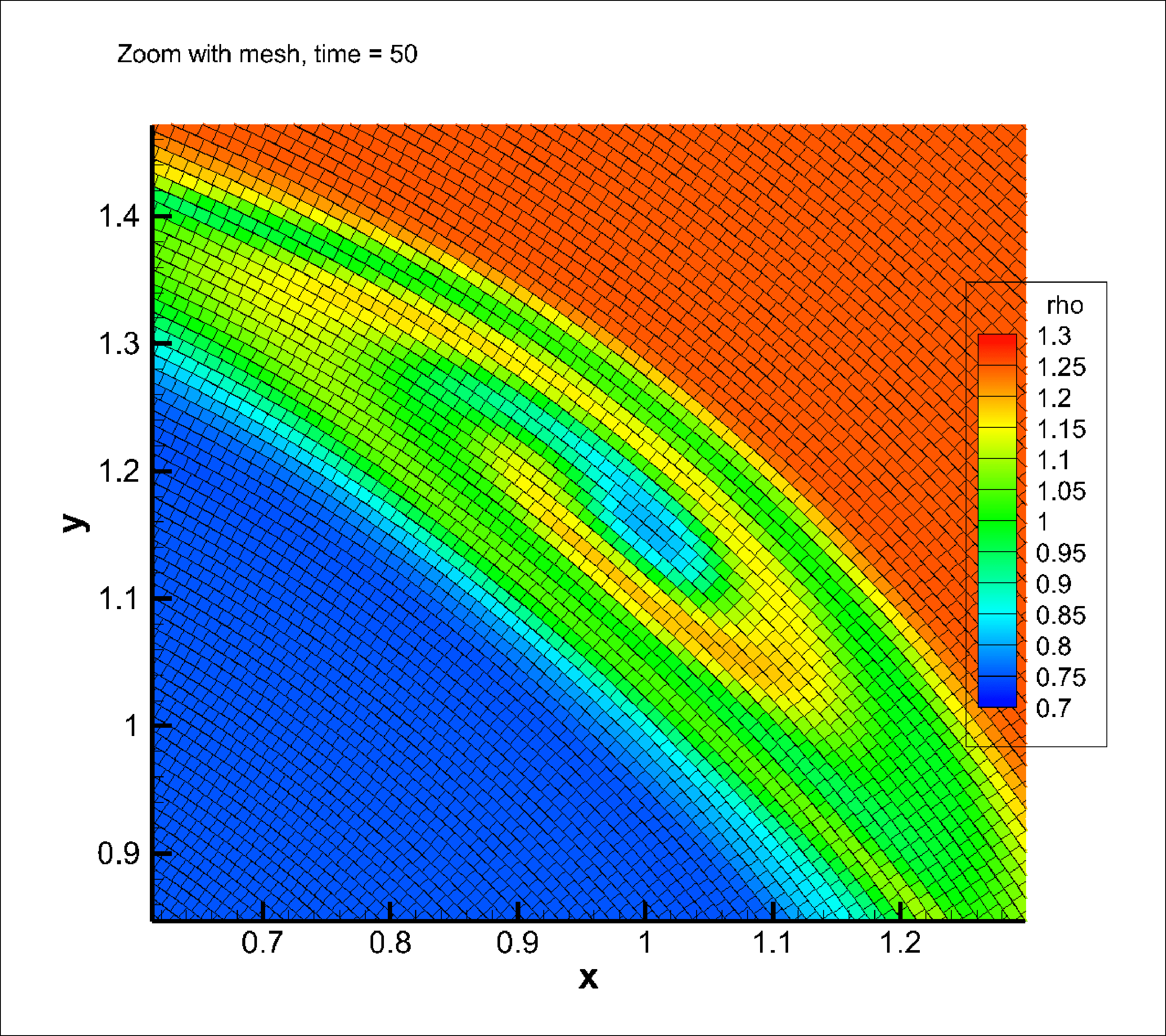}
	\caption{Kelvin-Helmholtz instabilities I. In the panel we show the evolution of the imposed periodic perturbations at different times. The results have been obtained with our second order well balanced ALE Osher-Romberg scheme over a grid with $100 \times 200$ control volumes. }
	\label{fig.KH_sector}
\end{figure*}

As confirmed by the previous tests, our well balanced ALE scheme is able to maintain the equilibrium up to machine precision for very long computational times.
So we can study with high accuracy the evolution of perturbations added to the density, the radial velocities and the pressure prescribed by the following initial condition 
\be
\begin{cases} 
\rho = \rho_E + A \rho_0  \sin(k \varphi) \text{exp} \left( - \frac{(r - r_m)^2}{ s}   \right ), \\
u = u_E +  A \sin(k \varphi) \text{exp} \left( - \frac{(r - r_m)^2}{ s }   \right ), \quad  
v = v_E, \\
P = P_E + A \sin(k \varphi) \text{exp} \left( - \frac{(r - r_m)^2}{s}   \right ), \\
\end{cases} 
\ee
with $ A= 0.1$, $k=8$, $ s = 0.005$. The computational results are depicted in Figure~\ref{fig.KH_sector}. In particular, for this flow configuration with physical perturbations in all flow quantities we observe the appearance of Kelvin-Helmholtz  instabilities and a very good resolution of the developing vortices, which is achieved thanks to the ALE technique and despite the rather coarse mesh of $100 \times 200$ elements used here. 

\begin{figure*}	
\textbf{ALE-WB} \qquad \qquad \qquad \qquad \qquad \qquad \qquad \textbf{EUL-WB} \qquad \qquad \qquad \qquad \qquad \qquad \qquad \textbf{ALE - noWB} \\[3pt]	
\includegraphics[width=0.3062\linewidth]{KH_sector_t175}  \quad \
\includegraphics[width=0.3062\linewidth]{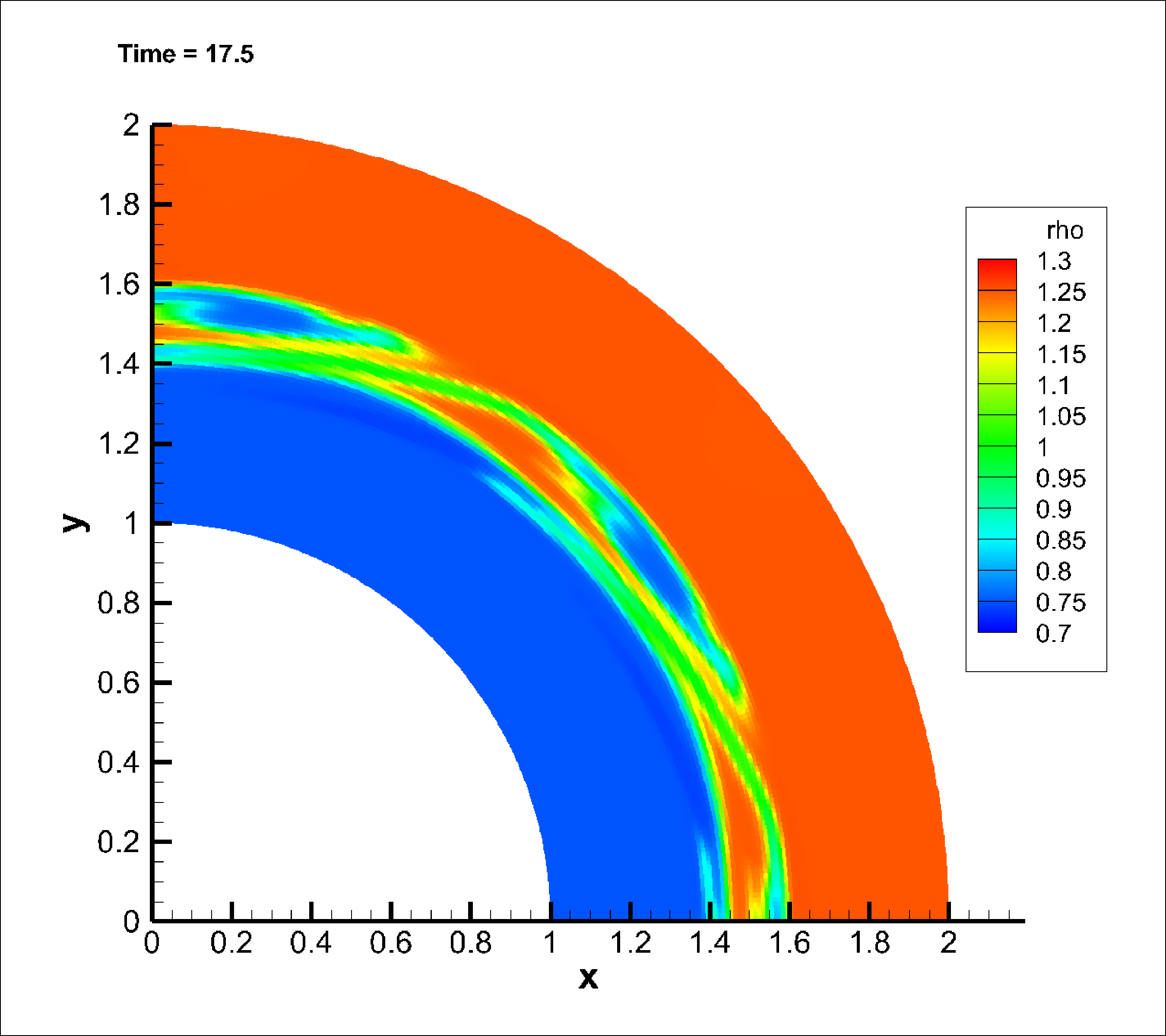} \quad \
\includegraphics[width=0.3062\linewidth]{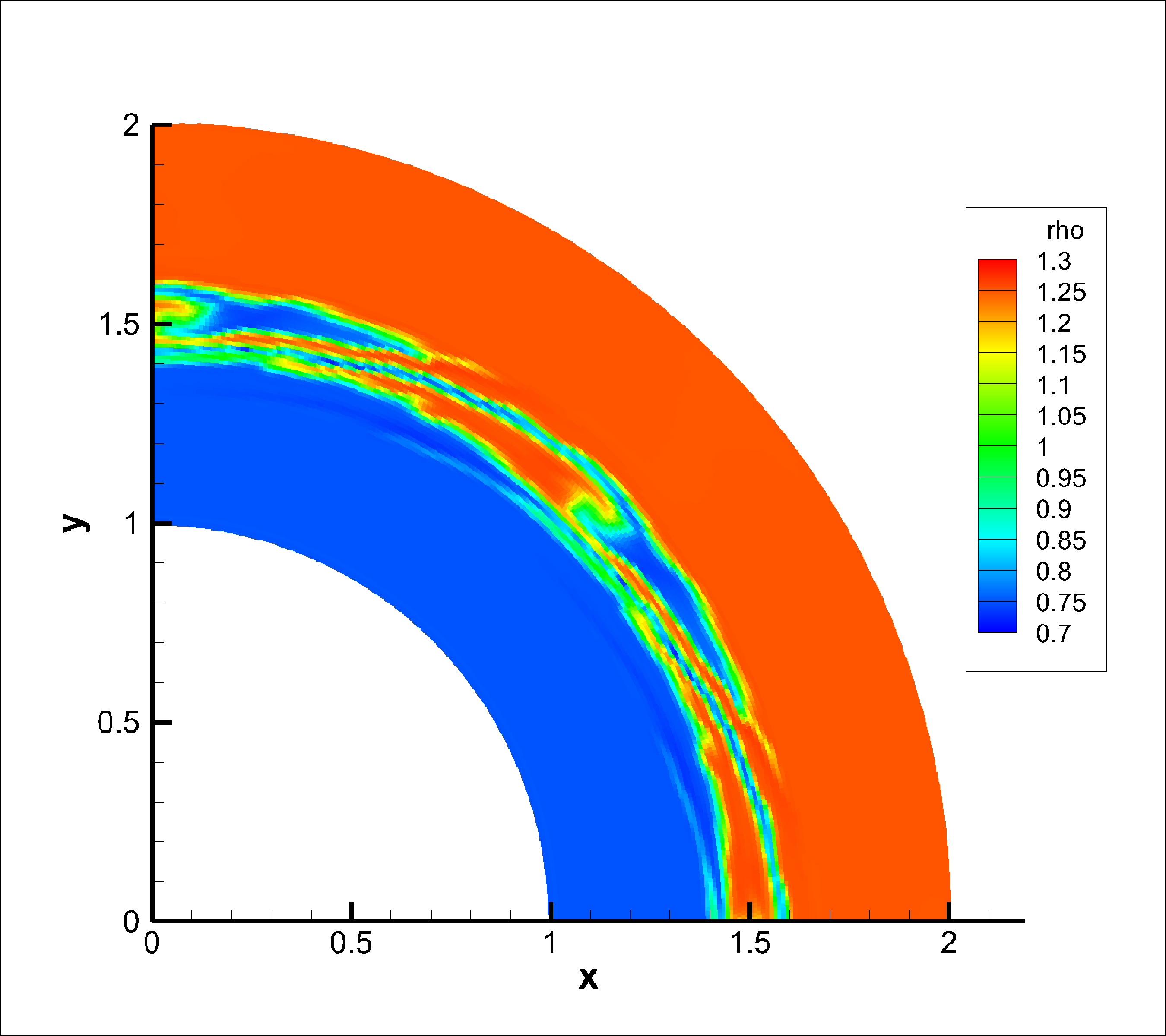} \\[3pt]
\ \includegraphics[width=0.3062\linewidth]{KH_sector_t250} \quad \
\includegraphics[width=0.3062\linewidth]{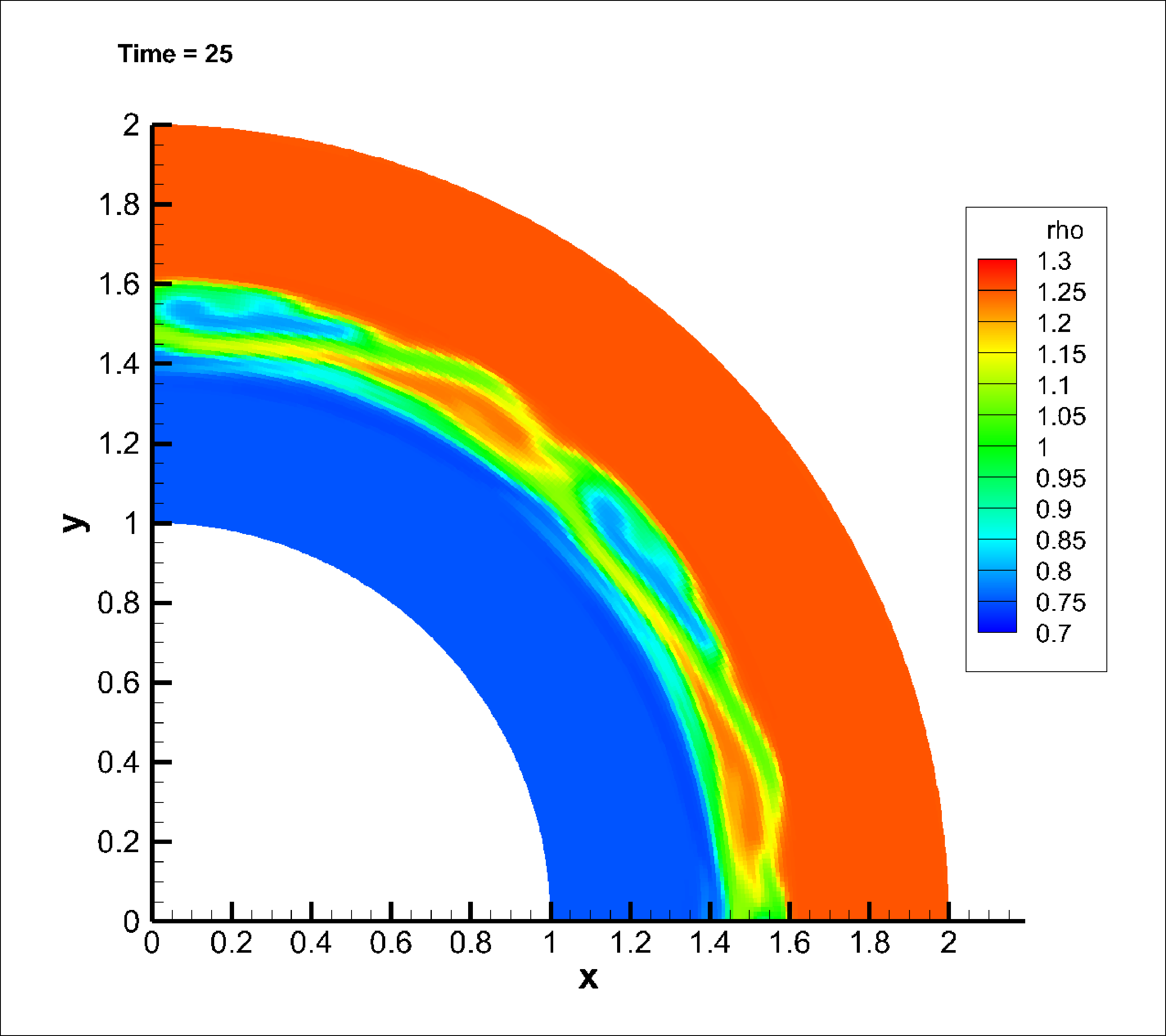}  \quad \
\includegraphics[width=0.3062\linewidth]{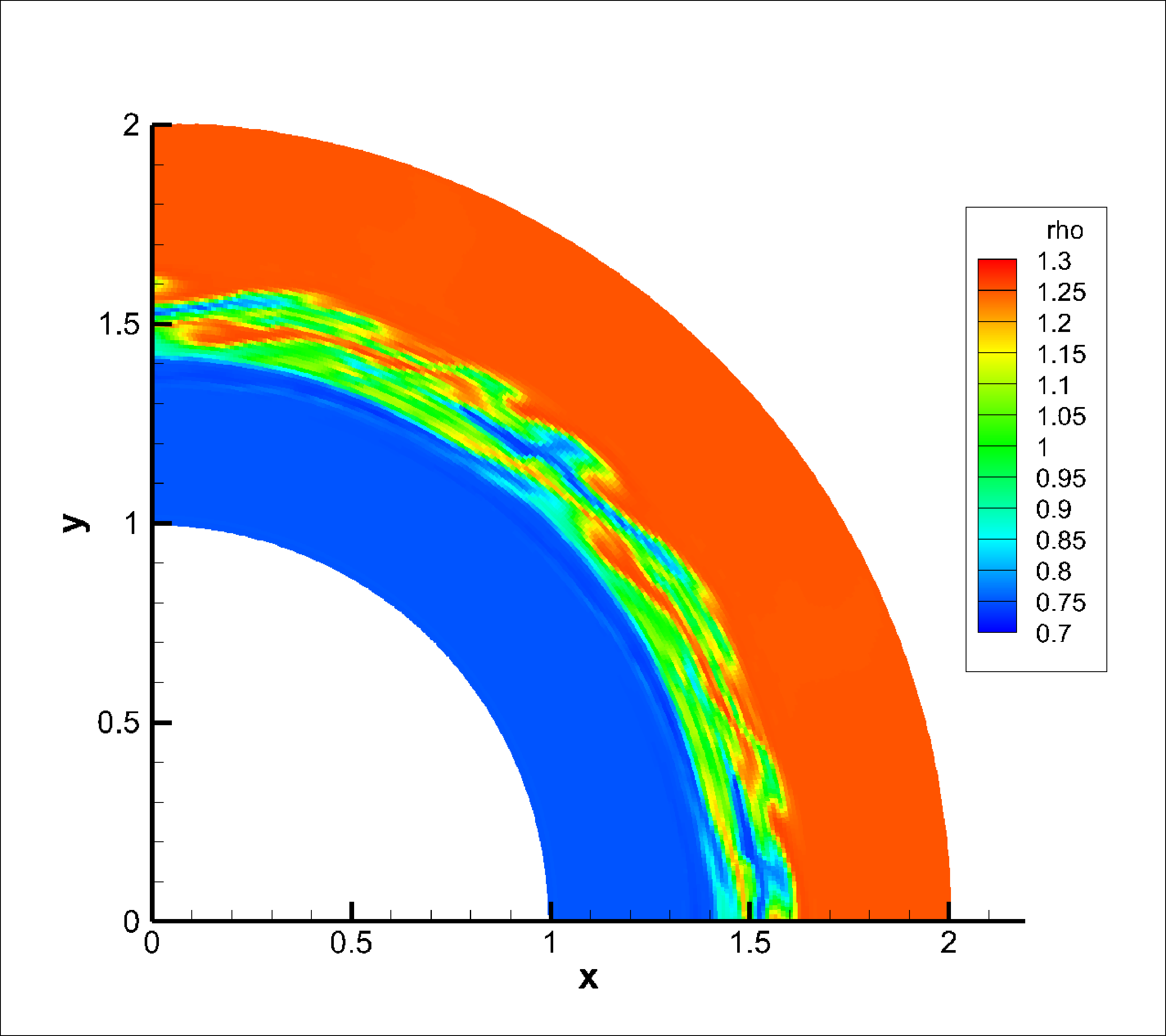} 
	\caption{Kelvin-Helmholtz instabilities I. In the panel we show the obtained solution for the density profile at time $t=17.5$ (first row) and time $t=25$ (second row). The results presented in the first column have been obtained using the Osher-Romberg well balanced ALE scheme. The ones in the second column have been obtained using a zero velocity mesh (Eulerian case) and the well balanced Osher-Romberg scheme. The third column is obtained with a standard nonconforming ALE scheme (i.e. using the ALE Osher type flux without well balancing). One can apreciate that it is really the coupling between the ALE and the well balancing that allows to achieve this high resolution. }
		\label{fig.comparisonMethosforKH}
\end{figure*}

Moreover we have compared our well balanced ALE scheme with a well balanced \textit{Eulerian} method on a fixed grid, which appears to be quite diffusive, and a \textit{not} well balanced ALE scheme, which produces visible spurious oscillations in the density profile. The results are presented in Figure \ref{fig.comparisonMethosforKH} and, once again, they show that it is indeed the coupling between the well balanced techniques and the moving mesh framework that allows to achieve a high resolution on small perturbations around an equilibrium solution for very long computational times. 

\begin{figure*}
	\includegraphics[width=0.24\linewidth]{KH_sector_t375} 	
	\includegraphics[width=0.24\linewidth]{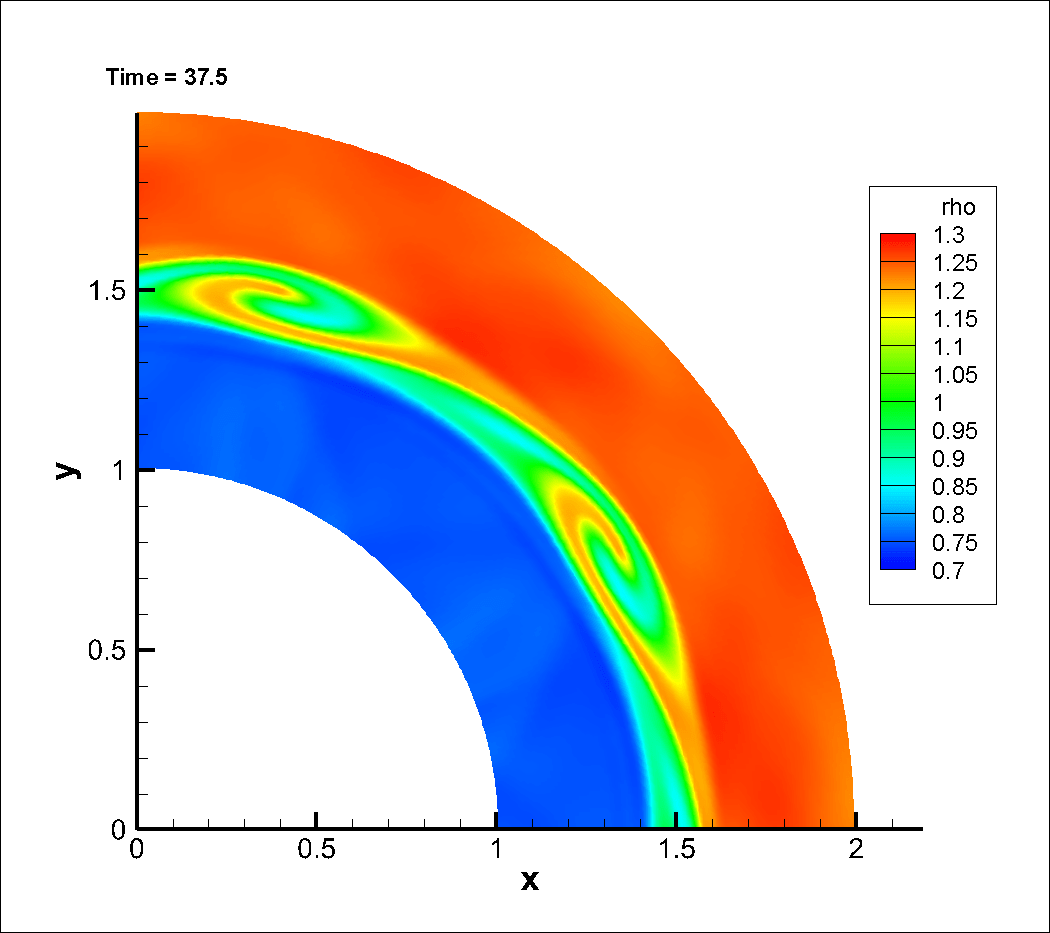} 
		\includegraphics[width=0.24\linewidth]{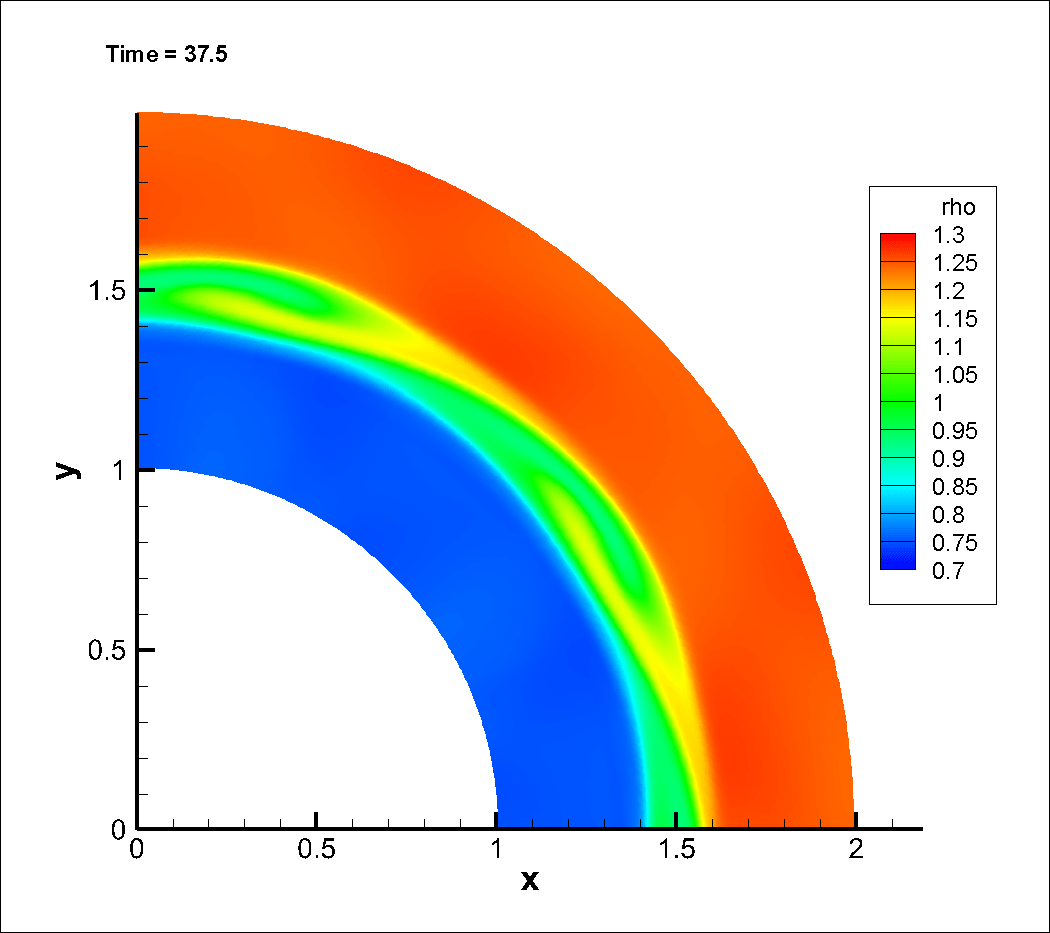}
	\includegraphics[width=0.24\linewidth]{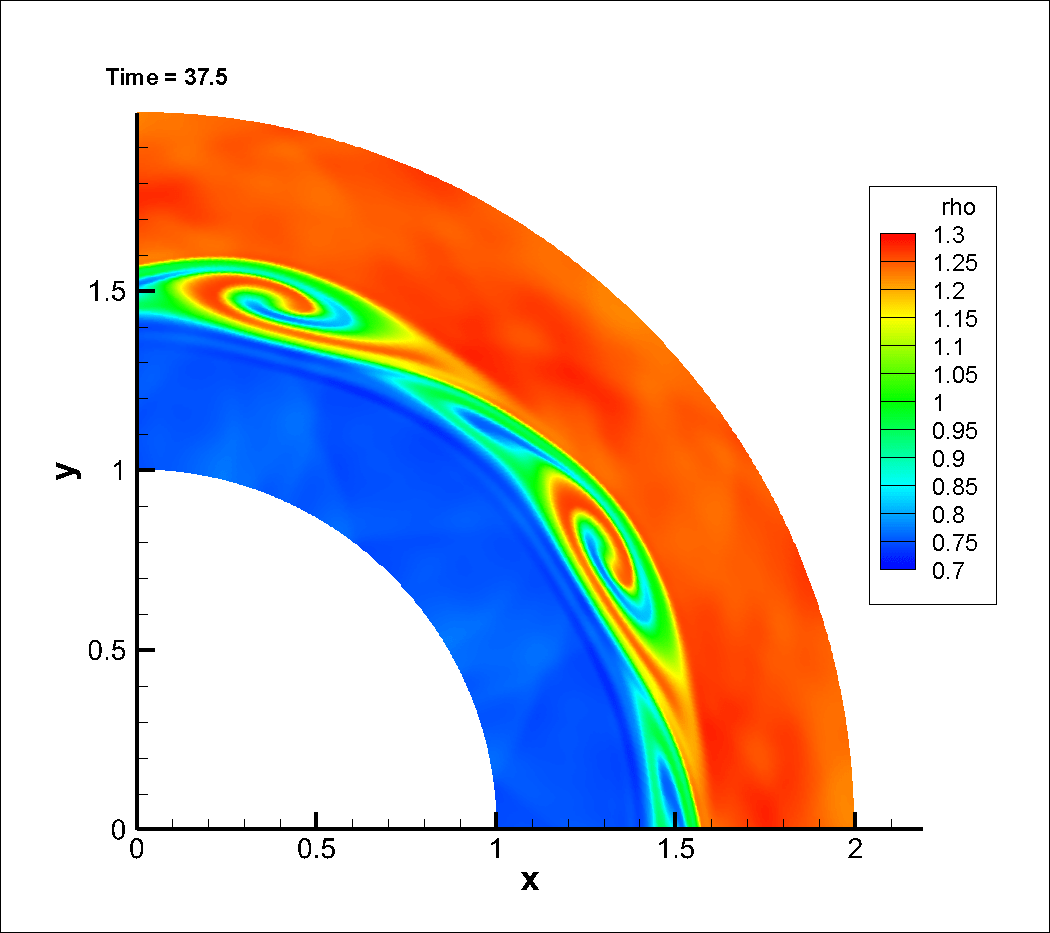} 
	\caption{{Kelvin-Helmholtz instabilities I. Method comparison at time $t=37.5$. The first image is obtained with our code and $100x200$ elements. The second and the third one with PLUTO, 
	$100\times200$ elements and respectively a second order scheme with mc\_lim limiter and a third order scheme and minmod\_lim limiter. The last image is obtained with the third order version of PLUTO  using mc\_lim and $200\times400$ elements. All images are drawn with the same color map.}}
	\label{fig.KeplerianDisk_PLUTO_KH}
\end{figure*}

{We also compare our numerical results at time $t=37.5$ with those obtained by PLUTO, see Figure \ref{fig.KeplerianDisk_PLUTO_KH}. 
In order to obtain the same accuracy of our new second order well balanced Osher Romberg ALE scheme (left image of the panel) one needs the third order version 
of PLUTO on a finer mesh (last image of the panel).}

\subsection{Keplerian disc with Kelvin-Helmholtz instabilities II}
\begin{figure*}	
	\includegraphics[width=0.24\linewidth]{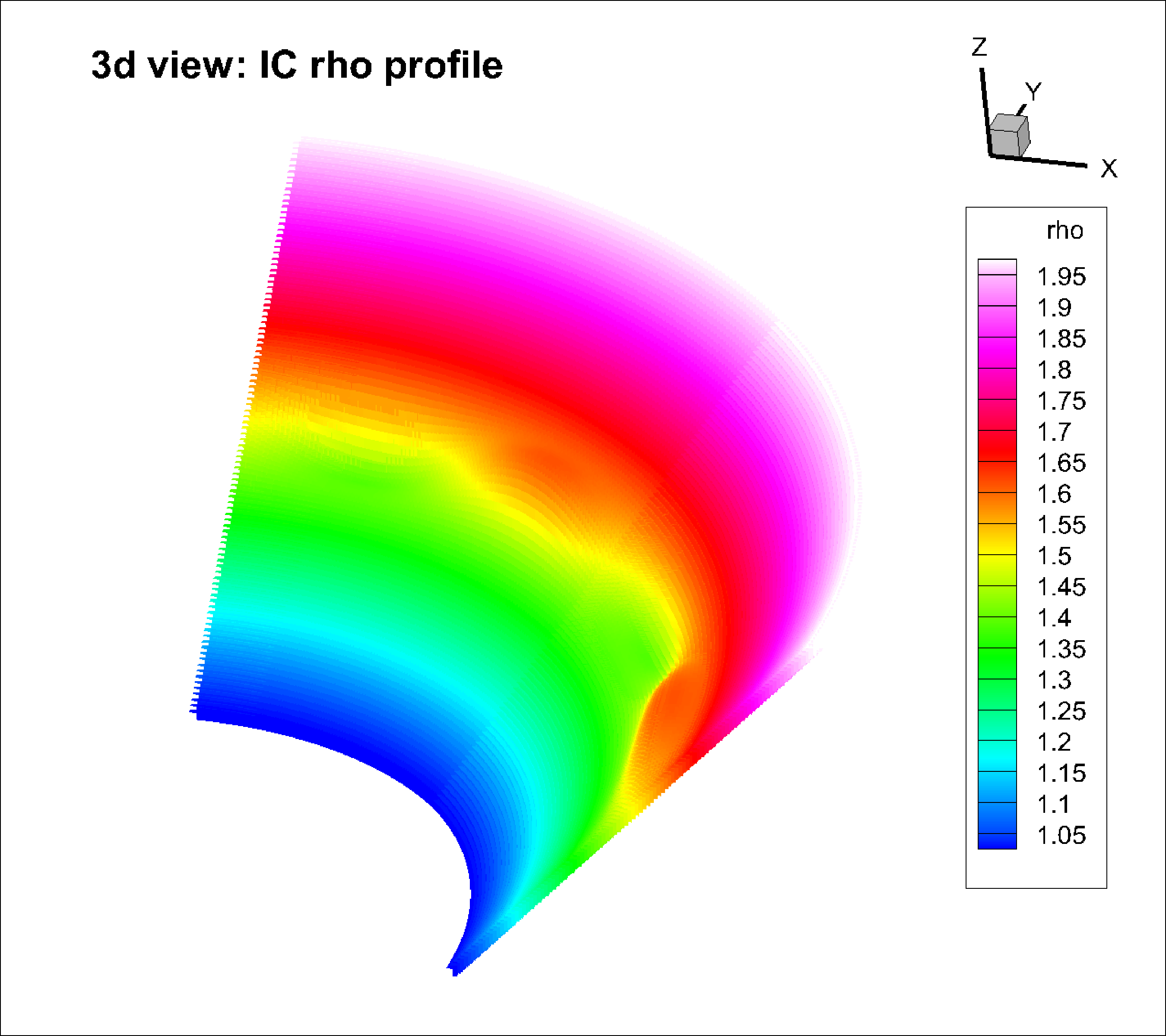} 		
	\includegraphics[width=0.24\linewidth]{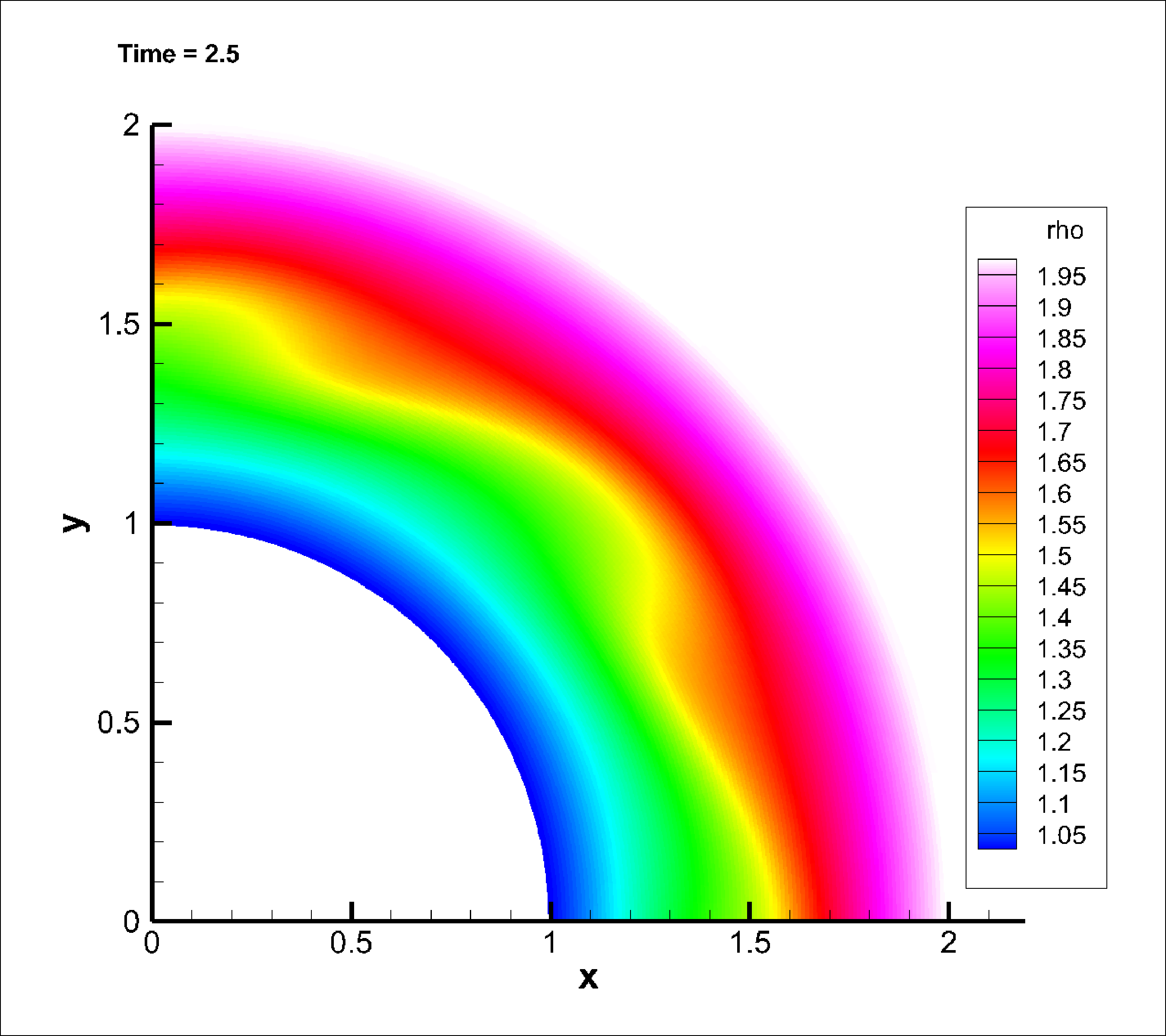} 		
	\includegraphics[width=0.24\linewidth]{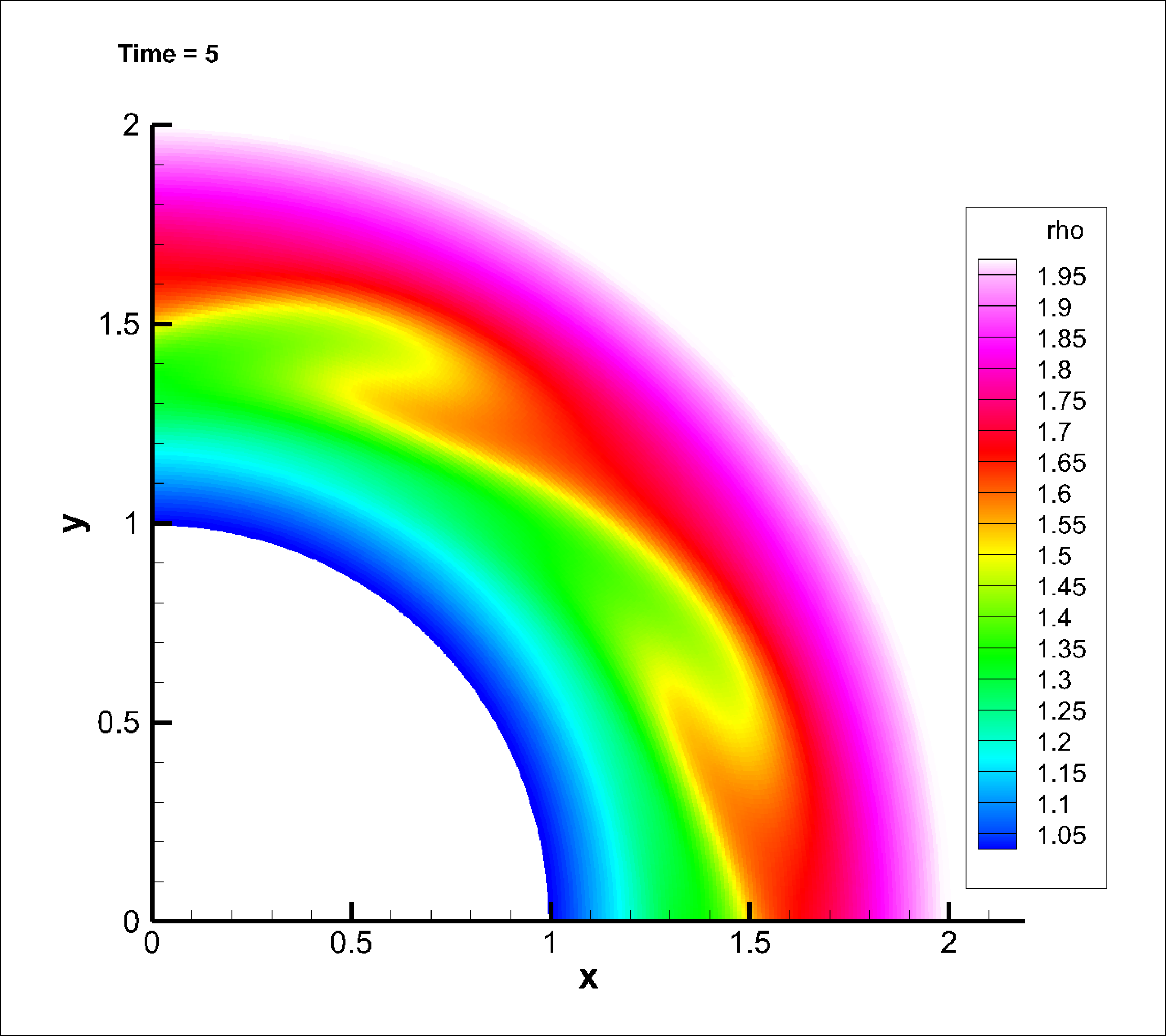} 		
	\includegraphics[width=0.24\linewidth]{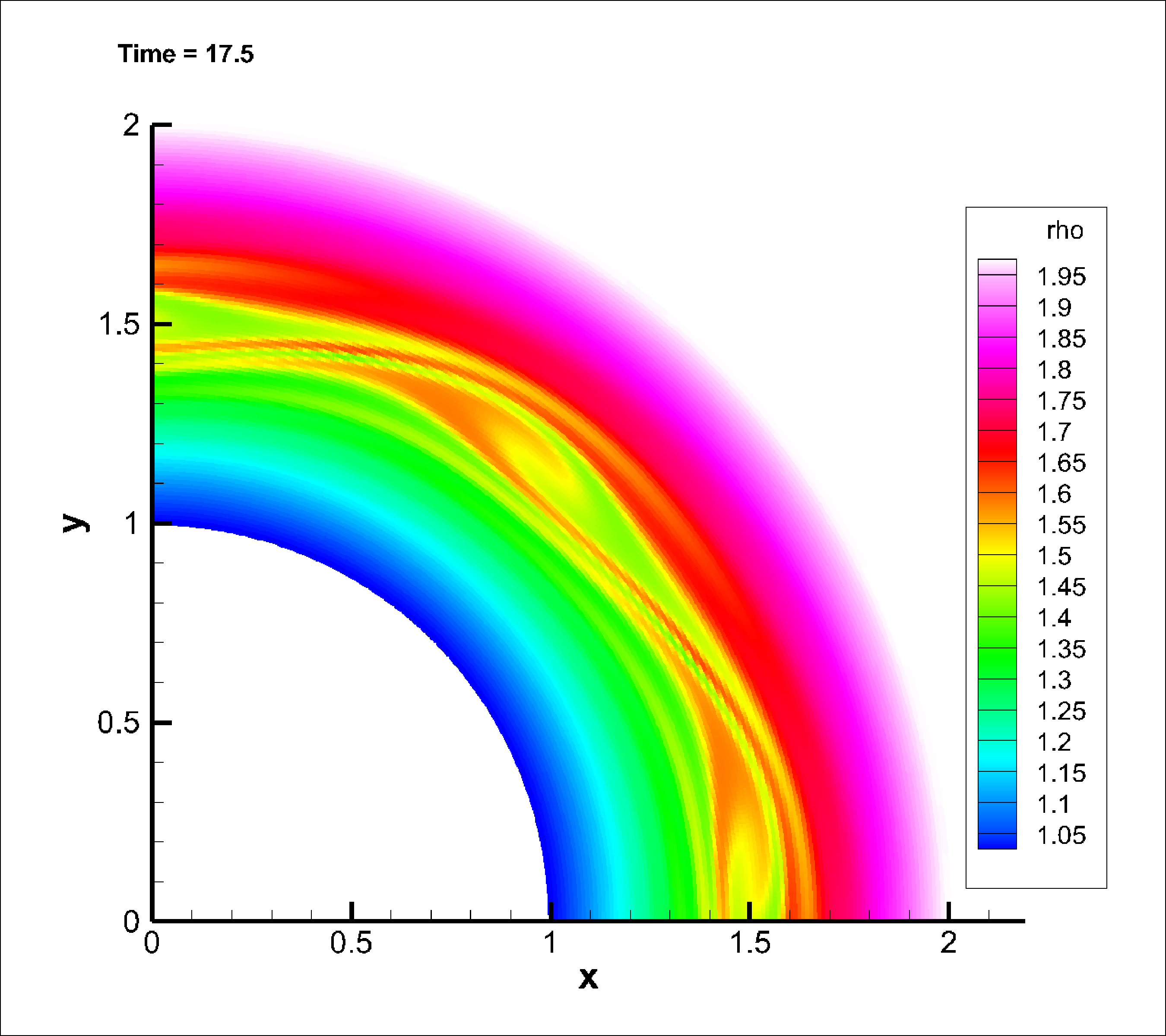} 	 \\
	\ \includegraphics[width=0.24\linewidth]{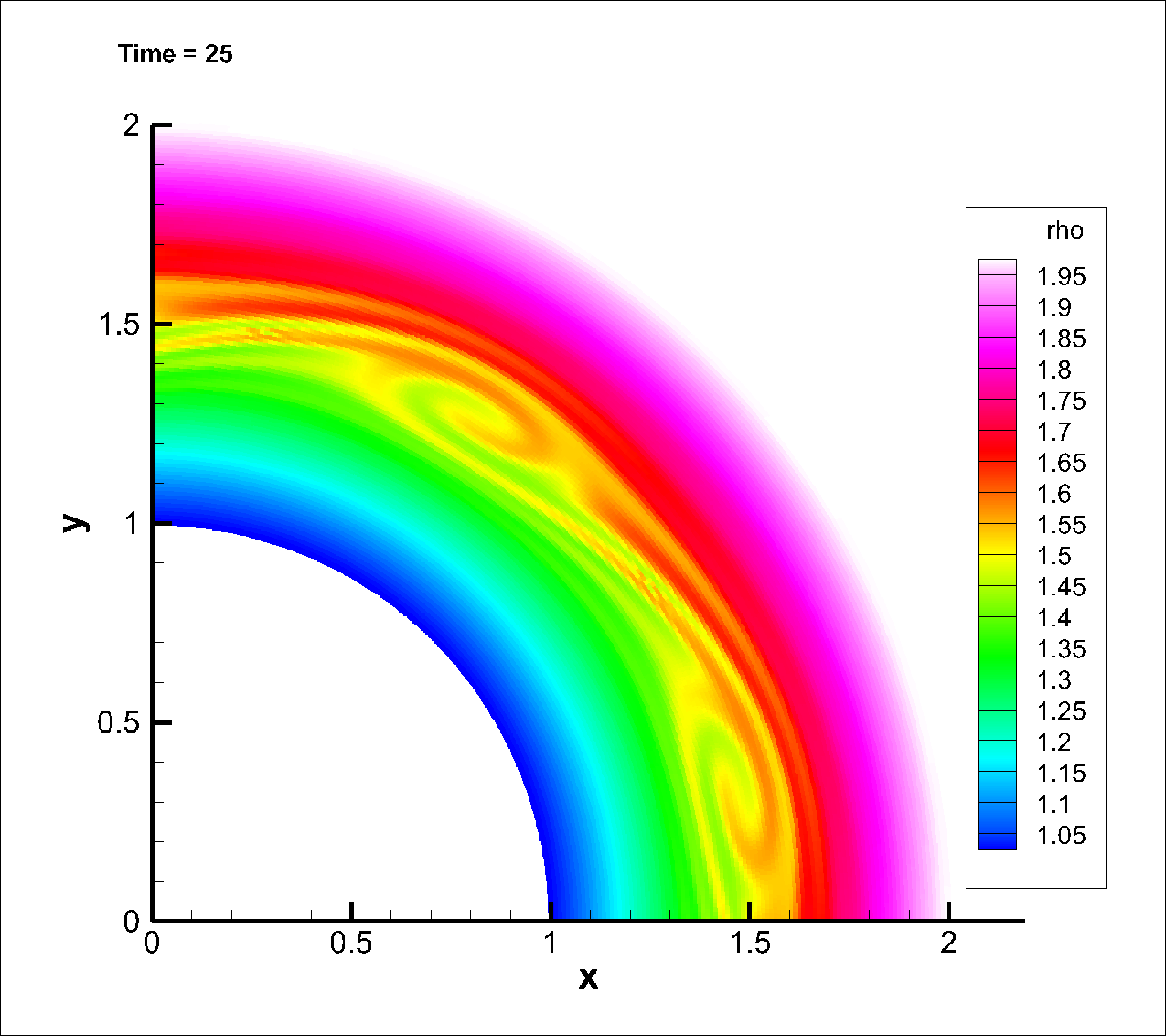}
	\includegraphics[width=0.24\linewidth]{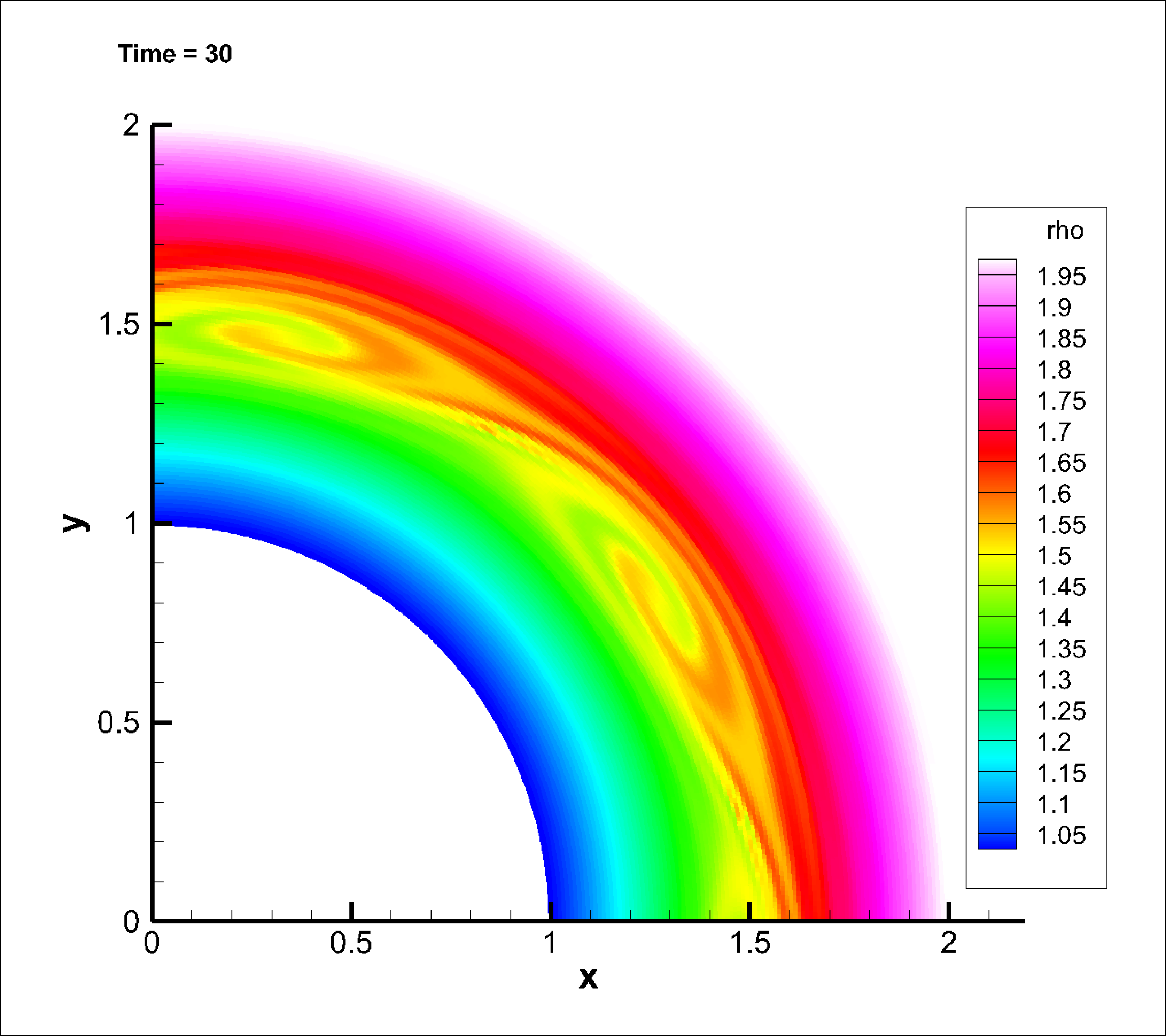}
	\includegraphics[width=0.24\linewidth]{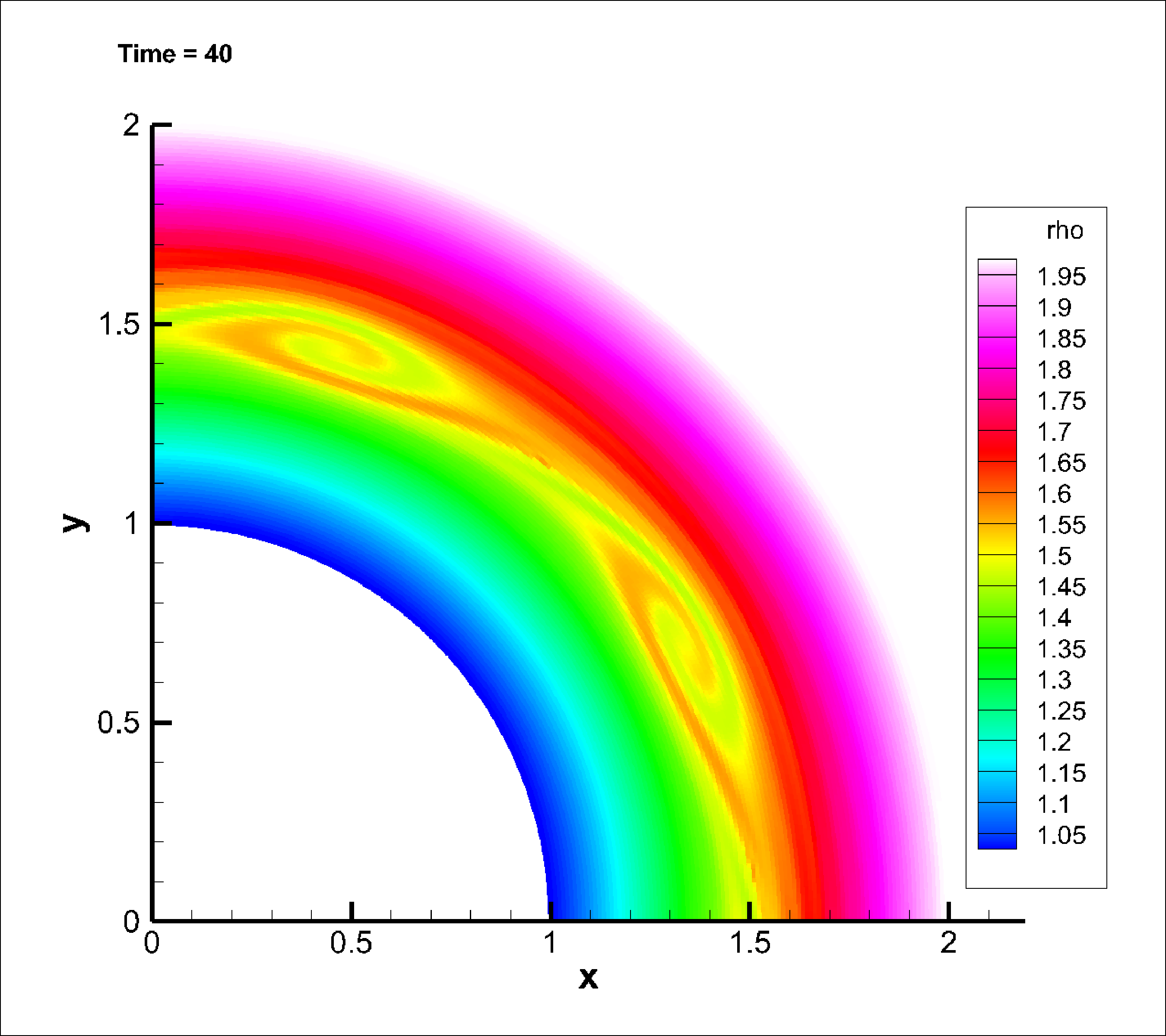} 		
	\includegraphics[width=0.24\linewidth]{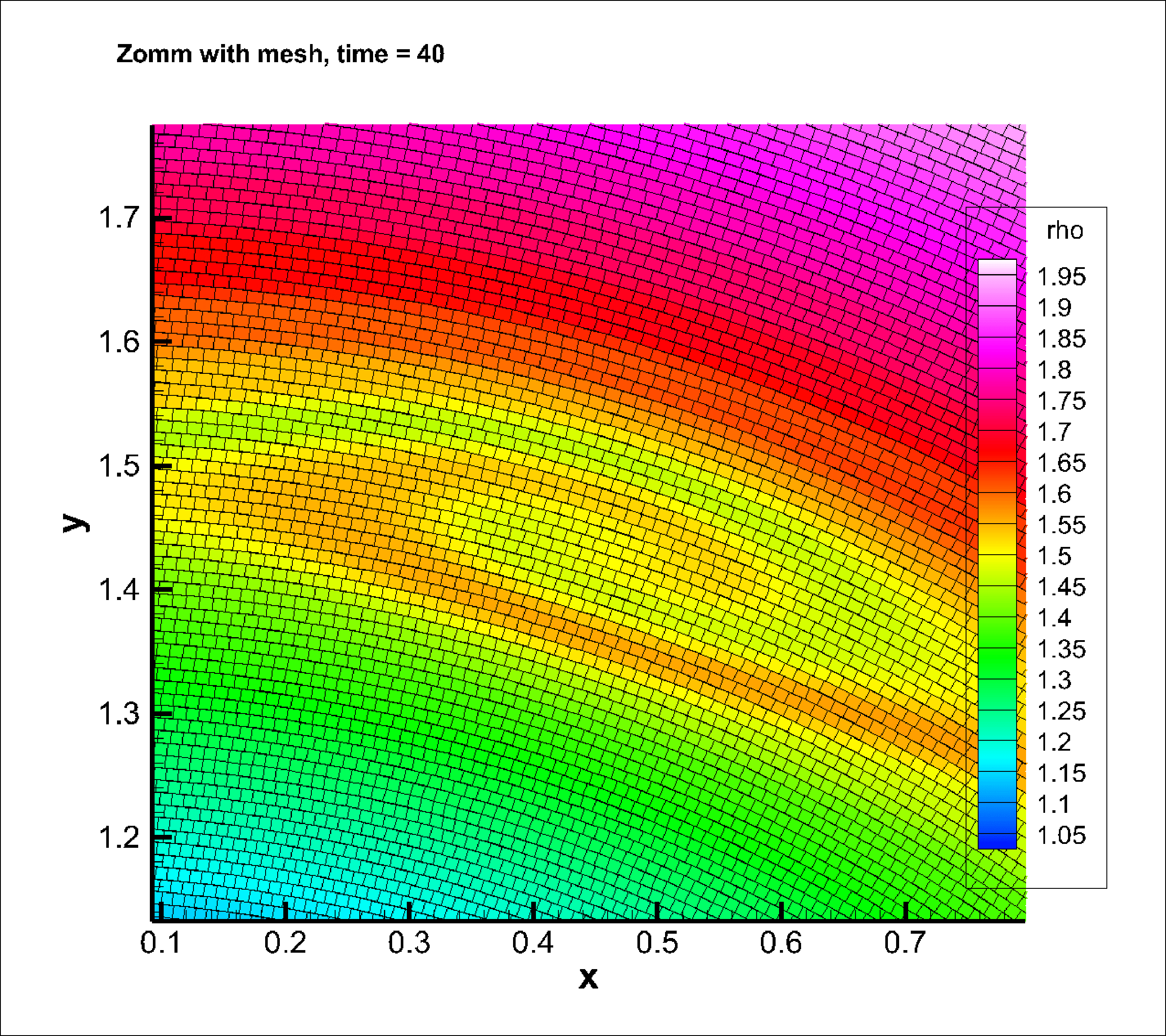}
	\caption{Kelvin-Helmholtz instabilities II. In the panel we show the evolution of the imposed periodic perturbations at different times. The results have been obtained with our second order Osher-Romberg scheme over a grid with $100 \times 200$ control volumes. }
	\label{fig.KH_sector2}
\end{figure*}

We finally consider another equilibrium solution which satisfies the equilibrium constraints in \eqref{eq.EquilibriaConstraint}-\eqref{eq.PressureAndGravForces} and which reads 
\be
	\label{eq.equilibriumForKH2}
	\rho_E = r,  \quad 
	u_E = 0, \quad 
	v_E = \sqrt{ \frac{Gm_s}{r} }, \quad 
	P_E = 1,  
\ee
with $G = 1$, $m_s = 1$ and $r_m= 1.5$. 
With respect to the previous example, here the density profile is linear. However, also in this example we expect the Kelvin-Helmholtz instabilities to arise if some perturbations are added to the  stationary profile.
The computational domain and the boundary conditions are chosen as before. 
The initial condition used in this test problem reads 
\be
\begin{cases} 
	\rho = \rho_E + A   \sin(k \varphi) \text{exp} \left( - \frac{(r - r_m)^2}{ s}   \right ), \\
	u = u_E +  A \sin(k \varphi) \text{exp} \left( - \frac{(r - r_m)^2}{ s }   \right ), \quad 
	v = v_E, \\
	P = P_E + A \sin(k \varphi) \text{exp} \left( - \frac{(r - r_m)^2}{s}   \right ), \\
\end{cases} 
\ee
with $ A= 0.1$, $k=8$, $ s = 0.005$, i.e. we are again solving a problem that is close to an equilibrium and therefore difficult to solve with standard numerical techniques that are not well balanced. 
The computational results are depicted in Figure~\ref{fig.KH_sector2}. Again we observe the appearance of Kelvin-Helmholtz instabilities that are well resolved also on a rather coarse mesh, without any  visible spurious numerical oscillations. 

\begin{figure*}
\includegraphics[width=0.24\linewidth]{KH_sector_LinearRHO_t250}
\includegraphics[width=0.24\linewidth]{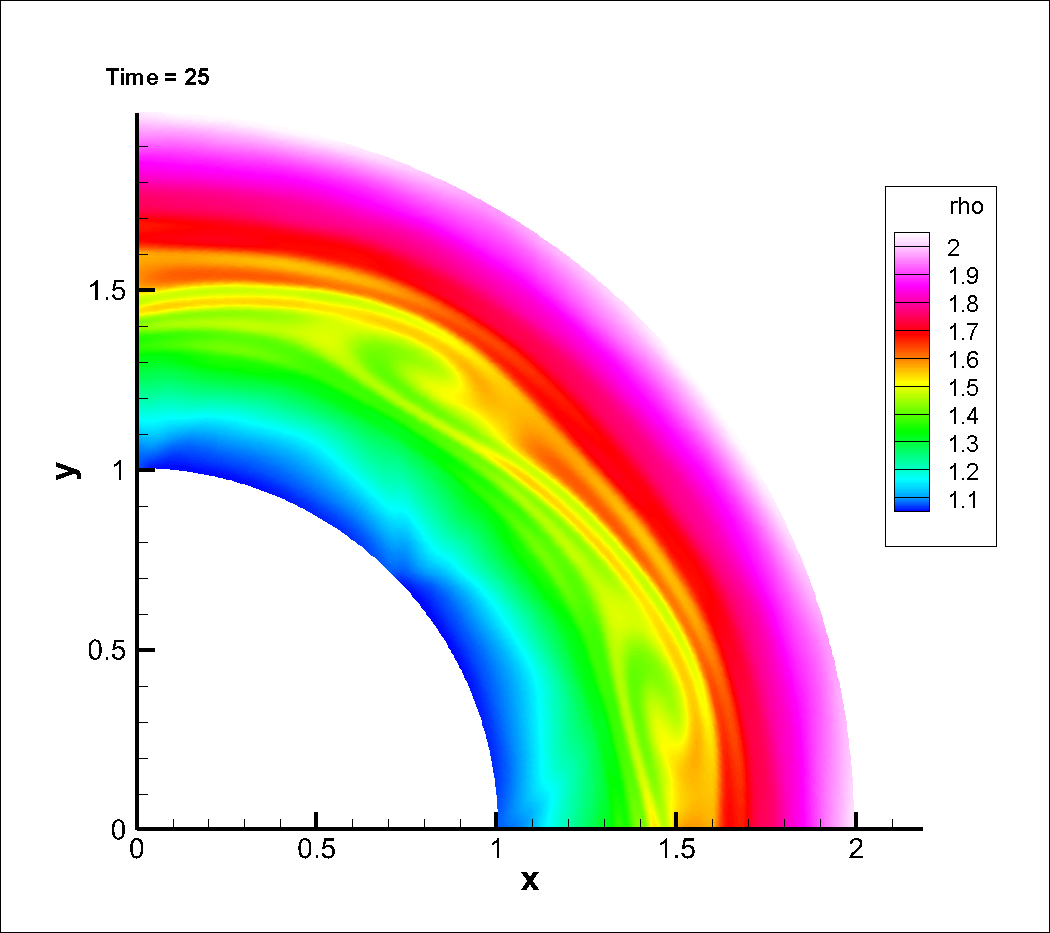} 
\includegraphics[width=0.24\linewidth]{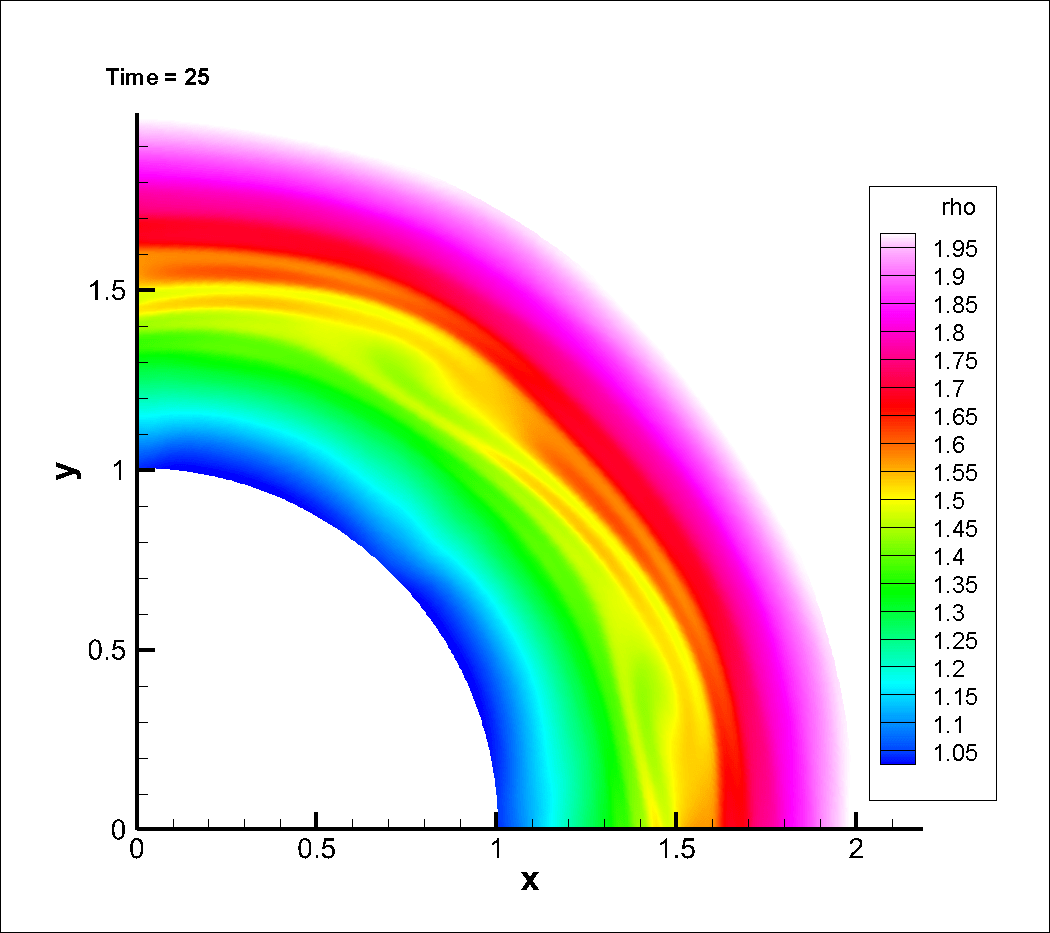} 
\includegraphics[width=0.24\linewidth]{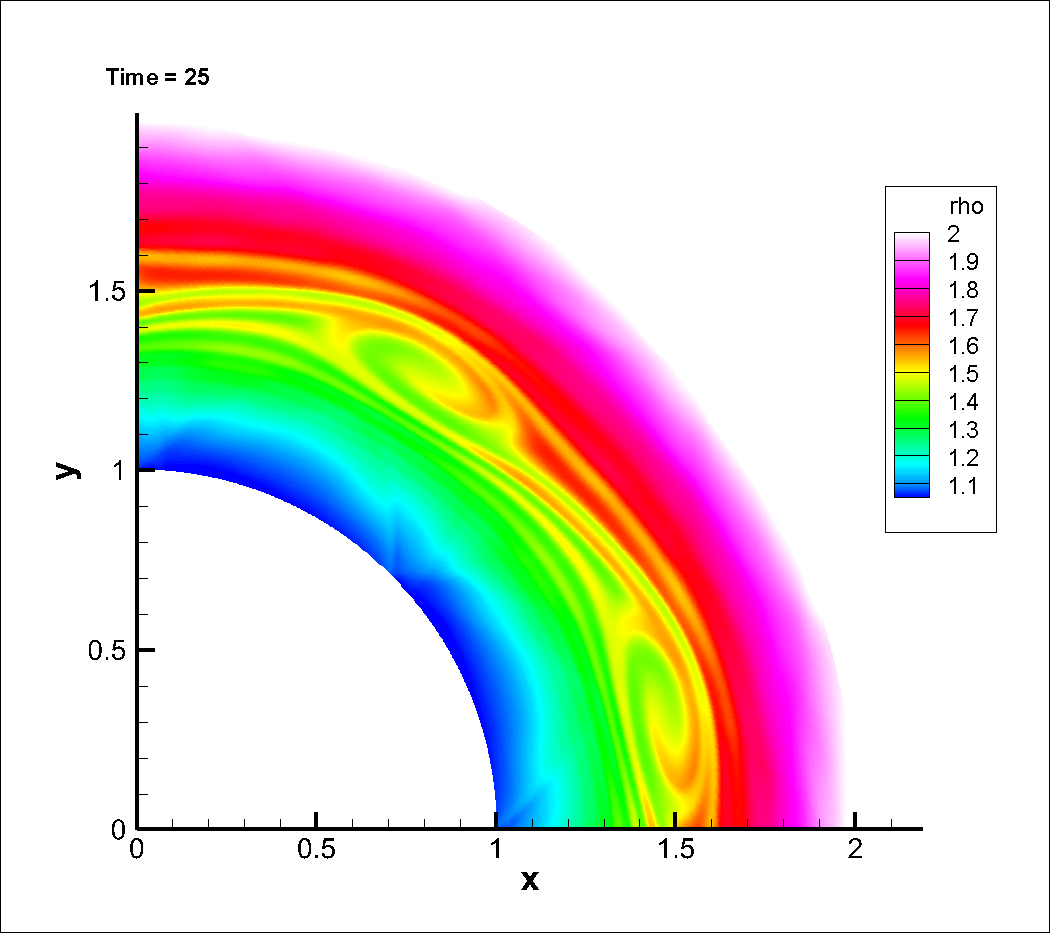} \\[2pt]
\ \includegraphics[width=0.24\linewidth]{KH_sector_LinearRHO_t400}
\includegraphics[width=0.24\linewidth]{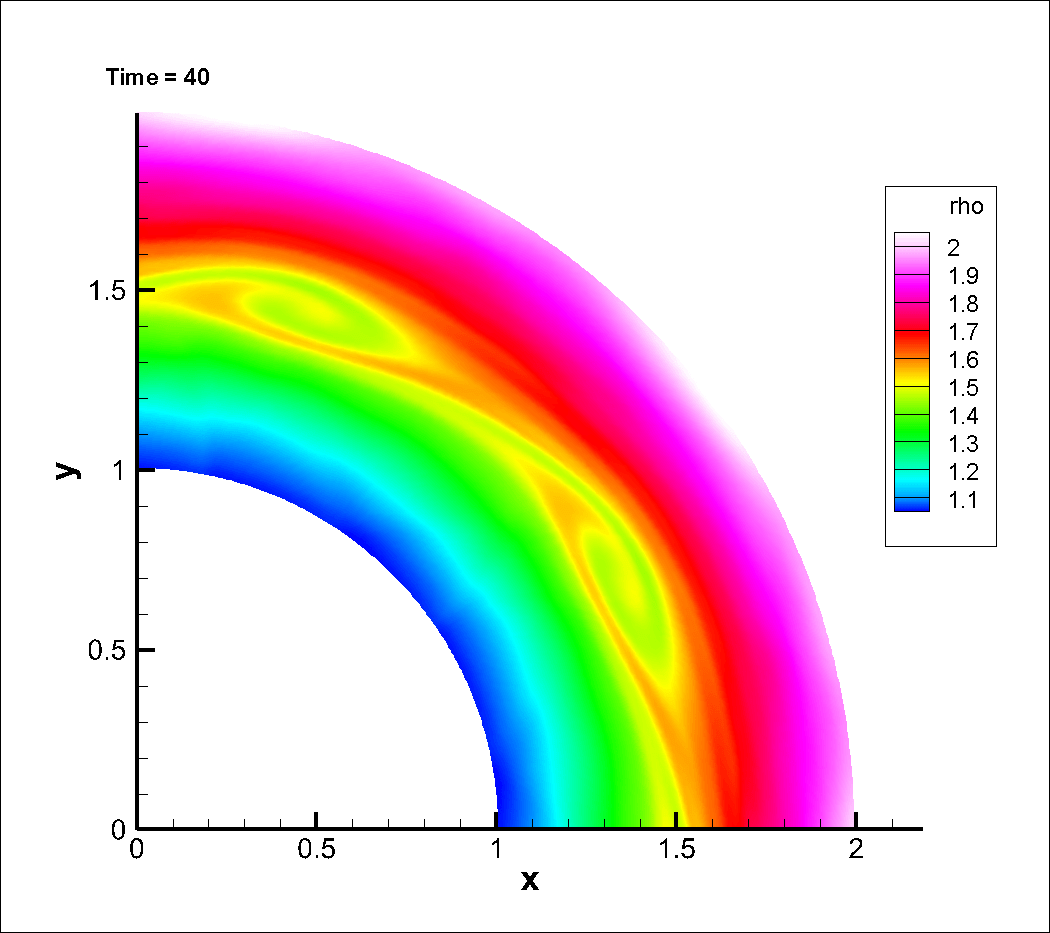} 
\includegraphics[width=0.24\linewidth]{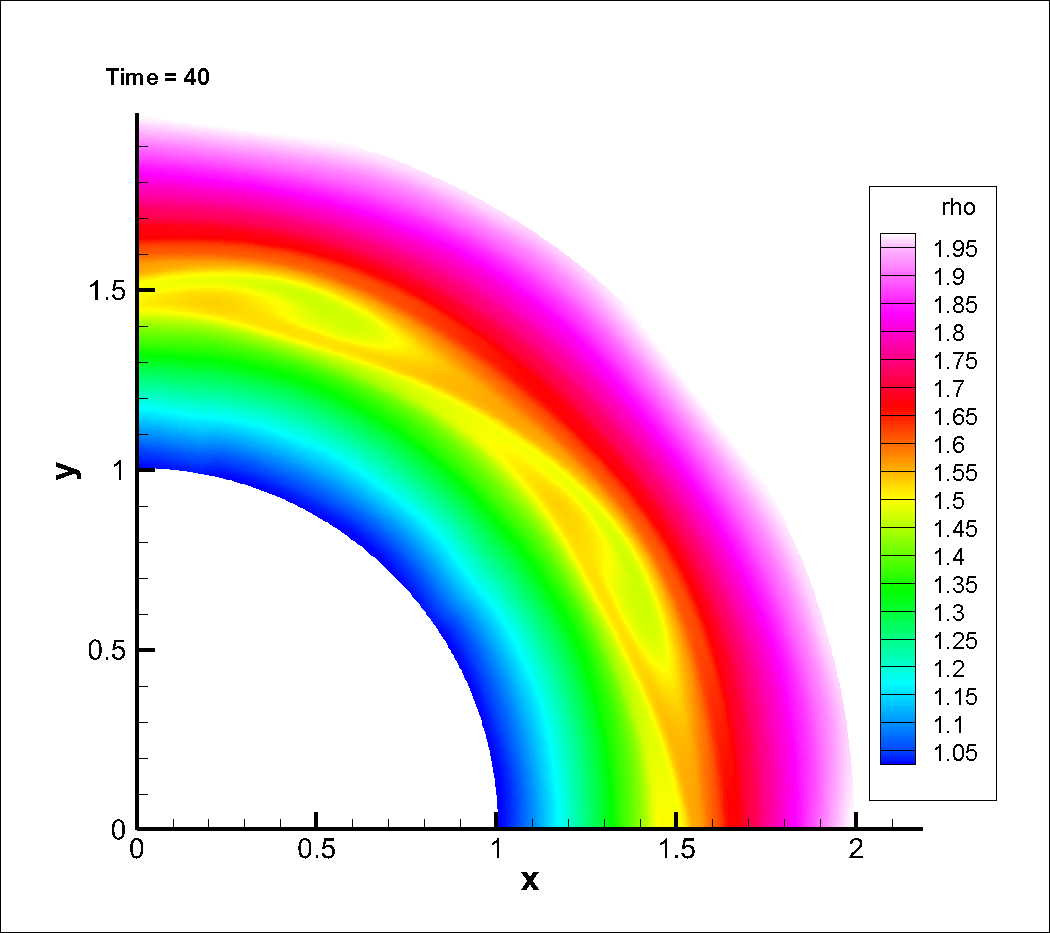} 
\includegraphics[width=0.24\linewidth]{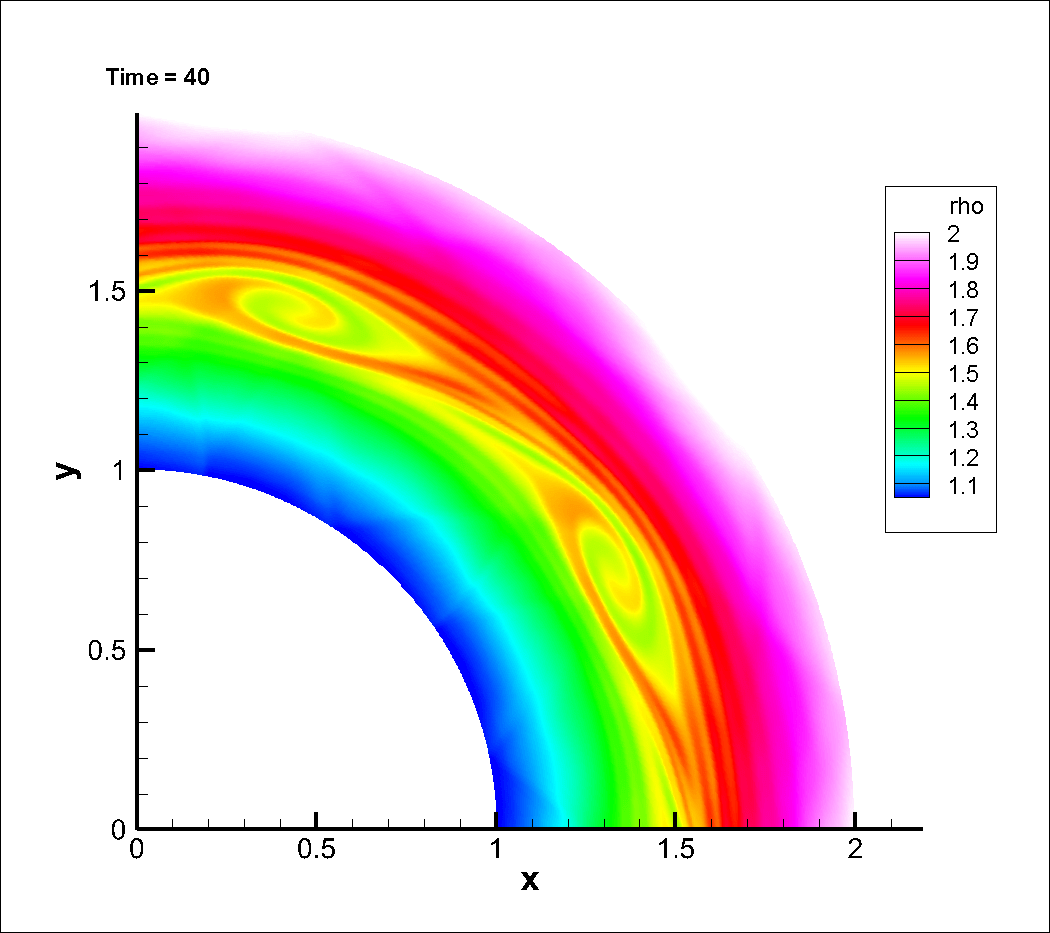} 
	\caption{{Kelvin-Helmholtz instabilities II. Method comparison at time $t=25$ (first row) and at time $t=40$ (second row). The first images are obtained with our code and $100x200$ elements. The second and the third ones with PLUTO, $100\times200$ elements and using respectively a second order (mc\_lim) and a third order (minmod\_lim) scheme. The last images are obtained with 
	PLUTO using a third order scheme (mc\_lim) and $200\times400$ elements. All images are drawn with the same color map. The vortices have a similar resolution in the leftmost and rightmost images.}} 
	\label{fig.KeplerianDisk_KHlinear_resultsPLUTO}
\end{figure*}

{Finally, we compare once again our code with results obtained with PLUTO, refer to Figure \ref{fig.KeplerianDisk_KHlinear_resultsPLUTO}. A similar resolution of the vortices is obtained with our second order code and the third order version of PLUTO with a finer mesh (refer to Section \ref{ssec.PLUTO_setting} for the details on the PLUTO configuration we have chosen). In this case we want to underline also that our code avoids other oscillations that instead can be noticed in the images obtained with PLUTO.}

\section{Conclusions}
In this work we have developed a new and highly accurate well balanced path-conservative finite volume scheme for the Euler equations with gravity by proposing two specifically designed numerical fluxes  and a quite general reconstruction procedure. 
We underline that the novelties introduced in the algorithm are based on the following \textit{key idea}: the construction of a path which directly exploits the known stationary solution (and so the scheme is exact on it), and treats in a approximate way only the fluctuations around the equilibrium.

The proposed method is innovative already in one space dimension, since to the knowledge of the authors, it is the first time that the little dissipative path-conservative Osher scheme proposed by \cite{OsherNC} is modified in order to be well balanced for non-trivial equilibria of the Euler equations of gas dynamics with gravity. In particular, the way in which the absolute value of the Jacobian matrix is  rewritten in order to exploit even in the viscosity part of the scheme the same well balanced strategy that was already established for the non-dissipative part is original. 
Moreover, it is the use of the Romberg quadrature formula (instead of the Gaussian one) that provides the good properties to maintain both the desired order of accuracy and the well balancing.

Furthermore, the method has been carefully extended in a non trivial way to the two dimensional framework preserving the well balancing even for moving domains (with only few constraints on the mesh construction). In particular the coupling with modern nonconforming ALE techniques enables the resolution even of complex shear flows with differential rotation in an effective way. At this point it is 
noteworthy to stress again that standard \textit{conforming} Lagrangian schemes will \textit{crash} after \textit{finite times} for any vortex flow with differential rotation due to \textit{mesh tangling}. 
Indeed the reduced dissipation characterizing the Lagrangian methods, together with the high mesh quality provided by the nonconforming treatment of sliding lines, and the increased accuracy near the equilibria given by the well balanced techniques, allow us to obtain significant improvements compared to the existing state of the art. 
The major benefits are achieved with our new class of schemes when studying physical phenomena that arise close to a stationary equilibrium solution, where standard discretizations would hide the 
flow physics by spurious oscillations and excessive numerical dissipation. 

{We furthermore have provided a thorough comparison of our new numerical method with the results that can be obtained with the PLUTO code, which is based on finite volumes and therefore 
is rather close to the scheme proposed in this paper.} 

{While the moving nonconforming mesh treatment proposed here is rather invasive and probably quite difficult to introduce in existing astrophysical codes, our new path-conservative finite  volume scheme that achieves the well balancing of the method at the level of the Riemann solver is instead \textit{straightforward} to implement in existing schemes and computer codes based on Riemann solvers,  i.e. those using classical first or second order Godunov-type finite volume and finite difference methods. All that is needed is to replace the conventional algebraic source term by our new well balanced path-conservative approximate Riemann solver, which interprets the gravity source term as a nonconservative product. Nevertheless, using the novel ideas on well balanced SPH methods very recently  presented in \cite{PCSPH}, it seems also possible to extend the new well balanced approach for the Euler equations with gravity presented here to Smooth Particle Hydrodynamics. However, this is beyond the scope of the present paper and its feasibility will be subject to further investigations. } 

Future research will consider the application to more complex systems of hyperbolic PDE, such as the unified model of continuum mechanics presented in \cite{PeshRom2014,HPRmodel,HPRmodelMHD}, an extension to three space dimensions as well as to more general classes of stationary solutions and an automatic detector of the equilibrium profiles in order to extend our method to situations in which the equilibrium is not known exactly \apriori. {Based on the high order path-conservative methods introduced in \cite{ADERGRMHD} we also plan to use the  algorithms developed in this paper in order to design exactly well balanced schemes for gravity driven equilibrium flows in \textit{general relativity}, where the use of well balanced methods appears to be still rather unknown. } 
We also plan to extend the presented method to better than second order of accuracy by extending the Lagrangian ADER-WENO and ADER-DG schemes proposed in \cite{Lagrange3D,ALELTS2D,ALEDG} 
to moving nonconforming unstructured meshes in a well balanced manner. Finally, we envisage to remove the mesh constraints and design a well balanced scheme for completely general moving nonconforming unstructured meshes.

\section*{Acknowledgments}

The research presented in this paper has been partially financed by the European Research Council (ERC) under the 
European Union's Seventh Framework Programme (FP7/2007-2013) with the research project \textit{STiMulUs}, 
ERC Grant agreement no. 278267.
This research has been also supported by the Spanish Government and FEDER
through the research project MTM2015-70490-C2-1-R and the Andalusian Government research projects
P11-FQM-8179 and P11-RNM-7069.
Moreover this project has received funding from the European Union's Horizon 2020 research and innovation Programme under the Marie Sklodowska-Curie grant agreement no. 642768.




\bibliographystyle{mnras}
\bibliography{references} 




\appendix

\section{Proof of well balancing for a general element in 2D}
\label{app.WellBalancing_1element}

\begin{figure}
	\includegraphics*[width=0.9\linewidth]{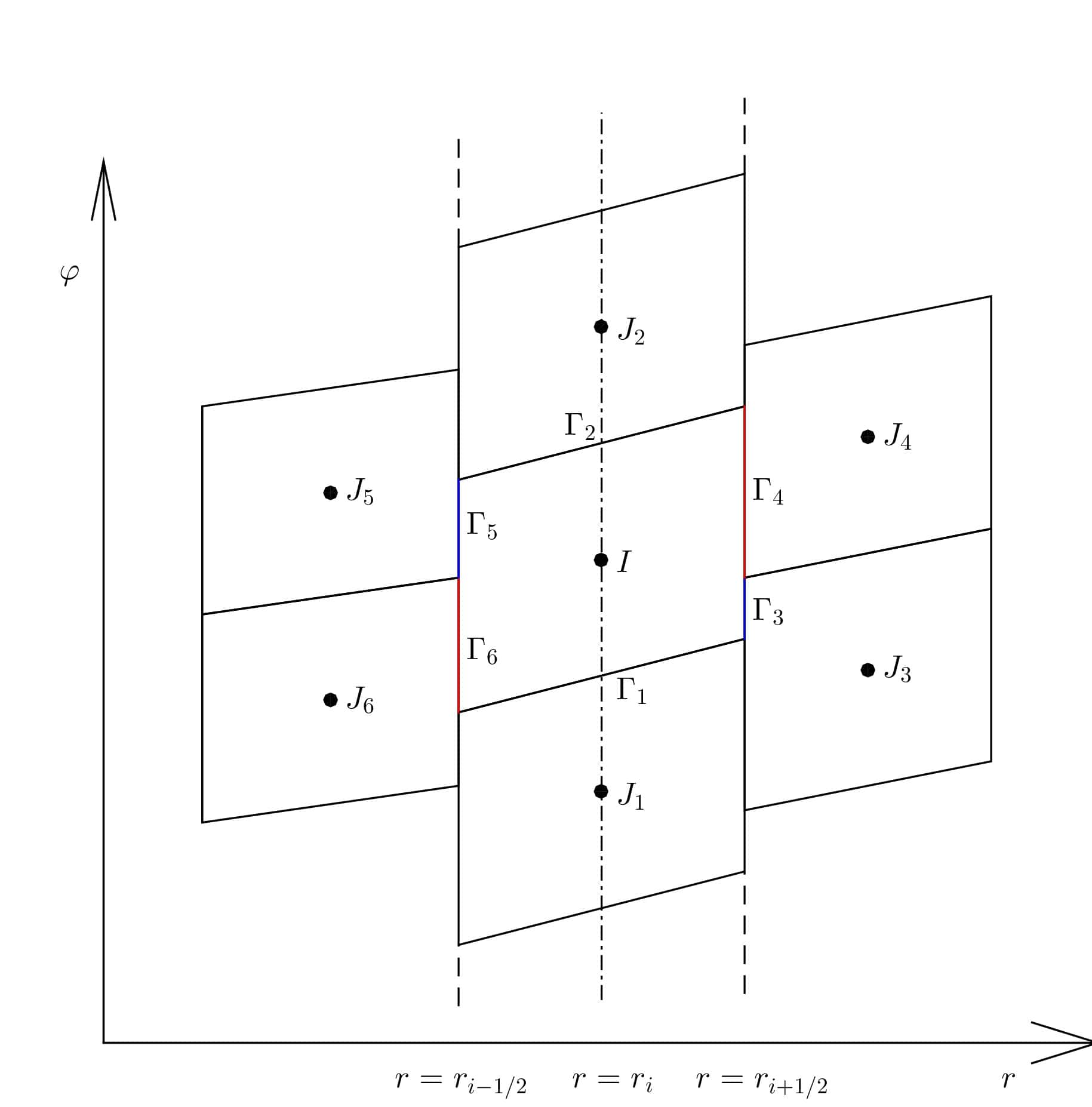}
	\caption{Portion of a general nonconforming mesh that satisfies the constraints in Section \ref{ssec.MovingDomainDiscretization}. We consider an element $I$, and its neighbors $J_i, i=1, \dots 6$, respectively through the edges $\Gamma_j, i=1,\dots,6$. In particular $\Gamma_{1,2}$ are parallel, $\Gamma_{3,4,5,6}$ lie on vertical straight lines and the barycenter of $I, J_1$ and $J_2$ have the same $r$ coordinate.}
	\label{fig.OneElementMeshForWBflux}
\end{figure}

In this section we recall the first order ALE one-step finite volume scheme in two space dimensions and we show that our formulation is well balanced for each element of a mesh  that satisfies the constraints stated at the beginning of Section \ref{ssec.MovingDomainDiscretization}.
Consider a generic element $I$ and its neighbors $J_i, i=1, \dots 6$, respectively through the edges $\Gamma_j, i=1,\dots,6$, as depicted in Figure \ref{fig.OneElementMeshForWBflux}.

As derived in Section \ref{sec.WBALE2d} our first order ALE scheme can be written as
\be
|T_I^{n+1}| \, \Q_I^{n+1} \!=  |T_I^n| \, \Q_I^n 
& \!-\! \sum \limits_{J_i} \int_0^1 \!\!\! \int_0^1 
| \partial C_{I, J_i\!}^n| \, \tilde{\mathbf{D}}_{I, J_i\!}  \cdot \mathbf{\tilde n}_{I, J_i} \,d\chi d\tau \!\!\!\!\!\!\!\!\!\!\!\! \!\!\!\!\\
\ee 
and a sufficient condition to be well balanced is that 
\be
\label{eq.SuffCond_2d}
\sum \limits_{J_i} \,\,\int_0^1 \!\! \int_0^1 
| \partial C_{I,J_i}^n| \ \tilde{\mathbf{D}}_{I,J_i}  \cdot \mathbf{\tilde n}_{I,J_i}\, d\chi d\tau = \0
\ee
when evaluated on equilibrium states.

Note that $\Gamma_{3,4,5,6}$ are parallel to the radial direction so the normal vectors are $\mbf{\tilde{n}} = (n_r, 0, 0)$, hence the flux across these edges is exactly the 1D flux, which has already been proven to be zero when evaluated on stationary solutions.

Therefore \eqref{eq.SuffCond_2d} reduces to 
\be
\int_0^1 \!\! \int_0^1 \left ( 
| \partial C_{I,J_1\!}^n| \ \tilde{\mathbf{D}}_{I,J_1\!} \! \cdot \! \mathbf{\tilde n}_{I,J_1\!} \, +\, 
| \partial C_{I,J_2\!}^n| \ \tilde{\mathbf{D}}_{I,J_2\!}  \!\cdot \! \mathbf{\tilde n}_{I,J_2\!}\,  d\chi d\tau \right) 
\ee
where, since $\Gamma_{1,2}$ are parallel and have the same length, 
\be
\mathbf{\tilde n}_{I,J_2\!} = - \mathbf{\tilde n}_{I,J_1\!} = (\tilde{n}_r, \tilde{n}_\varphi, \tilde{n}_t ) \ \text{ and } \ | \partial C_{I,J_2\!}^n| = | \partial C_{I,J_1\!}^n|,
\ee so we can rewrite
\be 
\int_0^1 \!\! \int_0^1  \!
| \partial C_{I,J_1\!}^n| \left ( \tilde{\mathbf{D}}_{I,J_1\!} \! \cdot \! \mathbf{\tilde n}_{I,J_1\!} \, -\, 
 \ \tilde{\mathbf{D}}_{I,J_2\!}  \!\cdot \! \mathbf{\tilde n}_{I,J_1\!}\,  d\chi d\tau \right).
\ee
Now, by exploiting \eqref{eq.WBspacetimeFlux_3}  the integrand can be rewritten as 
\be
& \left| \partial C_{I,J_1\!}^n\right | \Biggl ( \frac{1}{2}  \left  ( {\f(\q^{E}_{J_1}) + \f(\q^{E}_{I})} + { \mathcal{B}_{I\!,J_1\!} \left(  \q^{E}_{J_1} - \q^{E}_{I} \right)  } \right ) \tilde{n}_r \\
+ & \frac{1}{2} \!\left  ( {\g(\q^{E}_{J_1}) \!+\! \g(\q^{E}_{I})} \right ) \! \tilde{n}_\varphi 
+ \frac{1}{2} \!\left (\q^{E}_{J_1} \!+\! \q^{E}_{I}\right )\! \tilde{n}_t
- \frac{1}{2} \Vn_{I\!,J_1\!} \!\left (\q^{E}_{J_1}\!-\!\q^{E}_{I} \right ) \!\!\!\!\!\!\!\!\!\!\!  \\
- \frac{1}{2} & \left  ( {\f(\q^{E}_{J_2}) + \f(\q^{E}_{I})} + { \mathcal{B}_{I\!,J_2\!} \left(  \q^{E}_{J_1} - \q^{E}_{I} \right)  } \right ) \tilde{n}_r \\
- & \frac{1}{2} \!\left  ( {\g(\q^{E}_{J_2}) \!+\! \g(\q^{E}_{I})} \right ) \!\tilde{n}_\varphi 
- \frac{1}{2} \!\left (\q^{E}_{J_2} \!+\! \q^{E}_{I}\right ) \tilde{n}_t
+ \frac{1}{2} \!\Vn_{I\!,J_2\!} \left (\q^{E}_{J_2}\!-\!\q^{E}_{I} \right ) \Biggr ).  \!\!\!\!\!\!\!\!\!\!\!
\ee 
We already know that the component multiplied by $\tilde{n}_r$ vanishes at the equilibrium. 
Moreover, since the barycenters of $I, J_1, J_2$ are aligned along the same  straight line $r=r_i$,
\be
\label{eq.EqualValueAtEq}
\q^{E}_{J_2} = \q^{E}_{J_1},
\ee
and so the terms multiplied by $\tilde{n}_\varphi $ and $\tilde{n}_t$ cancel between them.
For what concerns the viscosity, in the case of the Osher-Romberg scheme we refer to 
\eqref{viscosity_n_osher2}-\eqref{eq.OsherRomberg2d_Rterm} that proves 
\be
\Vn_{I\!,J_i\!} \!\left (\q^{E}_{J_i}\!-\!\q^{E}_{I} \right ) = 0
\ee
provided that the rest of the scheme is well balanced (as we have just proven).
For the HLL-type flux we have 
\be
- \alpha^0_{I,J_1} \left (\tilde{I}_{\frac{I+J_1}{2}}\, n_r + I n_\varphi \right) \left (\q^{E}_{J_1}\!-\!\q^{E}_{I} \right ) - \alpha^1_{I,J_1} \mathcal{R}_{I,J_1} \\ 
+ \alpha^0_{I,J_2} \left (\tilde{I}_{\frac{I+J_2}{2}}\, n_r + I n_\varphi \right) \left (\q^{E}_{J_2}\!-\!\q^{E}_{I} \right ) + \alpha^1_{I,J_2} \mathcal{R}_{I,J_2}
\ee 
where $\tilde{I}_{\frac{I+J_i}{2}}$ vanishes as in the one dimensional case, $\mathcal{R}_{I,J_i}$ vanishes because we have already proven that the rest of the scheme vanishes, and the term multiplied by $n_\varphi$ cancels because of \eqref{eq.EqualValueAtEq}.


\bsp	
\label{lastpage}
\end{document}